\magnification=1200
\input amstex
\documentstyle{amsppt}
\hoffset=-0.5pc
\nologo
\vsize=57.2truepc
\hsize=38.5truepc
\spaceskip=.5em plus.25em minus.20em

\define\fiel{}
\define\fra{\frak}

\define\Bobb{\Bbb}

\define\Delt{\Delta}
\define\Del{D}
\define\do{d^{\nabla}}
\define\akmantwo {1}
\define\almolone {2}
\define\atiyathr {3}
\define\bangoone {4}
\define\bangotwo {5}
\define\batviltw {6}
\define\batvilfo {7}
\define\batavilk {8}
\define\ecartone {9}
\define\ecarttwo{10}
\define\drinftwo{11}
\define\fronione{12}
\define\gersthtw{13}
\define\getzlthr{14}
\define\helleron{15}
\define\poiscoho{16}
\define\souriau {17}
\define\extensta{18}
\define\bv      {19}
\define\duality {20}
\define\twilled {21}
\define\banach  {22}
\define\lradq   {23}
\define\berikas {24}
\define\kaehler {25}
\define\qr      {26}
\define\perturba{27}
\define\cohomolo{28}
\define\modpcoho{29}
\define\intecoho{30}
\define\huebkade{31}
\define\huebstas{32}
\define\bkellone{33}
\define\kikkaone{34}
\define\kinywein{35}
\define\kjesetwo{36}
\define\kosmatwe{37}
\define\kosmathr{38}
\define\kosmafou{39}
\define\kosmafte{40}
\define\koszulon{41}
\define\kravcthr{42}
\define\liuleone{43}
\define\liuletwo{44}
\define\mackone {45}
\define\maninbtw{46}
\define\molinobo{47}
\define\nijenhui{48}
\define\nomiztwo{49}
\define\palaione{50}
\define\breinhar{51}
\define\rinehone{52}
\define\roytbtwo{53}
\define\roytwein{54}
\define\sabimikh{55}
\define\sarkhone{56}
\define\sarkhtwo{57}
\define\stashnin{58}
\define\stashset{59}
\define\vanesthr{60}
\define\weinomni{61}
\define\yamaguti{62}
\topmatter
\title Higher homotopies and Maurer-Cartan algebras: \\ 
quasi-Lie-Rinehart, Gerstenhaber, and Batalin-Vilkovisky algebras
\endtitle
\author Johannes Huebschmann
\endauthor
\affil 
Universit\'e des Sciences et Technologies de Lille
\\
UFR de Math\'ematiques
\\
CNRS-UMR 8524
\\
F-59 655 VILLENEUVE D'ASCQ C\'edex, France
\\
Johannes.Huebschmann\@math.univ-lille1.fr
\endaffil
\date{February 25, 2004}
\enddate
\dedicatory
Dedicated to Alan Weinstein on his 60th birthday
\enddedicatory
\abstract{
Higher homotopy generalizations of Lie-Rinehart algebras,
Gerstenhaber, and Batalin-Vilkovisky algebras are explored.
These are defined in terms of various antisymmetric
bilinear operations satisfying weakened versions of the Jacobi identity,
as well as in terms of operations involving more than two variables
of the Lie triple systems kind.
A basic tool is the Maurer-Cartan algebra---the algebra of
alternating forms on a vector space so that Lie brackets
correspond to square zero derivations of this algebra---and
multialgebra generalizations thereof. The higher homotopies
are phrased in terms of these multialgebras.
Applications to foliations are discussed:
objects which serve as replacements for the Lie algebra
of vector fields on the \lq\lq space of leaves\rq\rq\ 
and for the algebra of multivector fields are developed,
and the spectral sequence of a foliation is shown to arise
as a special case of a more general spectral sequence including
as well the Hodge-de Rham spectral sequence.}
\endabstract
\endkeywords
\endtopmatter
\document

\leftheadtext{Johannes Huebschmann}
\rightheadtext{Quasi-Lie-Rinehart, Gerstenhaber, and BV-algebras}

\beginsection Introduction 

In this paper we will explore, in the framework of Lie-Rinehart algebras
and suitable higher homotopy generalizations thereof, various antisymmetric
bilinear operations satisfying weakened versions of the Jacobi identity,
as well as similar operations involving more than two variables; such 
operations have recently arisen in algebra, differential geometry, and 
mathematical physics but are lurking already behind a number of classical 
developments. 
Our aim is to somewhat unify these structures by means of the 
relationship between Lie-Rinehart, Gerstenhaber, and Batalin-Vilkovisky 
algebras which we first observed in our paper \cite\bv. This will be, perhaps, 
a first step towards {\it taming the bracket zoo that arose recently in 
topological field theory\/}, cf. what we wrote in the introduction to \cite\bv.
The notion of Lie-Rinehart algebra and its generalization are likely to
provide a good conceptual framework for that purpose.
It will also relate new notions like those of Gerstenhaber and 
Batalin-Vilkovisky algebra, and generalizations thereof, with classical ones 
like those of connection, curvature, and torsion,  
as well as with less classical ones like Yamaguti's triple product
\cite\yamaguti\ and operations of the kind introduced in \cite\kinywein; 
it will  connect new developments
with old results due to E. Cartan \cite\ecartone\  and Nomizu \cite\nomiztwo\ 
describing the geometry of Lie groups and of reductive homogeneous spaces
and, more generally, with more recent results in the geometry of Lie loops
\cite{\kikkaone,\,\sabimikh}. We will see that the new structures have 
incarnations in mathematical nature, e.~g. in the theory of foliations.
The higher homotopies which are exploited below are of a special kind, though,
where only the first of an (in general) infinite family is non-zero.
\smallskip
Let $R$ be a commutative ring with 1. A {\it Lie-Rinehart algebra\/}
$(A,L)$ consists of a commutative $R$-algebra $A$, an $R$-Lie algebra $L$,
an $A$-module structure on $L$, and an action $L \otimes_R A \to A$ of $L$ 
on $A$ by derivations. These are required to satisfy suitable 
compatibility conditions which arise by abstraction from the pair
$(A,L) = (C^{\infty}(M),\roman{Vect}(M))$ consisting of the smooth functions
$C^{\infty}(M)$ and smooth vector fields $\roman{Vect}(M)$ on a smooth  
manifold $M$. In a series of papers \cite{\poiscoho--\twilled},
we studied these objects and variants thereof and used them to solve various 
problems in algebra and geometry. See \cite\lradq\ for a survey and leisurely 
introduction. In differential geometry, a special case of a Lie-Rinehart 
algebra arises from the space of sections of a Lie algebroid.
\smallskip
In \cite{\bv,\,\twilled,\,\banach} we have shown that certain Gerstenhaber 
and Batalin-Vilkovisky algebras admit natural interpretations in terms of 
Lie-Rinehart algebras. The starting point was the following observation:
It is nowadays well understood that a skew-symmetric bracket on a vector
space $\fra g$ is a Lie-bracket (i.~e. satisfies the Jacobi identity) if and
only if the coderivation $\partial$ on the graded exterior coalgebra
$\Lambda'[s\fra g]$ corresponding to the bracket on $\fra g$ has square
zero, i.~e. is a differential; this coderivation is then the ordinary
Lie algebra homology operator. This kind of characterization is not available
for a general Lie-Rinehart algebra:
Given a commutative algebra $A$ and an $A$-module $L$,
a Lie-Rinehart structure on $(A,L)$ cannot be characterized in terms of a
coderivation on $\Lambda_A[sL]$ with reference to a suitable coalgebra
structure on $\Lambda_A[sL]$
(unless the $L$-action on $A$ is trivial);
in fact, in the Lie-Rinehart context, a certain {\sl dichotomy
between $A$-modules and chain complexes which are merely 
defined over $R$ persists thoughout\/};
cf. e.~g. the Remark 2.5.2 below.
On the other hand,
Lie-Rinehart algebra structures 
on $(A,L)$ correspond to Gerstenhaber algebra structures
on the exterior $A$-algebra $\Lambda_A[sL]$; cf. e.~g.
\cite\kosmathr.
In particular, when $A$ is the ground ring
and $L$  just an ordinary Lie algebra $\fra g$,
under the obvious identification of
$\Lambda[s\fra g]$ and
$\Lambda'[s\fra g]$ as graded $R$-modules,
the (uniquely determined) generator of the Gerstenhaber bracket
on $\Lambda[s\fra g]$ is exactly the Lie algebra homology operator
on $\Lambda'[s\fra g]$.
Given a general commutative algebra $A$ and an $A$-module $L$,
the interpretation of Lie-Rinehart algebra structures 
on $(A,L)$ in terms of Gerstenhaber algebra structures
on $\Lambda_A[sL]$
provides, among other things, a link between 
Gerstenhaber's and Rinehart's papers \cite\gersthtw\  and \cite\rinehone\ 
(which seems to have been completely missed in the literature). In the present 
paper, we will extend this link to suitable higher homotopy notions which we 
refer to by the attribute \lq\lq quasi\rq\rq; we will introduce Lie-Rinehart 
triples, quasi-Lie-Rinehart algebras, and certain
quasi-Gerstenhaber algebras and quasi-Batalin-Vilkovisky algebras, 
and we will explore the various
relationships between these notions.
Below we will  comment on the relationship with notions of 
quasi-Gerstenhaber and quasi-Batalin-Vilkovisky algebras
already in the literature.
\smallskip
When an algebraic structure (e.~g. a commutative algebra, Lie algebra, etc.) is 
\lq\lq resolved\rq\rq\ by an object, which we here somewhat vaguely refer to
as a \lq\lq resolution\rq\rq\ (free, 
or projective, or variants thereof) 
having the given structure as its zero-th homology, on the resolution, the
algebraic structure is in general defined only up to higher homotopies; 
likewise, an $A_{\infty}$ structure is defined in terms of a bar construction 
or variants thereof, cf. e.~g. \cite\huebkade,\ \cite\huebstas\ and the 
references there. Exploiting higher homotopies of this kind, in a series of 
articles \cite{\perturba\--\intecoho} we constructed small free resolutions 
for certain classes of groups from which we then were able to do explicit 
calculations in group cohomology which until today still cannot be done by 
other methods. A historical overview related with $A_{\infty}$-structures
may be found in the Addendum to \cite\bkellone; cf. also
\cite\berikas\ and \cite\huebkade\ for more historical comments.
\smallskip
In the present paper, we will explore 
a certain higher homotopy 
related with Lie-Rinehart algebras and variants thereof. 
A Lie algebra up to higher homotopies
(equivalently: $L_{\infty}$-algebra)
on an $R$-chain complex $\fra h$
may be defined in terms of a coalgebra perturbation
of the differential
on the graded symmetric coalgebra on the suspension of $\fra h$;
alternatively, it may be defined in terms of a suitable
Maurer-Cartan algebra (see below). Since a genuine Lie-Rinehart 
structure on $(A,L)$ cannot be characterized in terms of a
coderivation on $\Lambda_A[sL]$,
the first alternative breaks down for
a general Lie-Rinehart algebra.
The higher homotopies we will explore in the present paper
do {\it not\/} live on an object close to a resolution 
of the above kind or close to a symmetric coalgebra; they 
may conveniently be phrased in terms of an {\it object of a rather
different nature\/} which,
extending terminology introduced by van Est \cite\vanesthr, we refer to 
as a {\it Maurer-Cartan\/} algebra.
A special case thereof arises in the following fashion: Given a 
finite dimensional vector space
$\fra g$ over a field $\fiel k$,
skew symmetric brackets
on $\fra g$ correspond bijectively to degree $-1$
derivations of the graded algebra 
of alternating forms on $\fra g$
(with reference to multiplication of forms),
and those brackets which satisfy the Jacobi identity
correspond to square zero derivations, i.~e. differentials.
This observation generalizes to 
Lie-Rinehart algebras of the kind $(A,L)$
under the assumption that $L$ be a finitely generated projective $A$-module;
see Theorem 2.2.16 below. For an ordinary Lie algebra
$\fra g$ over a field $\fiel k$,
in \cite\vanesthr, the resulting differential graded algebra
$\roman{Alt}(\fra g, \fiel k)$ 
(which calculates the cohomology of $\fra g$)
has been called {\it Maurer-Cartan algebra\/}.
The main point of this paper is
that higher homotopy variants of the notion of 
Maurer-Cartan algebra provide the correct framework to phrase certain higher 
homotopy versions of Lie-Rinehart-, Gerstenhaber, and Batalin-Vilkovisky 
algebras to which we will refer as {\it quasi-Lie-Rinehart-, 
quasi-Gerstenhaber, and quasi-Batalin-Vilkovisky  algebras\/}.

The differential graded algebra of alternating forms on a Lie algebra occurs, 
at least implicitly, in \cite\ecarttwo\ 
and has a long history of use since then, cf. \cite\koszulon, and once I learnt 
in a talk by van Est that this algebra has been used by E.~Cartan in the 
1930's to characterize the structure of Lie groups and Lie algebras. 
\smallskip
For the reader's convenience, we will explain briefly and somewhat informally
a special case of a quasi-Lie-Rinehart algebra at the present stage: Let 
$(M,\Cal F)$ be a 
foliated manifold, the foliation being written as 
$\Cal F$, let $\tau_{\Cal F}$ be the tangent bundle of the foliation $\Cal F$,
and choose a complement $\zeta$ of $\tau_{\Cal F}$ so that the tangent bundle 
$\tau_M$ of $M$ may be written as $\tau_M = \tau_{\Cal F} \oplus \zeta$. Let 
$L_{\Cal F}\subseteq \roman{Vect}(M)$ be the Lie algebra of smooth vector 
fields tangent to the foliation $\Cal F$, and let $Q$ be the 
$C^{\infty}(M)$-module $\Gamma(\zeta)$ of smooth sections of $\zeta$. The Lie 
bracket in $\roman{Vect}(M)$ induces a left $L_{\Cal F}$-module structure on 
$Q$---the Bott connection@--- and the space $Q^{L_{\Cal F}}$ of invariants, 
that is, of vector fields on $M$ which are horizontal (with respect to the 
decomposition $\tau_M = \tau_{\Cal F} \oplus \zeta$) and constant on the
leaves inherits a Lie bracket. 
The standard complex $\Cal A$ 
arising from a fine resolution of the sheaf of germs of functions on $M$ which 
are constant on the leaves acquires a differential graded algebra structure and 
has $\roman H^0(\Cal A)$ equal to the algebra of functions on $M$ which are 
constant on the leaves, and the Lie algebra $Q^{L_{\Cal F}}$ of invariants
arises as $\roman H^0(\Cal Q)$ where $\Cal Q$ is the complex coming 
from a fine resolution of the sheaf $\Cal V_Q$ of germs of vector fields on 
$M$ which are horizontal (with respect to the decomposition 
$\Gamma(\tau_M) = L_{\Cal F} \oplus Q$) and constant on the leaves.
In a sense, $Q^{L_{\Cal F}}$ is the  {\sl Lie
algebra of vector fields on the \lq\lq space of leaves\rq\rq\/}, 
that is, the space of 
sections of a certain {\sl geometric object which may be seen 
as a replacement for 
the in general non-existant tangent bundle of the \lq\lq space of 
leaves\rq\rq\/}.
Within our approach, this philosophy is pushed further in the following 
fashion: The pair 
$(\Cal A,\Cal Q)$ acquires what we will call a {\it quasi-Lie-Rinehart
structure\/} in an obvious fashion; see (4.12) and (4.15) below for the 
details. We view $\Cal A$ as the {\it algebra of generalized 
functions\/} and  $\Cal Q$ as the {\it generalized Lie algebra of vector 
fields\/} for the foliation. The pair $(\roman H^0(\Cal A),\roman H^0(\Cal Q))$ 
is necessarily a Lie-Rinehart algebra, and the entire cohomology 
$(\roman H^*(\Cal A),\roman H^*(\Cal Q))$ acquires a graded Lie-Rinehart 
algebra structure. As a side remark, we note that here the resolution of 
the sheaf $\Cal V_Q$ is by no means a projective one; indeed, it is a fine 
resolution of that sheaf, the bracket on $\Cal Q$ is 
not an ordinary Lie(-Rinehart) bracket, in particular, does not satisfy the 
Jacobi identity, and the entire additional structure is encapsulated in certain 
homotopies which may conveniently be phrased in terms of a suitable 
Maurer-Cartan algebra which here arises from the de Rham algebra of $M$.
When the foliation does {\it not\/} come from a fiber bundle,
the structure of the graded Lie-Rinehart algebra 
$(\roman H^*(\Cal A),\roman H^*(\Cal Q))$
will in general be more complicated than that for the case when the foliation
comes from a fiber bundle. Thus the cohomology of a quasi-Lie-Rinehart algebra
involves an ordinary Lie-Rinehart algebra in degree zero
but in general contains considerably more information.
In particular, in the case of a foliation it contains more than just
\lq\lq functions and vector fields on the space of leaves\rq\rq; the additional
information partly includes the {\it history\/} of the 
\lq\lq space of leaves\rq\rq, 
that is, it includes information as to how this space
arises from the foliation, how the leaves sit inside the ambient space,
about singularities, etc.
In Section 6 we will show that,
when the foliation is transversely orientable with a 
{\it basic\/} transverse volume form $\omega$, a
corresponding quasi-Batalin-Vilkovisky algebra isolated in Theorem 6.10 below
has an underlying quasi-Gerstenhaber algebra
which, in turn, yields a kind of generalized Schouten algebra  
(generalized algebra of multivector fields) for the foliation;
the cohomology  of this quasi-Gerstenhaber algebra may then be viewed as the 
Schouten algebra for the \lq\lq space of leaves\rq\rq. 
See (6.15) below for details.
\smallskip
Thus our approach will provide new insight, for example, into the geometry of 
foliations; see in particular (1.12), (2.10), (4.15), (6.15)
below. The formal structure behind foliations which we will phrase in terms of
quasi-Lie-Rinehart algebras and its offspring does not seem to have been 
noticed in the literature before---indeed, it involves, among a number of other
things, a suitable grading which seems unfamiliar in the literature
on quasi-Gerstenhaber and quasi-Batalin-Vilkovisky algebras, cf.
(6.17) below---, nor the formal connections with Yamaguti's 
triple product and with Lie loops.
\smallskip
A simplified version of the question we will examine is this: Given a Lie 
algebra $\fra g$ with a decomposition $\fra g =\fra h \oplus \fra q$ where 
$\fra h$ is a Lie subalgebra, what kind of structure does $\fra q$ then inherit?
Variants of this question and possible answers may be found at a number of 
places in the literature, cf. e.~g. \cite{\ecartone,\,\nomiztwo} where,
in particular, in a global situation, an answer is given
for reductive homogeneous spaces. In the 
framework of Lie-Rinehart algebras, 
this issue does not seem to have been raised 
yet, not even for the special case of Lie algebroids.
\smallskip
As a byproduct, we find a certain formal relationship between Yamaguti's triple 
product and certain forms $\Phi^*_*$ which may be found in \cite\koszulon.
In particular, the failure of a quasi-Gerstenhaber bracket to satisfy 
the Jacobi identity 
is measured by an additional piece of structure which we refer to as 
an $h$-{\it Jacobiator\/}; an $h$-Jacobiator, in turn, is defined in terms of 
Koszul's forms $\Phi^3_*$.
Likewise the quadruple and quintuple products studied in Section 3 below are 
related with Koszul's forms, and these, in turn, are related with certain 
higher order operations which may be found e.~g. in \cite\sabimikh.
We do not pursue this here; we hope to eventually come back 
to it in another article.
\smallskip
A Courant algebroid has been shown in \cite\roytwein\ 
to acquire an $L_{\infty}$-structure, that is, a Lie algebra structure
up to higher homotopies.
The present paper paves, perhaps, the way towards
finding a higher homotopy Lie-Rinehart or higher homotopy
Lie algebroid structure on a Courant algebroid incorporating
the Courant algebroid structure.
\smallskip
Graded quasi-Batalin-Vilkovisky algebras
have been explored already in \cite\getzlthr.
Our notions of quasi-Gerstenhaber and quasi-Batalin-Vilkovisky algebra,
while closely related,
do {\it not\/} coincide with those in 
\cite\bangoone,\,
\cite\bangotwo,\,
\cite\getzlthr,\, 
\cite\kosmafte,\,
\cite \roytbtwo. In particular,
our algebras are bigraded while those in the quoted references
are ordinary graded algebras; the appropriate totalization 
(forced, as noted above,
by our application of the newly developed algebraic structure
to foliations and written in Section 6 below as the functor $\roman{Tot}$)
of our bigraded objects leads to differential graded objects 
which are {\it not\/} equivalent to those in the quoted references.
See Remark 6.17 below for more details on the relationship
between the various notions.
Also the approaches differ in motivation;
the guiding idea behind \cite\getzlthr\ and \cite\kosmafte\  seems to be 
Drinfeld's quasi-Hopf algebras.
Our motivation, as indicated above, comes from foliations
and the search for appropriate algebraic notions
encapsulating the infinitesimal structure 
of the \lq\lq space of leaves\rq\rq\ and its history,
as well as the search for a corresponding Lie-Rinehart generalization
of the operations on a reductive homogenous space
isolated by Nomizu  and elaborated upon by Yamaguti
(mentioned earlier) and taken up again by
M. Kinyon and A. Weinstein in \cite\kinywein. 
Indeed, the present 
paper was prompted by the preprint versions of \cite\kinywein\ and 
\cite\weinomni. It is a pleasure to dedicate it to Alan Weinstein.
Throughout this work I have been stimulated
by M. Kinyon via some e-mail correspondence at an early stage
of the project as well as by 
M. Bangoura, P. Michor, D. Roytenberg and Y. Kosmann-Schwarzbach.
I am indebted to J. Stasheff 
and to the referees
for a number of comments on a draft of the 
manuscript which helped improve the exposition.
\smallskip
This work was partly carried out and presented 
during two stays at the Erwin Schr\"odinger Institute at Vienna.
I wish to express my gratitude for hospitality and support.

\beginsection 1. Lie-Rinehart triples

Let $R$ be a commutative ring with 1,
not necessarily a field; $R$ could be, for example,
the algebra of smooth functions on a smooth manifold, cf. \cite\duality. 
The problem we wish to explore is this:
\newline\noindent
{\smc Question 1.1.}
Given  a Lie-Rinehart algebra $(A,L)$ and an $A$-module direct sum decomposition
$L= H \oplus Q$ inducing an $(R,A)$-Lie algebra structure on $H$, what kind of 
structure does then $Q$ inherit, and by what additional structure are $H$ and 
$Q$ related?
\newline\noindent
{\smc Question 1.2.}
Given an $(R,A)$-Lie algebra structure on $H$ and the (new) structure 
(which we will isolate below) on $Q$, what kind of additional structure
turns the $A$-module direct sum $L= H \oplus Q$ into an $(R,A)$-Lie algebra
in such a way that the latter induces the given structure on $H$ and $Q$?
\smallskip\noindent
{\smc Example 1.3.1.}
Let $\fra g$ be an ordinary $R$-Lie algebra with a decomposition 
$\fra g =\fra h \oplus \fra q$ where $\fra h$ is a Lie subalgebra.
Recall that the  decomposition of $\fra g$ is said to be {\it reductive\/}
\cite\nomiztwo\ provided $[\fra h,\fra q] \subseteq \fra q$. Such a reductive 
decomposition arises from a reductive homogeneous space 
\cite{\ecartone,\,\kikkaone,\,\nomiztwo,\,\sabimikh,\,\yamaguti}.
For example, every homogeneous space of a compact Lie group or,
more generally, of a reductive Lie group, is reductive.
Nomizu has shown that, on such a reductive homogeneous space, the torsion and 
curvature of the \lq\lq canonical affine connection of the second kind\rq\rq\ 
(affine connection with parallel torsion and curvature) yield a bilinear and 
a ternary operation which, at the identity, come down to a certain bilinear and 
ternary operation on the constituent $\fra q$ \cite\nomiztwo, and Yamaguti 
gave an algebraic characterization of pairs of such operations \cite\yamaguti.
\smallskip\noindent
{\smc Example 1.3.2.} A quasi-Lie bialgebra $(\fra h, \fra q)$,  cf. 
\cite\kosmatwe, consists of a
(real or complex) Lie algebra $\fra h$ and a
(real or complex) vector space $\fra q$
with suitable 
additional structure where $\fra q = \fra h^*$, so that
$\fra g =\fra h \oplus \fra h^*$ is an ordinary Lie algebra; the pair 
$(\fra g,\fra h)$ is occasionally referred to in the literature as a 
{\it Manin pair\/}. Quasi-Lie bialgebras arise as classical limits of 
quasi-Hopf algebras; these, in turn, were introduced by Drinfeld \cite\drinftwo.
\smallskip\noindent
{\smc Example 1.4.1.}
Let $R$ be the field $\Bobb R$ of real numbers, let $(M,\Cal F)$ be a 
foliated manifold, let $\tau_{\Cal F}$ be the tangent bundle of the foliation
$\Cal F$, and choose a complement $\zeta$ of $\tau_{\Cal F}$ so that the 
tangent bundle $\tau_M$ of $M$ may be written as 
$\tau_M = \tau_{\Cal F} \oplus \zeta$. Thus, as a vector bundle, $\zeta$ is 
canonically isomorphic to the normal bundle of the foliation. Let 
$(A,L)$ be the Lie-Rinehart algebra $(C^{\infty}(M),\roman{Vect}(M))$, let 
$L_{\Cal F}\subseteq L$ be the $(\Bobb R, A)$-Lie algebra of smooth vector 
fields tangent to the foliation $\Cal F$, and let $Q$ be the $A$-module  
$\Gamma(\zeta)$ of smooth sections of $\zeta$. Then $L= L_{\Cal F} \oplus Q$
is a $A$-module direct sum decomposition of the $(\Bobb R, A)$-Lie algebra 
$L$, and the question arises what kind of Lie structure $Q$ carries.
This question, in turn, may be subsumed under the more general
question to what extent the \lq\lq space of leaves\rq\rq\ 
can be viewed as a smooth
manifold. This more general question is not only of academic interest since, 
for example, in interesting physical situations, the true classical state space 
of a constrained system is the \lq\lq space of leaves\rq\rq\ of a foliation 
which is in
general not fibrating, and the Noether theorems are conveniently phrased 
in the framework of foliations.

\smallskip\noindent
{\smc Example 1.4.2.}
Let $R$ be the field $\Bobb C$ of complex numbers, $M$ a smooth complex 
manifold $M$, $A$ the algebra of smooth complex functions on $M$, $L$ the 
$(\Bobb C, A)$-Lie algebra of smooth complexified vector fields, and let $L'$ 
and $L''$ be the spaces of smooth sections of the holomorphic and 
antiholomorphic tangent bundle of $M$, respectively. Then $L'$ and $L''$ are
$(\Bobb C,A)$-Lie algebras, and $(A,L',L'')$ is a {\it twilled Lie-Rinehart 
algebra\/} in the sense of \cite{\twilled,\, \banach}. Adjusting the notation 
to that in (1.4.1), let $H=L'$ and $Q=L''$. Thus, in this 
particular case, $Q=L''$ is in fact an ordinary $(R,A)$-Lie algebra, and the 
additional structure relating $H$ and $Q$ is encapsulated in the notion of 
{\it twilled Lie-Rinehart algebra\/}. 
The integrability condition for an almost complex structure may be
phrased in term of the twilled Lie-Rinehart axioms; see 
\cite{\twilled,\,\banach} for details.
\smallskip
The situation of Example 1.4.1 is somewhat more general than that of Example 
1.4.2 since in Example 1.4.1 the constituent $Q$ carries a structure which is 
more general than that of an ordinary $(R,A)$-Lie algebra. Another example 
for a decomposition of the kind spelled out in Questions 1.1 and 1.2 above
arises from combining the situations of Example 1.4.1 and of Example 1.4.2, 
that is, from a smooth manifold foliated by holomorphic manifolds, and yet 
another example arises
from a holomorphic foliation. Abstracting from these examples, 
we isolate the notion of {\it Lie-Rinehart triple\/}. For ease of exposition,
we also introduce the weaker concepts of 
{\it almost pre-Lie-Rinehart triple\/} and {\it pre-Lie-Rinehart triple\/}.
Distinguishing between these three notions may appear 
pedantic but will clarify the statement of Theorem 2.7 below.
See also Remark 2.8.4 below.
As for the terminology we note that our notion of {\it triple\/}
is {\it not\/} consistent with the usage of {\it Manin triple\/}
in the literature. 
However, a Lie-Rinehart algebra involves a pair consisting of an algebra
and a Lie algebra, and in this context, 
it is also common in the literature to refer to this structure as
a {\it pair\/} which, in turn, 
is not consistent with the notion of {\it Manin pair\/}.
We therefore prefer to use our terminology
Lie-Rinehart triple etc.
\smallskip
Let $A$ be a commutative $R$-algebra. Consider two $A$-modules $H$ and $Q$,
together with 
\newline\noindent
--- skew-symmetric $R$-bilinear brackets of the kind $(1.5.1.H)$ and 
$(1.5.1.Q)$ below, not necessarily Lie brackets;
\newline\noindent
--- $R$-bilinear operations of the kind (1.5.2.$H$), (1.5.2.$Q$), (1.5.3), 
(1.5.4)  below; and 
\newline\noindent
--- a skew-symmetric $A$-bilinear pairing $\delta$ of the kind (1.5.5) below:
$$
\alignat1
[\,\cdot\,,\cdot\,]_H
&\colon H \otimes_R H
@>>>
H,
\tag1.5.1.$H$
\\
[\,\cdot\,,\cdot\,]_Q
&\colon Q \otimes_R Q
@>>>
Q,
\tag1.5.1.$Q$
\\
H \otimes_R A &@>>> A,
\quad x\otimes_R a \mapsto x(a),
\quad x\in H,\  a \in A,
\tag1.5.2.$H$
\\
Q \otimes_R A &@>>> A,
\quad \xi\otimes_R a \mapsto \xi(a),
\quad \xi \in Q,\  a \in A,
\tag1.5.2.$Q$
\\
\cdot\,\colon &H \otimes_R Q @>>> Q,
\tag1.5.3
\\
\cdot\,\colon &Q \otimes_R H @>>> H,
\tag1.5.4
\\
\delta\colon &Q \otimes_A Q @>>> H .
\tag1.5.5
\endalignat
$$
We will say that the data $(A,H,Q)$ constitute an {\it almost pre-Lie-Rinehart 
triple\/} provided they satisfy (i), (ii), and (iii) below.
\smallskip\noindent
(i) The values of the adjoints $H @>>> \roman{End}_R(A)$ and
$Q @>>> \roman{End}_R(A)$ of (1.5.2.$H$) and (1.5.2.$Q$) respectively
lie in $\roman{Der}_R(A)$;
\smallskip\noindent
(ii) (1.5.1.$H$), (1.5.2.$H$) and the $A$-module structure on $H$
and, likewise, (1.5.1.$Q$), (1.5.2.$Q$) and the $A$-module structure on $Q$,
satisfy the following Lie-Rinehart axioms (1.5.6.$H$), (1.5.7.$H$), 
(1.5.6.$Q$), (1.5.7.$Q$):
$$
\alignat1
(a x) (b) &= a (x (b)),
\quad a, b \in A,\  x\in H,
\tag1.5.6.$H$
\\
[x, a y]_H &= x (a) y+ a [x, y]_H,
\quad a \in A,\ x,y \in H,
\tag1.5.7.$H$
\\
(a \xi) (b) &= a (\xi (b)),
\quad a, b \in A,\  \xi\in Q,
\tag1.5.6.$Q$
\\
[\xi, a \eta]_Q &= \xi (a) \eta+ a [\xi, \eta]_Q,
\quad a \in A,\ \xi,\eta \in Q;
\tag1.5.7.$Q$
\endalignat
$$
\noindent
(iii) (1.5.3) and (1.5.4)
behave  like connections, that is,
for $a \in A,\ x \in H, \ \xi \in Q$, the identities
$$
\alignat1
x \cdot (a\xi) &= (x(a)) \xi + a (x\cdot \xi),
\tag1.5.8
\\
(ax) \cdot \xi &=  a(x\cdot \xi),
\tag1.5.9
\\
\xi \cdot (a x) &= (\xi(a)) x + a (\xi\cdot x),
\tag1.5.10
\\
(a\xi) \cdot x &=  a(\xi\cdot x),
\tag1.5.11
\endalignat
$$
are required to hold.
\smallskip
We will say that an almost pre-Lie-Rinehart triple $(A,H,Q)$ is a
{\it pre-Lie-Rinehart triple\/} provided that (i) $(A,H)$, endowed with 
the operations {\rm (1.5.1.$H$)} and  {\rm (1.5.2.$H$)}, is a Lie-Rinehart 
algebra---equivalently, the bracket {\rm (1.5.1.$H$)} satisfies the Jacobi 
identity@---, 
and that (ii) the operation {\rm (1.5.3)} turns $Q$ into a left 
$(A,H)$-module, that is, the \lq\lq connection\rq\rq\ given by this operation 
is \lq\lq flat\rq\rq, i.~e. satisfies the identity
$$
[x,y]_H \cdot \xi = x \cdot (y\cdot \xi) - y \cdot (x\cdot \xi),
\quad
x,y \in H,\ \xi \in Q.
\tag1.5.12
$$
{\smc (1.5.13)\/}
{\sl Thus a pre-Lie-Rinehart triple $(A,H,Q)$ consists of a Lie-Rinehart 
algebra $(A,H)$\/} (the structure of which is given by (1.5.1.$H$), 
(1.5.2.$H$)) {\sl and a left $(A,H)$-module\/} $Q$ (given by the operation 
(1.5.3) which, in turn, is required to satisfy the axioms (1.5.8) and (1.5.9)) 
{\sl together with the additional structure\/} (1.5.1.$Q$), (1.5.2.$Q$), 
(1.5.4), (1.5.5) {\sl subject to the axioms\/} (1.5.6.$Q$), (1.5.7.$Q$), 
(1.5.10), (1.5.11).
\smallskip
Given an almost pre-Lie-Rinehart triple $(A,H,Q)$, let $L=H\oplus Q$ be the  
$A$-module direct sum, and define an $R$-bilinear skew-symmetric bracket
$$
[\,\cdot\,,\cdot\,]
\colon L \otimes_R L
@>>>
L
\tag1.6.1
$$
by means of the formula
$$
[(x,\xi),(y,\eta)]
=
[x,y]_H +[\xi,\eta]_Q
+\delta(\xi,\eta)
+ x \cdot \eta
-\eta \cdot x
+ \xi \cdot y
-y \cdot \xi
\tag1.6.2
$$
and, furthermore, an operation
$$
L \otimes_R A @>>> A
\tag1.6.3
$$
in the obvious way, that is, by means of the association
$$
(\xi,x) \otimes_Ra \mapsto \xi(a) +x(a),
\quad x \in H,\, \xi \in Q,\, a \in A.
\tag1.6.4
$$
By construction, the values of the adjoint of (1.6.3)  then lie in
$\roman{Der}_R(A)$, that is, this adjoint is then of the form
$$
L = H \oplus Q @>>> \roman{Der}_R(A) .
\tag1.6.5
$$
An almost pre-Lie-Rinehart triple $(A,H,Q)$ will be said to be a
{\it Lie-Rinehart triple\/} if {\rm (1.6.1)} and  {\rm (1.6.3)} turn 
$(A,L)$ where $L=H\oplus Q$ into a Lie-Rinehart algebra. A Lie-Rinehart triple 
$(A,H,Q)$ where $\delta$ is zero is a {\it twilled Lie-Rinehart algebra\/}
\cite{\twilled,\,\banach}.
Thus  Lie-Rinehart triples generalize twilled Lie-Rinehart algebras.
\smallskip
A direct sum decomposition $L= H \oplus Q$ of an $(R,A)$-Lie algebra $L$
such that $(A,H)$ inherits a Lie-Rinehart structure yields a Lie-Rinehart triple
$(A,H,Q)$ in an obvious fashion: The brackets (1.5.1.$H$) and (1.5.1.$Q$)
result from restriction and projection; the operations (1.5.2.$H$) and 
(1.5.2.$Q$)
are  obtained by restriction as well; further, the requisite operations (1.5.3) 
and (1.5.4) are given by the composites 
$$ 
\cdot\,\colon H\otimes_R Q 
@>{[\,\cdot\,,\cdot\,]|_{H\otimes_R Q}}>> H\oplus Q @>>\roman{pr_{Q}}> Q 
\tag1.7.1 
$$ 
and 
$$ 
\cdot\,\colon Q\otimes_R H @>{[\,\cdot\,,\cdot\,]|_{Q\otimes_R H}}>> H\oplus Q 
@>>\roman{pr_{H}}> H 
\tag1.7.2 
$$ 
where, for $M= H\otimes_R Q$ and $M=Q\otimes_R H$, $[\,\cdot\,,\cdot\,]|_{M}$ 
denotes the restriction of the Lie bracket to $M$. The pairing (1.5.5) is the 
composite
$$
\delta
\colon Q\otimes_A Q @>{[\,\cdot\,,\cdot\,]|_{Q\otimes_R Q}}>> L=H\oplus Q 
@>>\roman{pr_H}> H ;
\tag1.7.3 
$$ 
at first it is only $R$-bilinear but is readily seen to be $A$-bilinear. The 
formula (1.6.2) is then merely a decomposition of the initially given bracket 
on $L$ into components according to the direct sum decomposition of $L$ into 
$H$ and $Q$, and (1.6.3) is accordingly a decomposition of the $L$-action on 
$A$. Furthermore, given $x,y \in H$ and $\xi \in Q$, in $L$ we have the identity
$$
\align
[x,y] \cdot\xi - \xi \cdot [x,y] &=
[[x,y],\xi] = [x,[y,\xi]] - [y,[x,\xi]]  
\\
&= 
x\cdot(y \cdot \xi)
-
(y\cdot \xi)\cdot x 
-
[x,\xi\cdot y] 
\\
&\quad
-y\cdot (x \cdot \xi)
+(x\cdot \xi)\cdot y -
[\xi\cdot x,y]    
\endalign
$$
which at once implies (1.5.12).

\smallskip
\noindent
{\smc Remark 1.8.1.}
Thus we see that, in particular, if an almost pre-Lie-Rinehart triple $(A,H,Q)$
is a Lie-Rinehart triple, it is necessarily a pre-Lie-Rinehart triple, 
cf. (1.5.13).
\smallskip
\noindent
{\smc Remark 1.8.2.} In the situation of Example 1.3.2, when $\fra g$ arises 
from a quasi-Lie bialgebra (so that $\fra q = \fra h^*$), in the literature, 
the piece of structure $\delta$ is often written as an element of 
$\Lambda^3 \fra h$.

\proclaim{Theorem 1.9}
A pre-Lie-Rinehart triple $(A,H,Q)$ is a genuine Lie-Rinehart triple, that is,
the bracket $[\cdot,\cdot]$, cf. {\rm (1.6.1)}, and the operation
 {\rm (1.6.3)} turn $(A,L)$ where $L=H\oplus Q$ into a 
Lie-Rinehart algebra, if and only if the brackets $[\,\cdot\,,\cdot\,]_H$ and
$[\,\cdot\,,\cdot\,]_Q$ on $H$ and $Q$, respectively, and the operations 
{\rm (1.5.3)}, {\rm (1.5.4)}, and {\rm (1.5.5)}, are related by
$$
\alignat1
\xi (x(a))- x(\xi(a)) &= (\xi \cdot x)(a) -(x \cdot \xi)(a) 
\tag1.9.1
\\
x\cdot[\xi,\eta]_Q
&=
[x\cdot\xi,\eta]_Q
+
[\xi,x\cdot \eta]_Q
-(\xi \cdot x)\cdot \eta
+(\eta \cdot x)\cdot \xi
\tag1.9.2
\\
\xi\cdot[x,y]_H
&=
[\xi\cdot x,y]_H
+
[x,\xi\cdot y]_H
-(x \cdot \xi)\cdot y
+(y \cdot \xi)\cdot x ,
\tag1.9.3
\\
\xi (\eta(a))- \eta(\xi(a)) &= [\xi,\eta]_Q(a) + (\delta(\xi,\eta))(a) 
\tag1.9.4
\\
[\xi,\eta]_Q \cdot x 
=
\xi\cdot(\eta \cdot x) &-\eta\cdot (\xi \cdot x)
- \delta(x\cdot\xi,\eta)
-\delta(\xi,x\cdot \eta)
+[x,\delta(\xi,\eta)]_H 
\tag1.9.5
\\
\sum_{(\xi,\eta,\vartheta)\ \roman{cyclic}}\left(
[[\xi,\eta]_Q,\vartheta]_Q\right. &
+\left. (\delta(\xi,\eta))\cdot\vartheta \right)
= 0
\tag1.9.6
\\
\sum_{(\xi,\eta,\vartheta)\ \roman{cyclic}}\delta([\xi,\eta]_Q,\vartheta)
&=
\sum_{(\xi,\eta,\vartheta)\ \roman{cyclic}}\xi \cdot \delta(\eta,\vartheta)
\tag1.9.7
\endalignat
$$
where $a\in A,\ x,y \in H,\ \xi,\eta,\vartheta \in Q$.
\endproclaim

Recall \cite\extensta\ that,
given a commutative algebra $A$ and Lie-Rinehart algebras
$(A,L')$, $(A,L)$ and $(A,L'')$ where
$L'$ is an ordinary $A$-Lie algebra,
an extension of Lie-Rinehart algebras
$$
0 @>>> L' @>>> L @>>> L'' @>>> 0
$$
is an extension of $A$-modules which is also
an extension of ordinary Lie algebras
so that the projection from $L$ to $L''$ is a morphism of Lie-Rinehart algebras.
Theorem 1.9 entails at once the following.

\proclaim{Corollary 1.9.8}
Given a Lie-Rinehart triple $(A,H,Q)$,
the left $(A,H)$-module structures {\rm (1.5.2.$H$)} on $A$ 
and {\rm (1.5.3)} on $Q$ are trivial if and only if
$(A,Q)$ is a Lie-Rinehart algebra in such a way that the projection
from $E=H \oplus Q$ to $Q$ fits into an extension
$$
0 @>>> H @>>> E @>>> Q @>>> 0
$$
of Lie-Rinehart algebras. \qed
\endproclaim

Thus Lie-Rinehart triples $(A,H,Q)$ having trivial left $(A,H)$-module 
structures  on $A$ and $Q$ 
and extensions of Lie-Rinehart algebras of the kind $(A,L)$ 
together with an $A$-module section of the projection map
are equivalent notions.

\demo{Proof of Theorem {\rm (1.9)}}
The bracket (1.6.1) is plainly skew-symmetric. Hence the proof comes down to 
relating the Jacobi identity in $L$ and the Lie-Rinehart compatibility 
properties with (1.9.1)--(1.9.7).
\smallskip
Thus, suppose that the bracket $[\,\cdot\,,\cdot\,]$ on $L=H\oplus Q$ given by 
(1.6.1) and the operation $L \otimes_R A \to A$ given by (1.6.3) turn $(A,L)$
into a Lie-Rinehart algebra. Given $\xi \in Q$ and $x \in H$, we have 
$[\xi,x] = \xi \cdot x -x \cdot \xi$; since $L$ acts on $A$ by derivations, 
for $a \in A$, we conclude
$$
\xi (x(a))- x(\xi(a)) 
=[\xi,x] (a)
= (\xi \cdot x)(a) -(x \cdot \xi)(a),
$$
that is, (1.9.1) holds. Likewise, given $\xi,\eta \in Q$, 
$[\xi,\eta]= [\xi,\eta]_Q + \delta(\xi,\eta) \in L$ whence, for $a \in A$,
$$
\xi (\eta(a))- \eta(\xi(a)) 
=[\xi,\eta] (a) = [\xi,\eta]_Q(a) + (\delta(\xi,\eta))(a) ,
$$
that is, (1.9.4) holds. Next, since $L$ is a Lie algebra, its bracket 
satisfies the Jacobi identity. Hence, given $x \in H$ and $\xi,\eta \in Q$,
$$
\align
x\cdot [\xi,\eta]_Q-[\xi,\eta]_Q \cdot x 
&= [x,[\xi,\eta]_Q] =[x,[\xi,\eta]]- [x,\delta(\xi,\eta)] 
\\
&= [[x,\xi],\eta] + [\xi,[x,\eta]] - [x,\delta(\xi,\eta)] 
\\
&= [x\cdot\xi -\xi \cdot x,\eta] + [\xi,x\cdot \eta - \eta \cdot x]
- [x,\delta(\xi,\eta)]_H
\\
&=  [x\cdot\xi,\eta] 
   -[\xi \cdot x,\eta]
   +[\xi,x\cdot \eta]
   -[\xi,\eta \cdot x]
- [x,\delta(\xi,\eta)]_H
\\
&= [x\cdot\xi,\eta] + [\xi,x\cdot \eta]
\\
&\quad
  -(\xi\cdot x)\cdot\eta +\eta\cdot (\xi \cdot x)
  +(\eta\cdot x)\cdot \xi - \xi\cdot(\eta \cdot x)
\\
&\quad - [x,\delta(\xi,\eta)]_H
\\
&= [x\cdot\xi,\eta]_Q + \delta(x\cdot\xi,\eta)
 + [\xi,x\cdot \eta]_Q +\delta(\xi,x\cdot \eta)
\\
&\quad
  -(\xi\cdot x)\cdot\eta +\eta\cdot (\xi \cdot x)
  +(\eta\cdot x)\cdot \xi - \xi\cdot(\eta \cdot x)
\\
&\quad - [x,\delta(\xi,\eta)]_H
\endalign
$$
whence, comparing components in $H$ and $Q$, we conclude
$$
\align
x\cdot [\xi,\eta]_Q
&=
 [x\cdot\xi,\eta]_Q + [\xi,x \cdot \eta]_Q
  -(\xi\cdot x)\cdot\eta +(\eta\cdot x)\cdot \xi 
\\
[\xi,\eta]_Q \cdot x 
&=
 \xi\cdot(\eta \cdot x) -\eta\cdot (\xi \cdot x)
- \delta(x\cdot\xi,\eta)
-\delta(\xi,x\cdot \eta)
+[x,\delta(\xi,\eta)]_H
\endalign
$$
that is, (1.9.2) and (1.9.5) hold.
\smallskip
Likewise given $\xi \in Q$ and $x,y \in H$,
$$
\align
\xi\cdot [x,y]_H-[x,y]_H \cdot \xi 
&= [\xi,[x,y]_H] 
\\
&= [[\xi,x],y] + [x,[\xi,y]] 
\\
&= [\xi\cdot x -x \cdot \xi,y] + [x,\xi\cdot y - y \cdot \xi]
\\
&=  [\xi\cdot x,y] 
   -[x \cdot \xi,y]
   +[x,\xi\cdot y]
   -[x,y \cdot \xi]
\\
&= [\xi\cdot x,y] + [x,\xi\cdot y]
\\
&\quad
  -(x\cdot \xi)\cdot y +y\cdot (x \cdot \xi)
  +(y\cdot \xi)\cdot x - x\cdot(y \cdot \xi)
\\
&= [\xi\cdot x,y]_H  + [x,\xi\cdot y]_H 
\\
&\quad
  -(x\cdot \xi)\cdot y +y\cdot (x \cdot \xi)
  +(y\cdot \xi)\cdot x - x\cdot(y \cdot \xi)
\endalign
$$
whence, comparing components in $Q$ and $H$, we conclude
$$
\align
\xi\cdot [x,y]_H
&=
 [\xi\cdot x,y]_H + [x,\xi \cdot y]_H
  -(x\cdot \xi)\cdot y +(y\cdot \xi)\cdot x
\\
[x,y]_H \cdot \xi 
&=
 x\cdot(y \cdot \xi) -y\cdot (x \cdot \xi)
\endalign
$$
that is, (1.9.3) and (1.5.12) hold; notice that (1.5.12) holds already by 
assumption.
\smallskip
Next, given $\xi,\eta,\vartheta \in Q$, 
$$
\align
[[\xi,\eta],\vartheta]
&=
[[\xi,\eta]_H,\vartheta]
+
[[\xi,\eta]_Q,\vartheta]
\\
&=
[\delta(\xi,\eta),\vartheta]
+
[[\xi,\eta]_Q,\vartheta]
\\
&=
[\delta(\xi,\eta),\vartheta]_H
+
[\delta(\xi,\eta),\vartheta]_Q
+
[[\xi,\eta]_Q,\vartheta]_H
+
[[\xi,\eta]_Q,\vartheta]_Q
\\
&=
(\delta(\xi,\eta))\cdot\vartheta 
-\vartheta \cdot \delta(\xi,\eta)
+
\delta([\xi,\eta]_Q,\vartheta)
+
[[\xi,\eta]_Q,\vartheta]_Q
\endalign
$$
Hence
$$
\align
[[\xi,\eta],\vartheta] + [[\eta,\vartheta],\xi] + [[\vartheta,\xi],\eta]
&=
(\delta(\xi,\eta))\cdot\vartheta 
+(\delta(\eta,\vartheta))\cdot\xi 
+(\delta(\vartheta,\xi))\cdot\eta 
\\
&\quad
+
[[\xi,\eta]_Q,\vartheta]_Q
+ [[\eta,\vartheta]_Q,\xi]_Q + [[\vartheta,\xi]_Q,\eta]_Q
\\
&
\quad
-
\xi \cdot \delta(\eta,\vartheta)
-
\eta \cdot \delta(\vartheta,\xi)
-
\vartheta \cdot \delta(\xi,\eta)
\\
&
\quad
+
\delta([\xi,\eta]_Q,\vartheta)
+
\delta([\eta,\vartheta]_Q,\xi)
+
\delta([\vartheta,\xi]_Q,\eta)
\endalign
$$
Thus the Jacobi identity implies
$$ 
\align
[[\xi,\eta]_Q,\vartheta]_Q+ 
[[\eta,\vartheta]_Q,\xi]_Q + [[\vartheta,\xi]_Q,\eta]_Q 
+
(\delta(\xi,\eta))\cdot\vartheta 
+(\delta(\eta,\vartheta))\cdot\xi 
+(\delta(\vartheta,\xi))\cdot\eta 
&= 0
\\
\delta([\xi,\eta]_Q,\vartheta)
+
\delta([\eta,\vartheta]_Q,\xi)
+
\delta([\vartheta,\xi]_Q,\eta)
-
\xi \cdot \delta(\eta,\vartheta)
-
\eta \cdot \delta(\vartheta,\xi)
-
\vartheta \cdot \delta(\xi,\eta)
&=0,
\endalign
$$
that is, (1.9.6) and (1.9.7) are satisfied.
\smallskip
Conversely, suppose that the brackets $[\,\cdot\,,\cdot\,]_H$ and 
$[\,\cdot\,,\cdot\,]_Q$
on $H$ and $Q$, respectively, and the operations 
{\rm (1.5.3)}, {\rm (1.5.4)}, and
{\rm (1.5.5)}, are related by (1.9.1)--(1.9.7). We can then read the above 
calculations backwards and conclude that the bracket (1.6.1) on $L$ satisfies 
the Jacobi identity and that the operation (1.6.3) yields a Lie algebra action 
of $L$ on $A$ by derivations. The remaining Lie-Rinehart algebra axioms hold by 
assumption. Thus $(A,L)$ is then a Lie-Rinehart algebra. \qed
\enddemo

\noindent
{\smc Remark 1.10.} Under the circumstances of Example 1.3, the requirements
(1.5.6.Q), (1.5.6.H), (1.5.7.Q), (1.5.7.H), (1.5.8)--(1.5.11) are vacuous,
and so are (1.9.1) and (1.9.4) as well.
\smallskip
Given  an $(R,A)$ Lie algebra $L$ and an $(R,A)$ Lie subalgebra $H$, the 
invariants $A^H\subseteq A$ constitute a subalgebra of $A$; we will then 
denote the normalizer of $H$ in $L$ in the sense of Lie algebras by $L_H$, 
that is, $L_H$ consists of all $\alpha \in L$ having the property that
$[\alpha,\beta] \in H$ whenever $\beta \in H$.

\proclaim{Corollary 1.11}
Given a Lie-Rinehart triple $(A,H,Q)$, the corresponding $(R,A)$-Lie algebra 
being written as $L = H \oplus Q$, the intersection $Q \cap L_H$ coincides with 
the invariants $Q^H$ under the $H$-action on $Q$ (given by the corresponding 
operation {\rm (1.5.3)}), the pair $(A^H,Q^H)$ acquires  a Lie-Rinehart algebra 
structure, and the projection from $L_H$ to $Q^H$ fits into an extension
$$
0
@>>>
H
@>>>
L_H
@>>>
Q^H
@>>>
0
\tag1.11.1
$$
of $(R,A^H)$-Lie algebras. Furthermore, the restriction of $\delta$ to 
$Q^H$ is a cocycle for this extension, that is, it yields  the curvature of the
connection for the extension determined by the $A^H$-module direct sum 
decomposition $L_H = H \oplus Q^H$.
\endproclaim

Notice that $H$ is here viewed as an ordinary $A^H$-Lie algebra, the 
$H$-action on $A^H$ being trivial by construction.

\demo{Proof}
Indeed, given $\alpha \in Q$ and $\beta \in H$,
$$
[\alpha,\beta] = \alpha \cdot \beta - \beta \cdot \alpha
\in L
$$
whence $[\alpha,\beta] \in H$ for every $\beta \in H$ if and only if
$\beta \cdot \alpha = 0 \in Q$ for every $\beta \in H$, that is, if and only 
if $\alpha$ is invariant under the $H$-action on $Q$. The rest of the claim is 
an immediate consequence of Theorem 1.9. \qed
\enddemo

\noindent
{\smc 1.12. Illustration.}
Under the circumstances of Example 1.4.1, Corollary 1.11 obtains, with 
$H= L_{\Cal F}$. Now $A^H= A^{L_{\Cal F}}\subseteq A$ is the algebra of smooth 
functions which are constant on the leaves, that is, the algebra of functions 
on the \lq\lq space of leaves\rq\rq, 
and $L_H$ consists of the vector fields which 
\lq\lq project\rq\rq\ to the \lq\lq space of leaves\rq\rq. 
Indeed, given a function $f$ 
which is constant on the leaves and vector fields $X \in L_H$ and
$Y \in L_{\Cal F}$, necessarily $Y(Xf) =[Y,X] f + X(Yf) = 0$ whence $Xf$ is 
constant on the leaves as well. Thus we may view $Q^H$ as the  {\sl Lie
algebra of vector fields on the \lq\lq space of leaves\rq\rq \/}, 
that is, as the space of 
sections of a certain {\sl geometric object which serves as a replacement for 
the in general non-existant tangent bundle of the 
\lq\lq space of leaves\rq\rq\ \/}.

\smallskip
\noindent
{\smc Remark 1.13.} In analogy to the deformation theory of
complex manifolds, given a Lie-Rinehart triple $(A,H,Q)$,
we may view $H$ and $Q$ as what corresponds to the antiholomorphic
and holomorphic tangent bundle, respectively, and accordingly study
deformations of the Lie-Rinehart triple via morphisms 
$\vartheta \colon H \to Q$ and spell out the resulting infinitesimal 
obstructions. This will include a theory of deformations of foliations.
Details will be given elsewhere.

\beginsection 2. Lie-Rinehart triples and Maurer-Cartan algebras

In this section we will explore the relationship between
Lie-Rinehart triples and suitably defined Maurer-Cartan algebras. 
In particular, we will show that, under an additional assumption,
the two notions are equivalent; see Theorem 2.8.3 below for details.
As an application we will explain how the spectral sequence of a foliation 
and the Hodge- de Rham spectral sequence arise as special cases of a single 
conceptually simple construction. More applications will be given
in subsequent sections.
\smallskip\noindent
{\smc 2.1. Maurer-Cartan algebras.\/} Given an $A$-module $L$ and an 
$R$-derivation $d$ of degree $-1$ on the graded $A$-algebra 
$\roman{Alt}_A(L,A)$, we will refer to $(\roman{Alt}_A(L,A),d)$ as a 
{\it Maurer-Cartan\/} algebra (over $L$) provided $d$ has square zero, 
i.~e. is a differential.
\smallskip
Recall that a {\it multicomplex\/} (over $R$) is a bigraded $R$-module
$\{M^{p,q}\}_{p,q}$ together with an operator
$d_j \colon M^{p,q} \to M^{p+j,q-j+1}$ for every $j \geq 0$ such that the sum 
$d = d_0 + d_1 + \dots $ is a differential, i.~e. $dd=0$, cf. \cite\liuleone,
\cite\liuletwo. The idea of multicomplex occurs already in \cite\helleron\ 
and was exploited at various places in the literature including 
\cite\modpcoho, \cite\intecoho.
We note that an infinite sequence of the kind $(d_2, d_3, ...)$
is a system of higher homotopies. We will refer to a multicomplex 
$(M;d_0,d_1,d_2,\dots)$ whose underlying bigraded object $M$ is endowed with a 
bigraded algebra structure such that the operators $d_j$ are derivations with 
respect to this algebra structure as a {\it multi\/} $R$-{\it algebra\/}. 
\smallskip
Given $A$-modules $H$ and $Q$, consider the bigraded $A$-algebra
$(\roman{Alt}_A(Q, \roman{Alt}_A(H,A))$; we will refer to a multi $R$-algebra 
structure (beware: {\it not\/} multi $A$-algebra structure) on this bigraded 
$A$-algebra having at most $d_0, d_1, d_2$ non-zero as a {\it Maurer-Cartan 
algebra\/} structure. The resulting multi $R$-algebra will then be written as
$$ 
(\roman{Alt}_A(Q, \roman{Alt}_A(H,A));d_0, d_1, d_2)
\tag2.1.1
$$
and referred to as a (multi) {\it Maurer-Cartan algebra\/} (over $(Q,H)$). 
Usually we will discard \lq\lq multi\rq\rq\  and more simply refer to a 
{\it Maurer-Cartan algebra\/}. We note that, for degree reasons, when 
(2.1.1) is a Maurer-Cartan algebra, the operator $d_2$ is necessarily an 
$A$-derivation (since $d_2(a) = 0$ for every $a \in A
\cong \roman{Alt}_A^0(Q, \roman{Alt}_A^0(H,A))$).

\smallskip
\noindent{\smc Remark 2.1.2.}
In this definition, we could allow for non-zero derivations of the kind $d_j$ 
for $j \geq 3$ as well. This would lead to a more general notion of multi 
Maurer-Cartan algebra not studied here. The presence of a non-zero operator at 
most of the kind $d_2$ is an instance of a higher homotopy of a special kind 
which suffices to explain the \lq\lq quasi\rq\rq\ 
structures explored lated in the paper.

\smallskip
\noindent{\smc Remark 2.1.3.} Given a (multi) Maurer-Cartan algebra of the 
kind (2.1.1), the sum $d=d_0+d_1+d_2$ turns $\roman{Alt}_A(Q \oplus H,A)$
into a Maurer-Cartan algebra. However, not every
Maurer-Cartan structure on $\roman{Alt}_A(Q \oplus H,A)$
arises in this fashion, that is, a 
multi Maurer-Cartan algebra structure captures 
additional structure of interaction between $A$, $Q$, and $H$,
indeed, it captures essentially a Lie-Rinehart triple structure.
The purpose of the present section is to make this precise.

\smallskip
For later reference, we spell out the following, the proof of which is
immediate.

\proclaim{Proposition 2.1.4}
Given the three derivations $d_0, d_1, d_2$,
$$
(\roman{Alt}_A^*(Q, \roman{Alt}_A^*(H,A)),
d_0, d_1, d_2)
$$
is a (multi) Maurer-Cartan algebra if and only if the following identities are 
satisfied.
$$
\align
d_0d_0 &= 0
\tag2.1.4.1
\\
d_0d_1 + d_1 d_0 &= 0
\tag2.1.4.2
\\
d_0 d_2 +d_1d_1 + d_2 d_0&= 0
\tag2.1.4.3
\\
d_1 d_2 + d_2 d_1&= 0
\tag2.1.4.4
\\
d_2 d_2 &= 0 . \quad \qed
\tag2.1.4.5
\endalign
$$
\endproclaim

\smallskip\noindent
{\smc 2.2. Lie-Rinehart and Maurer-Cartan algebras.\/}
Let $A$ be a commutative $R$-algebra and $L$ an $A$-module, together with a
skew-symmetric $R$-bilinear bracket
$$
[\,\cdot\,,\cdot\,]_L \colon L \otimes_R L @>>> L
\tag2.2.1
$$
and an operation
$$
L \otimes_R A @>>> A, \quad x\otimes_R a \mapsto x(a),
\quad x\in H,\  a \in A
\tag2.2.2
$$
such that the values of the adjoint $L @>>> \roman{End}_R(A)$ lie in 
$\roman{Der}_R(A)$ and that (2.2.1), (2.2.2) and the $A$-module structure 
on $L$ satisfy the Lie-Rinehart axioms (2.2.3) and (2.2.4) below: 
$$
\alignat1
(a x) (b) &= a (x (b)), \quad a, b \in A,\  x\in L,
\tag2.2.3
\\
[x, a y]_L &= x (a) y+ a [x, y]_L, \quad a \in A,\ x,y \in L.
\tag2.2.4
\endalignat
$$
Let $M$ be a graded $A$-module, together with an operation
$$
L \otimes_R M @>>> M, \quad x\otimes m \mapsto x(m),
\quad x\in L,\  m \in m
\tag2.2.5
$$
subject to the following requirement:
For $\alpha \in  L, \,a \in A,\,m \in M$
$$
\align
(a\,\alpha)(m) &= a(\alpha(m)),
\tag2.2.6
\\
\alpha(a\,m) &=   a\,\alpha(m) + \alpha(a)\,m.
\tag2.2.7
\endalign
$$
We refer to an operation of the kind (2.2.5) as a {\it generalized\/}
$L$-{\it connection\/} on $M$.
Under these circumstances, the ordinary  {\smc Cartan-Chevalley-Eilenberg\/} 
(CCE) operator $d$ is defined, at first on the bigraded object 
$\roman{Alt}_R(L,M)$ of $M$-valued $R$-multilinear alternating forms on $L$.
Indeed, given an $R$-multilinear alternating function $f$ on $L$ of $n-1$ 
variables which is homogeneous, i.~e. the values of $f$ lie in a homogeneous 
constituent of $M$, the Cartan-Chevalley-Eilenberg (CCE) formula yields
$$
\aligned
(-1)^{|f|+1}(df)&(\alpha_1,\dots,\alpha_n)
=
\quad 
\sum_{i=1}^n (-1)^{(i-1)}
\alpha_i(f (\alpha_1, \dots\widehat{\alpha_i}\dots, \alpha_n))
\\
&\phantom{=}+\quad 
\sum_{1 \leq j<k \leq n} (-1)^{(j+k)}f(\lbrack \alpha_j,\alpha_k \rbrack,
\alpha_1, \dots\widehat{\alpha_j}\dots\widehat{\alpha_k}\dots,\alpha_n),
\endaligned
\tag2.2.8
$$
where $\alpha_1,\dots,\alpha_n \in L$ and where as usual 
\lq $\ \widehat {}\ $\rq \ indicates omission of the corresponding term.
We note that, when the values of the homogeneous alternating function $f$ on 
$L$ of $n-1$ variables lie in $M_q$, $|f| = q-n+1$. Here and below our 
convention is that, given graded objects $N$ and $M$,  a homogeneous morphism
$h \colon N_p \to M_q$ has  degree $|h| = q-p$. This is the standard grading 
on the Hom-functor for graded objects. The requirements (2.2.3), (2.2.4), 
(2.2.6), and (2.2.7) entail that the operator $d$ on $\roman{Alt}_R(L,M)$ 
passes to an $R$-linear operator on the (bi)graded $A$-submodule 
$\roman{Alt}_A(L,M)$ of $A$-multilinear functions, written here and 
henceforth as  
$$
d \colon \roman{Alt}_A(L,M) @>>> \roman{Alt}_A(L,M)
\tag2.2.9
$$
as well. The sign $(-1)^{|f|+1}$  in (2.2.8) is the appropriate one according 
to the customary Eilenberg-Koszul convention in differential homological 
algebra since (2.2.9) involves graded objects. For $M=A$, the operator $d$ is 
plainly a derivation on $\roman{Alt}_A(L,A)$.
\smallskip

Let $M_1$ and $M_2$ be graded $A$-modules endowed with generalized 
$L$-connections of the kind (2.2.5), and let 
$$
\langle\cdot,\cdot\rangle \colon M_1 \otimes_AM_2 @>>> M
$$ 
be an  $A$-module pairing which is compatible with the generalized 
$L$-connections in the sense that 
$$
x(\langle m_1,m_2 \rangle) =(x (m_1),m_2) +  (m_1,x (m_2)),
\quad 
x\in L,\  m_1 \in M_1,\ m_2 \in M_2.
$$
This pairing induces a (bi)graded pairing
$$
\roman{Alt}_A(L,M_1) \otimes_R \roman{Alt}_A(L,M_2) @>>>
\roman{Alt}_A(L,M)
\tag2.2.10
$$
which is compatible with the generalized CCE operators.

\smallskip
An $A$-module $M$ will be said to have {\it property\/} P provided for 
$x\in M$, $\phi(x) = 0$ for every $\phi \colon M \to A$ implies that $x$ 
is zero. For example, a projective $A$-module has property P, or a reflexive 
$A$-module has this property as well or, more generally, any $A$-module $M$ 
such that the canonical map from $M$ into its double $A$-dual is injective.
On the other hand, for example, for a smooth manifold $X$, the 
$C^{\infty}(X)$-module $D$ of formal (= K\"ahler) differentials does 
{\it not\/} have property P: On the real line, with coordinate $x$, consider 
the functions $f(x) = \sin x$ and $g(x) = \cos x$. The formal differential 
$d f -g dx$ is non-zero in $D$; however, the $C^{\infty}(X)$-linear maps from 
$D$ to $C^{\infty}(X)$ are the smooth vector fields, whence every such 
$C^{\infty}(X)$-linear map annihilates the formal differential $d f -g dx$.

\proclaim{Lemma 2.2.11}
When $L$ has the property P, the pair $(A,L)$, endowed with the bracket 
$[\,\cdot\,,\cdot\,]_L$ (cf. {\rm (2.2.1)}) and operation {\rm (2.2.2)} is a 
Lie-Rinehart algebra, that is, the bracket $[x,y]_L$ satisfies the Jacobi 
identity and the adjoint of {\rm (2.2.2)} is a morphism of $R$-Lie algebras, 
if and only if $(\roman{Alt}_A(L,A),d)$ is a Maurer-Cartan algebra.
\endproclaim

\demo{Proof}
A familiar calculation shows that $d$ is a differential if and only if the 
bracket $[x,y]_L$ satisfies the Jacobi identity and if
the adjoint of {\rm (2.2.2)} is a morphism of $R$-Lie algebras. 
Cf. also 2.8.5(i) below. \qed
\enddemo

\noindent
{\smc Example 2.2.12.\/} The Lie algebra $L$ of derivations of a polynomial
algebra $A$ in infinitely many indeterminates (over a field) has property
P as an $A$-module but is {\it not\/} a projective $A$-module. To include
this kind of example and others, 
it is necessary to build up the theory for modules 
having property P rather than just projective ones or even
finitely generated projective modules.

\smallskip
Let now $(A,L)$ be an (ungraded) Lie-Rinehart algebra, and let
$(\roman{Alt}_A(L,A),d)$ be the corresponding Maurer-Cartan algebra;
notice that the operator $d$ is not $A$-linear unless $L$ acts trivially on 
$A$. For reasons explained in \cite\lradq\ 
we will refer to this operator as {\it Lie-Rinehart\/} differential. We will 
say that the graded $A$-module $M$, endowed with the operation (2.2.5), is a 
{\it graded\/} ({\it left\/}) $(A,L)$-{\it module\/} provided this operation 
is an ordinary Lie algebra action on $M$. When $M$ is concentrated in degree 
zero, we simply refer to $M$ as a ({\it left\/}) $(A,L)$-{\it module\/}. In 
particular, with the obvious $L$-module structure, the  algebra $A$ itself is 
a (left) $(A,L)$-module. The proof of the following is straightforward and 
left to the reader.

\proclaim{Lemma 2.2.13}
When $(A,L)$ is a Lie-Rinehart algebra and when $M$ has the property P, the 
operation {\rm (2.2.5)} turns  $M$ into a left $(A,L)$-module if and only 
if the operator $d$ on $\roman{Alt}_A(L,M)$ turns $(\roman{Alt}_A(L,M),d)$ 
into a differential graded $(\roman{Alt}_A(L,A),d)$-module
via {\rm (2.2.10)} (with $M_1=A$ and $M_2=M$). \qed
\endproclaim

Given  a graded $(A,L)$-module $M$, 
we will refer to the resulting (co)chain complex
$$
(\roman{Alt}_A(L,M), d)
\tag2.2.14
$$
as the {\it Rinehart complex of\/} $M$-{\it valued forms on\/} $L$; often we 
write this complex more simply in the form $\roman{Alt}_A(L,M)$. It inherits 
a differential graded $\roman{Alt}_A(L,A)$-module structure  
via (2.2.10).
\smallskip
We now spell out the passage 
from Maurer-Cartan algebras to Lie-Rinehart algebras.

\proclaim{Lemma 2.2.15} Let $L$ be a 
finitely generated projective $A$-module. Then an 
$R$-derivation $d$ on the graded $A$-algebra $\roman{Alt}_A(L,A)$ determines 
a skew-symmetric $R$-bilinear bracket $[\,\cdot\,,\cdot\,]_L$ on $L$ of the 
kind {\rm (2.2.1)} and an operation $L \otimes A \to A$ of the  kind 
{\rm (2.2.2)} such that the identities {\rm (2.2.3)} and {\rm (2.2.4)} 
are satisfied and that the corresponding CCE operator {\rm (2.2.8)} (for
$M=A$) coincides with $d$. Furthermore, $\roman{Alt}_A(L,A)$ is
then a Maurer-Cartan algebra if and only if $(A,L)$ is a Lie-Rinehart algebra.
\endproclaim

\demo{Proof} The operator
$$
d\colon \roman{Alt}_A^q(L,A) @>>> \roman{Alt}_A^{q+1}(L,A) \quad (q \geq 0)
$$
induces, for $q=0$, an operation $L \otimes_R A @>>> A$ of the kind (2.2.2) 
and, for $q=1$, a skew-symmetric $R$-bilinear bracket $[\,\cdot\,,\cdot\,]_L$ 
on $L$ of the kind (2.2.1). More precisely: Given $x \in L$ and $a \in A$, let
$$
x(a) = - (d(a))(x).
$$
This yields an operation of the kind (2.2.2). Given $x,y \in L$, using 
the hypothesis that $L$ is 
a finitely generated projective 
$A$-module,
identify $x$ and $y$ with their 
images in the double $A$-dual $L^{**}$ and define the value $[x,y]_L$ by
$$
[x,y]_L(\alpha) = x(\alpha(y)) -y(\alpha(x)) -d\alpha(x,y)
$$
where $\alpha \in L^*=\roman{Hom}_A(L,A)$. This yields a bracket of the kind 
(2.2.1), that is, an $R$-bilinear (beware: not $A$-bilinear) skew-symmetric 
bracket on $L$. Notice that, at this stage, the operation of the kind (2.2.2) 
is already defined whence the definition of the bracket makes sense. Since, by 
assumption, $d$ is a derivation on $\roman{Alt}_A(L,A)$, the identities 
{\rm (2.2.3)} and {\rm (2.2.4)} are satisfied.
By construction, the resulting CCE operator coincides with $d$ in degree $0$ 
and in degree $-1$ whence the two operators coincide. Since a 
finitely generated projective 
$A$-module has property P, Lemma 2.2.11 completes the proof. \qed
\enddemo
                                               
Combining Lemma 2.2.13 and 2.2.15,  we arrive at the following.

\proclaim{Theorem 2.2.16} Given a finitely generated
projective $A$-module $L$, Lie-Rinehart 
algebra structures on $(A,L)$ and Maurer-Cartan algebra structures on 
$\roman{Alt}_A(L,A)$ are equivalent notions. \qed
\endproclaim

\smallskip\noindent
{\smc 2.3. Connections.\/} Let $(A,L)$ be a Lie-Rinehart algebra.
Given a graded $A$-module $M$, a 
degree zero operation $L\otimes_R M \to M$, not necessarily a graded left 
$L$-module structure but still satisfying (2.2.6) and (2.2.7), is referred to 
as an $(A,L)$-{\it connection\/}, cf. \cite{\poiscoho,\,\bv} or, somewhat more 
precisely, as a {\it graded left\/} $(A,L)$-{\it connection\/}; in this 
language, a (graded) $(A,L)$-module structure is a (graded) {\it flat\/} 
$(A,L)$-{\it connection\/}. Given a graded $A$-module $M$, together with a 
graded $(A,L)$-connection, we extend the definition of the Lie-Rinehart 
operator to an operator
$$
d\colon  \roman{Alt}_A(L,M) @>>> \roman{Alt}_A(L,M)
\tag2.3.1
$$ 
by means of the formula (2.2.8).
The resulting operator $d$ is well defined; it is a 
differential if and only if the $(A,L)$-connection on $M$ is flat, 
i.~e. an ordinary $(A,L)$-module structure.
\smallskip
\noindent
{\smc 2.4. From Lie-Rinehart triples to Maurer-Cartan algebras.\/} Let 
$(A,H,Q)$ be an almost pre-Lie-Rinehart triple. Consider the bigraded $A$-module
$$
\roman{Alt}_A^{*,*}(Q \oplus H,A)
\cong
\roman{Alt}_A^*(Q, \roman{Alt}_A^*(H,A)).
\tag2.4.1
$$
Henceforth we spell out a particular homogeneous constituent of bidegree 
$(p,q)$ (according to the conventions used below,
such a homogenous constituent will be of bidegree $(-p,-q)$ but for the moment
this usage of negative degrees is of no account) in the form
$$
\roman{Alt}_A^p(Q, \roman{Alt}_A^q(H,A)).
\tag2.4.2
$$
The operations (1.5.3) and (1.5.4) induce degree zero operations
$$
\align
H \otimes_R\roman{Alt}_A^*(Q,A) &@>>> \roman{Alt}_A^*(Q,A)
\tag2.4.3
\\
Q\otimes_R \roman{Alt}_A^*(H,A) &@>>> \roman{Alt}_A^*(H,A)
\tag2.4.4
\endalign
$$
on $\roman{Alt}_A^*(Q,A)$ and $\roman{Alt}_A^*(H,A)$, respectively, when 
(1.5.3) and (1.5.4) are treated like connections. 
By evaluation of the expression given
on the right-hand side of (2.2.8), with (1.5.1.$H$) and (1.5.1.$Q$) instead 
of (2.2.1),
and with 
(2.5.1) and (2.5.2) instead of (2.2.5),
these operations, in turn, induce two operators
$$
\align
d_0 \colon \roman{Alt}_A^p(Q, \roman{Alt}_A^q(H,A))
&@>>>
\roman{Alt}_A^p(Q, \roman{Alt}_A^{q+1}(H,A))
\tag2.4.5
\\
d_1 \colon \roman{Alt}_A^p(Q, \roman{Alt}_A^q(H,A))
&@>>>
\roman{Alt}_A^{p+1}(Q, \roman{Alt}_A^q(H,A)).
\tag2.4.6
\endalign
$$
A little thought reveals that, in view of (1.5.6.$H$), (1.5.6.$Q$), (1.5.7.$H$),
(1.5.7.$Q$), (1.5.8)--(1.5.11), these operators, which are {\sl at first\/} 
defined only on the $R$-multilinear alternating functions, in fact {\sl pass 
to operators on $A$-multilinear alternating functions\/}. Furthermore, the 
skew-symmetric $A$-bilinear pairing $\delta$, cf. (1.5.5), induces an operator
$$
d_2 \colon \roman{Alt}_A^p(Q, \roman{Alt}_A^q(H,A)) @>>>
\roman{Alt}_A^{p+2}(Q, \roman{Alt}_A^{q-1}(H,A)).
\tag2.4.7
$$
Hence, 
when $(A,H,Q)$ is a Lie-Rinehart triple, that is, when (1.6.1) and 
(1.6.3) turn $(A,H\oplus Q)$ into a Lie-Rinehart algebra,
${(\roman{Alt}_A(Q, \roman{Alt}_A(H,A)); d_0, d_1, d_2)}$
is a (multi) Maurer-Cartan algebra.
\smallskip\noindent
{\smc 2.5. Explicit description of the operators $d_0,d_1,d_2$:\/} 
Let $f$ be an alternating $A$-multilinear function on $Q$ of $p$ variables 
with values in $\roman{Alt}_A^q(H,A)$, so that $|f| = -q-p$ and 
$(-1)^{|f|+1} =(-1)^{p+q+1}$. Let $\xi_1,\dots,\xi_{p+2} \in Q$ and 
$x_1,\dots,x_{q+1} \in H$. 
\smallskip
\noindent{\smc The operator $d_0$:\/}
$$
\aligned
{}&(-1)^{p+q+1} \left((d_0f)(\xi_1,\dots,\xi_p)\right)(x_1,\dots,x_{q+1})
=
\\
&\phantom{=+} 
\sum_{j=1}^{q+1} 
(-1)^{p+j-1} x_j\left(
\left(f(\xi_1,\dots,\xi_p)\right)
(x_1,\dots\widehat{x_j}\dots,x_{q+1})
                \right)
\\
&\phantom{=}+ 
\sum_{1 \leq j<k\leq q+1} 
(-1)^{p+j+k} 
\left(f(\xi_1,\dots,\xi_p)\right)
(\lbrack x_j,x_k \rbrack_H, x_1,
\dots\widehat{x_j}\dots\widehat{x_k}\dots,x_{q+1})
\\
&\phantom{=}+ 
\sum_{j=1}^p \sum_{k=1}^{q+1} 
(-1)^{j+k+p+1} 
\left(f(x_k\cdot \xi_j, \xi_1,\dots\widehat{\xi_j}\dots,\xi_p)\right)
(x_1,\dots\widehat{x_k}\dots, x_{q+1})
\endaligned
\tag2.5.1
$$
The last term involving the double summation
necessarily appears since, for 
$1 \leq j \leq p$ and $1 \leq k \leq q+1$, the bracket $[x_k,\xi_j]$ in 
$Q \oplus H$, cf. (1.6.2), is given by
$$
[x_k,\xi_j] = x_k \cdot \xi_j - \xi_j \cdot x_k.
$$

\noindent{\smc Remark 2.5.2.} A crucial observation is this:
The operator $d_0$ may be written as the sum
$$
d_0 = d_H + d_Q
$$
of certain operators $d_H$ and  $d_Q$ defined on 
$\roman{Alt}_R(Q, \roman{Alt}_R(H,A))$ by
$$
\aligned
{}&(-1)^{p+q+1} \left((d_Hf)(\xi_1,\dots,\xi_p)\right)(x_1,\dots,x_{q+1})
=
\\
&\phantom{=+} 
\sum_{j=1}^{q+1} 
(-1)^{p+j-1}
x_j\left(
\left(f(\xi_1,\dots,\xi_p)\right)
(x_1,\dots\widehat{x_j}\dots,x_{q+1})\right)
\\
&\phantom{=}+ 
\sum_{1 \leq j<k\leq q+1} 
(-1)^{j+k} 
\left(f(\xi_1,\dots,\xi_p)\right)(\lbrack x_j,x_k \rbrack_H, x_1,
\dots\widehat{x_j}\dots\widehat{x_k}\dots,x_{q+1})
\\
{}&(-1)^{p+q+1}\left((d_Qf)(\xi_1,\dots,\xi_p)\right)(x_1,\dots,x_{q+1})
=\\
&\phantom{=} 
\sum_{1 \leq j \leq p, 1 \leq k \leq q+1}
(-1)^{j+k+p+1}
\left(f(x_k\cdot \xi_j, \xi_1,\dots\widehat{\xi_j}\dots,\xi_p)\right)
(x_1,\dots\widehat{x_k}\dots, x_{q+1}) .
\endaligned
$$
However, even when $(A,H,Q)$ is a (pre-)Lie-Rinehart triple, the individual 
operators $d_H$ and $d_Q$ are well defined merely on 
$\roman{Alt}_R(Q, \roman{Alt}_R(H,A))$; {\it only their sum is well defined\/} 
on $\roman{Alt}_A(Q, \roman{Alt}_A(H,A))$.

\smallskip
\noindent{\smc The operator $d_1$:\/}
$$
\aligned
{}&(-1)^{p+q+1} \left((d_1f)(\xi_1,\dots,\xi_{p+1})\right)(x_1,\dots,x_q)
=
\\
&\phantom{=+} 
\sum_{j=1}^{p+1} 
(-1)^{j-1} 
\xi_j\left(
\left(f(\xi_1,\dots\widehat{\xi_j}\dots,\xi_{p+1})\right)
(x_1,\dots,x_q)
     \right)
\\
&\phantom{=}+ 
\sum_{1 \leq j<k\leq p} 
(-1)^{j+k} 
\left(f(\lbrack \xi_j,\xi_k \rbrack_Q,
\xi_1\dots\widehat{\xi_j}\dots\widehat{\xi_k}\dots,\xi_{p+1})\right)
(x_1,\dots, x_q)
\\
&\phantom{=}+ 
\sum_{j=1}^{p+1} \sum_{k=1}^q
(-1)^{j+k+1} 
\left(f(\xi_1,\dots\widehat{\xi_j}\dots,\xi_{p+1})\right)
(\xi_j \cdot x_k,
x_1,\dots\widehat{x_k}\dots, x_q)
\endaligned
\tag2.5.3
$$
The last term involving the double summation
necessarily appears in view of (1.6.4). 
With the generalized operation of {\it Lie-derivative\/}
$$
(\xi,\alpha)\longmapsto \xi(\alpha),
\quad \xi \in Q,\ \alpha \in \roman{Alt}_A^q(H,A)\quad (q \geq 0)
$$
which, for 
$x_1,\dots,x_q \in H$, is given by
$$
(\xi(\alpha))(x_1,\dots,x_q) = \xi(\alpha(x_1,\dots,x_q))
- \sum_{k=1}^q \alpha(x_1,\dots,x_{k-1}, \xi \cdot x_k, x_{k+1},\dots, x_q),
$$
the identity (2.5.3) may be written as
$$
\aligned
{}&(-1)^{p+q+1} (d_1f)(\xi_1,\dots,\xi_{p+1})
=(-1)^{|f|+1} (d_1f)(\xi_1,\dots,\xi_{p+1})
\\
& =  
\sum_{j=1}^{p+1} 
(-1)^{j-1} 
\xi_j \left(f(\xi_1,\dots\widehat{\xi_j}\dots,\xi_{p+1})\right)
\\
&\phantom{=}+ 
\sum_{1 \leq j<k\leq p} 
(-1)^{j+k} 
f(\lbrack \xi_j,\xi_k \rbrack_Q,
\xi_1\dots\widehat{\xi_j}\dots\widehat{\xi_k}\dots,\xi_{p+1}) .
\endaligned
\tag2.5.3$'$
$$

\noindent{\smc The operator $d_2$:\/}
$$
\aligned
{}&(-1)^{p+q+1}\left((d_2f)(\xi_1,\dots,\xi_{p+2})\right)(x_1,\dots,x_{q-1})
=
\\
&\phantom{=+} 
\sum_{1 \leq j<k\leq p+2} 
(-1)^{j+k+p}
\left(f(
\xi_1, \dots\widehat{\xi_j}\dots\widehat{\xi_k}\dots,\xi_{p+2})\right) 
(\delta (\xi_j,\xi_k), x_1,\dots, x_{q-1})
\endaligned
\tag2.5.4
$$
\noindent{\smc Remark 2.5.5.}
The operator $d_2$ does not involve the pieces of structure (1.5.1.$H$), 
(1.5.1.$Q$), (1.5.2.$H$), (1.5.2.$Q$), (1.5.3), (1.5.4). Hence, for an 
arbitrary $A$-module $M$, the formula (2.5.4) given above yields an operator
$$
d_2 \colon \roman{Alt}_A^p(Q, \roman{Alt}_A^q(H,M)) @>>>
\roman{Alt}_A^{p+2}(Q, \roman{Alt}_A^{q-1}(H,M))
\quad(p \geq 0,\ q \geq 1).
\tag2.4.7$'$
$$
We will use this observation in (5.8.7) and (5.8.8) below.

\smallskip\noindent
{\smc Remark 2.6.} Given an almost pre-Lie-Rinehart triple $(A,H,Q)$, the 
vanishing of $d_2 d_2$ is automatic, for the following reason: View $H$ and $Q$ 
as abelian $A$-Lie algebras and $H$ as being endowed with the trivial 
$Q$-module structure. Since $\delta$ is a skew-symmetric $A$-bilinear pairing, 
we may use it to endow the $A$-module direct sum $L= H \oplus Q$ with a 
nilpotent $A$-Lie algebra structure (of class two) by setting
$$
[(x,\xi),(y,\eta)] = (\delta(\xi,\eta),0),\quad \xi,\eta \in Q,\ x,y \in H.
$$
We write $L_{\roman{nil}}$ for this nilpotent $A$-Lie algebra. The ordinary 
CCE complex for calculating the Lie algebra cohomology 
$\roman H^*(L_{\roman{nil}},A)$ (with trivial $L_{\roman{nil}}$-action on $A$) 
is just $(\roman{Alt}_A(L,A),d_2)$. Thus the vanishing of $d_2 d_2$ is 
automatic.

\proclaim{Theorem 2.7}
An almost pre-Lie-Rinehart triple $(A,H,Q)$ such that $H$ and $Q$ have 
property {\rm P} is a Lie-Rinehart triple, that is, {\rm (1.6.1)} 
and {\rm (1.6.3)} then turn $(A,H\oplus Q)$ into a Lie-Rinehart algebra,
if and only if
${(\roman{Alt}_A(Q, \roman{Alt}_A(H,A));d_0, d_1, d_2)}$
is a (multi) Maurer-Cartan algebra.
\endproclaim

\demo{Proof} The direct $A$-module sum $L=Q \oplus H$ has the property P.
The sum $d=d_0+ d_1+ d_2$ is an $R$-derivation on $\roman{Alt}_A(L,A)$.
Hence the claim is an immediate consequence of Lemma 2.2.11. \qed
\enddemo

\smallskip
\noindent
{\smc 2.8. From Maurer-Cartan algebras to Lie-Rinehart triples.\/}
Let $H$ and $Q$ be 
finitely generated projective
$A$-modules, and let $d_0, d_1, d_2$ 
be homogeneous $R$-derivations of the bigraded $A$-algebra 
$\roman{Alt}_A(Q, \roman{Alt}_A(H,A))$ of the kind
$$
d_j\colon \roman{Alt}_A^p(Q, \roman{Alt}_A^q(H,A))
@>>>
\roman{Alt}_A^{p+j}(Q, \roman{Alt}_A^{q-j+1}(H,A)).
$$

\proclaim{Proposition 2.8.1}
The operators $d_0, d_1, d_2$ induce an almost pre-Lie-Rinehart triple
structure on $(A,H,Q)$.
\endproclaim

\demo{Proof} Write $L=Q \oplus H$. The sum $d=d_0+ d_1+ d_2$ 
is a derivation on $\roman{Alt}_A(L,A)$.
By Lemma 2.2.15, $d$  induces a bracket $[\,\cdot\,,\cdot\,]_L$ 
on $L$ (of the kind (2.2.1)) and an operation $L \otimes_R A @>>> A$
of the kind (2.2.2).
Taking homogeneous components with reference to the 
direct sum decomposition $L=Q \oplus H$, we obtain an 
almost pre-Lie-Rinehart triple structure of the kind
{\rm (1.5.1.$H$)}, {\rm (1.5.2.$H$)}, {\rm (1.5.1.$Q$)}, {\rm (1.5.2.$Q$)},
{\rm (1.5.3)}, {\rm (1.5.4)}, {\rm (1.5.5)}
on $(A,H,Q)$. The three almost pre-Lie-Rinehart triple axioms
are implied by the fact that the operators
$d_0, d_1, d_2$ are derivations of the 
bigraded algebra $\roman{Alt}_A(Q, \roman{Alt}_A(H,A))$.\qed
\enddemo

\proclaim{Theorem 2.8.2} 
The triple $(A,H,Q)$, endowed with the induced operations of the kind 
{\rm (1.5.1.$H$)}, {\rm (1.5.2.$H$)}, {\rm (1.5.1.$Q$)}, {\rm (1.5.2.$Q$)},
{\rm (1.5.3)}, {\rm (1.5.4)}, {\rm (1.5.5)} given in {\rm (2.8.1)} above, 
is a pre-Lie-Rinehart triple if and only if $d_0$ is a differential; 
$(A,H,Q)$ is a Lie-Rinehart triple if and only if
$(\roman{Alt}_A(Q, \roman{Alt}_A(H,A));d_0, d_1, d_2)$
is a Maurer-Cartan algebra.
\endproclaim

\demo{Proof} This is a consequence of Lemmata 2.2.13 and 2.2.15. \qed
\enddemo

Combining Theorem 2.7 and Theorem 2.8.2, we arrive at the following.

\proclaim{Theorem 2.8.3} Given 
finitely generated projective $A$-modules $H$ and $Q$, 
Lie-Rinehart triple structures on $(A,H,Q)$
and (multi) Maurer-Cartan algebra structures on
$\roman{Alt}_A(Q, \roman{Alt}_A(H,A))$
are equivalent notions. \qed
\endproclaim

\noindent
{\smc Remark 2.8.4.\/} Concerning the hypotheses and hence
the range of applications, cf. e.~g. Example 2.2.12 above,
Theorem 2.7 is somewhat more general
than Theorem 2.8.3. This justifies, hopefully, the terminology
\lq\lq almost-\rq\rq\ and \lq\lq pre-Lie-Rinehart triple\rq\rq,
admittedly a bit cumbersome.
In fact, it would be interesting and important to establish the statement of
Theorem 2.8.3 for $A$-modules more general than finitely generated
and projective.

\smallskip
\noindent
{\smc 2.8.5. Direct verification of the Lie-Rinehart triple structure.\/}
Let $(A,H,Q)$ be an almost pre-Lie-Rinehart triple such that $H$ and $Q$ 
have property {\rm P}, and suppose that 
$(\roman{Alt}_A(Q, \roman{Alt}_A(H,A));d_0, d_1, d_2)$
is a (multi) Maurer-Cartan algebra. It is then instructive to deduce directly 
that $(A,H,Q)$  is a Lie-Rinehart triple.
\smallskip\noindent
(i) Consider the operator
$$
d_0d_0
\colon
\roman{Alt}_A^j(H, A)
@>>>
\roman{Alt}_A^{j+2}(H, A)
$$
for $j=0$ and $j=1$. Notice that $\roman{Alt}_A^j(H, A)$ equals
$\roman{Alt}_A^j(H, \roman{Alt}_A^0(Q,A))$ and that
$\roman{Alt}_A^{j+2}(H, A)$ equals 
$\roman{Alt}_A^{j+2}(H, \roman{Alt}_A^0(Q,A))$. For $j=1$, given $x,y,z \in H$ 
and  $\phi \in \roman{Hom}_A(H,A)=\roman{Alt}_A^j(H, A)$, we find
$$
(d_0d_0 \phi)(x,y,z)
=
\phi(
[[x,y]_H,z]_H +
[[y,z]_H,x]_H +
[[z,x]_H,y]_H). 
$$
Since $H$ has property P, we conclude that the bracket on $H$ satisfies 
the Jacobi identity, that is, $H$ is an $R$-Lie algebra. Likewise, for $j=0$,
given $x,y\in H$ and $a \in A$, we find
$$
(d_0d_0 a)(x,y)
=
x (y (a))-y(x (a)) - [x,y](a).
$$
Consequently the adjoint $H \to \roman{Der}_R(A)$ of (1.5.2.$H$) is a morphism 
of $R$-Lie algebras. In view of (1.5.6.$H$) and (1.5.7.$H$), we conclude 
that 
(1.5.1.$H$) and (1.5.2.$H$) turn
$(A,H)$ into a Lie-Rinehart algebra.
\smallskip\noindent
(ii) Next,
consider the operator
$$
d_0d_0
\colon
\roman{Alt}_A^1(Q, \roman{Alt}_A^0(H,A))
@>>>
\roman{Alt}_A^1(Q, \roman{Alt}_A^2(H,A)).
$$
We note that
$\roman{Alt}_A^1(Q, \roman{Alt}_A^0(H,A)) = \roman{Alt}_A^1(Q,A)=
\roman{Hom}_A(Q,A)$.
Let $\xi\in Q$, $x,y \in H$, and $\phi \in \roman{Hom}_A(H,A)$.
A straightforward calculation gives
$$
((d_0d_0 \phi)(\xi))(x,y)
=
\phi(y \cdot (x \cdot \xi) - x \cdot (y \cdot \xi)
+[x,y]_H \cdot \xi).
$$
Since $H$ is assumed to have property P, we conclude that, for every 
$\xi\in Q,\, x,y \in H$,
$$
[x,y]_H \cdot \xi =  x \cdot (y \cdot \xi)-y \cdot (x \cdot \xi),
$$
that is, (1.5.3) is a left $(A,H)$-module structure on $Q$.
\smallskip\noindent
(iii)
Pursuing the same kind of reasoning, consider the operator
$$
d_0d_1+d_1d_0
\colon
A=\roman{Alt}_A^0(Q, \roman{Alt}_A^0(H,A))
@>>>
\roman{Alt}_A^1(Q, \roman{Alt}_A^1(H,A)).
$$
Let $a \in A,\,\xi \in Q,\,x \in H$. Again a calculation shows that
$$
((d_0d_1+d_1d_0)a)(\xi)(x)
=
x (\xi(a))- \xi(x(a)) - ((x \cdot \xi)(a) -(\xi \cdot x)(a))
$$
whence the vanishing of $d_0d_1+d_1d_0$ in bidegree $(0,0)$ entails the 
compatibility property (1.9.1). Likewise consider the operator
$$
d_0d_1+d_1d_0
\colon
\roman{Hom}_A(H,A) =\roman{Alt}_A^0(Q, \roman{Alt}_A^1(H,A))
@>>>
\roman{Alt}_A^1(Q, \roman{Alt}_A^2(H,A)).
$$
Again a calculation shows that, for 
$\xi \in Q,\,x,y \in H, \phi\in \roman{Hom}_A(H,A)$,
$$
\align
((d_0d_1+d_1d_0)\phi)(\xi)(x,y)
&=
\phi\left(\xi\cdot[x,y]_H \right .
\\
&
\quad
\left .
-([\xi\cdot x,y]_H
+
[x,\xi\cdot y]_H
-(x \cdot \xi)\cdot y
+(y \cdot \xi)\cdot x)\right)
\endalign
$$
whence the vanishing of $d_0d_1+d_1d_0$ in bidegree $(0,1)$ entails the 
compatibility property (1.9.3). Likewise, the vanishing of the operator
$d_0d_1+d_1d_0$ in bidegree $(1,0)$, that is, of 
$$
d_0d_1+d_1d_0
\colon
\roman{Alt}_A^1(Q, \roman{Alt}_A^0(H,A))
@>>>
\roman{Alt}_A^2(Q, \roman{Alt}_A^1(H,A)),
$$
entails the compatibility property (1.9.2). 
\smallskip
In the same vein:
\newline\noindent
(iv) The vanishing of the operator
$$
d_1d_1 +d_2d_0 =
d_0d_2+d_1d_1 +d_2d_0
\colon
\roman{Alt}_A^0(Q, \roman{Alt}_A^0(H,A))
@>>>
\roman{Alt}_A^2(Q, \roman{Alt}_A^0(H,A))
$$
entails the compatibility property (1.9.4). 
\newline\noindent
(v) The vanishing of the operator
$$
d_1d_1 +d_2d_0 =
d_0d_2+d_1d_1 +d_2d_0
\colon
\roman{Alt}_A^1(Q, \roman{Alt}_A^0(H,A))
@>>>
\roman{Alt}_A^3(Q, \roman{Alt}_A^0(H,A)),
$$
together with (1.9.4), entails the compatibility property (1.9.6),
the \lq\lq generalized Jacobi identity for the bracket 
$[\,\cdot\,,\cdot\,]_Q$\rq\rq. 
For intelligibility and later reference
(cf. (4.10) and (6.11) below), we sketch the argument:
Let $\alpha \in \roman{Alt}_A^1(Q, \roman{Alt}_A^0(H,A))$
and $\xi,\eta,\vartheta \in Q$. A straightforward calculation yields
$$
\align
(d_1 d_1\alpha)(\xi,\eta,\vartheta)
&=
-\sum_{(\xi,\eta,\vartheta)\ \roman{cyclic}}
\alpha([[\xi,\eta]_Q,\vartheta]_Q
\\
&
\quad
-\sum_{(\xi,\eta,\vartheta)\ \roman{cyclic}}
\left(
\xi(\eta(\alpha(\vartheta)))
-
\eta(\xi(\alpha(\vartheta)))
-[\xi,\eta]_H(\alpha(\vartheta))
\right)
\endalign
$$
Using (1.9.4), we substitute
$(\delta(\xi,\eta))(\alpha(\vartheta))$
for
$\xi(\eta(\alpha(\vartheta)))
-
\eta(\xi(\alpha(\vartheta)))
-[\xi,\eta]_H(\alpha(\vartheta))$
and
obtain
$$
(d_1 d_1\alpha)(\xi,\eta,\vartheta)
=
-\sum_{(\xi,\eta,\vartheta)\ \roman{cyclic}}
\alpha([[\xi,\eta]_Q,\vartheta]_Q
-\sum_{(\xi,\eta,\vartheta)\ \roman{cyclic}}
(\delta(\xi,\eta))(\alpha(\vartheta))
$$
Likewise,
a calculation gives
$$
(d_2 d_0 \alpha)(\xi,\eta,\vartheta)
=
\sum_{(\xi,\eta,\vartheta)\ \roman{cyclic}}
(\delta(\xi,\eta))(\alpha(\vartheta))
-
\sum_{(\xi,\eta,\vartheta)\ \roman{cyclic}}
\alpha(\delta(\xi,\eta) \cdot\vartheta)
$$
whence
the vanishing of the operator
$d_1d_1 +d_2d_0$ 
on $\roman{Alt}_A^1(Q, \roman{Alt}_A^0(H,A))$ implies
$$
\sum_{(\xi,\eta,\vartheta)\ \roman{cyclic}}
\alpha \left([[\xi,\eta]_Q,\vartheta]_Q
+(\delta(\xi,\eta)) \cdot\vartheta\right)
=0.
$$
For later reference we note that
$$
(d_2 d_0 \alpha(\vartheta))(\xi,\eta)
= (\delta(\xi,\eta)) (\alpha(\vartheta))
$$
whence
$$
\sum_{(\xi,\eta,\vartheta)\ \roman{cyclic}}
\alpha([[\xi,\eta]_Q,\vartheta]_Q)
=
(d_2 d_0 \alpha)(\xi,\eta,\vartheta)
+
\sum_{(\xi,\eta,\vartheta)\ \roman{cyclic}}
(d_2 d_0 \alpha(\vartheta))(\xi,\eta)
\tag2.8.6
$$
\newline\noindent
(vi) The vanishing of the operator
$$
d_0d_2+d_1d_1 +d_2d_0
\colon
\roman{Alt}_A^0(Q, \roman{Alt}_A^1(H,A))
@>>>
\roman{Alt}_A^2(Q, \roman{Alt}_A^1(H,A))
$$
entails the compatibility property (1.9.5),
the \lq\lq generalized $Q$-module structure on $H$\rq\rq.
\newline\noindent
(vii) The vanishing of the operator
$$
d_1d_2+d_2d_1
\colon
\roman{Alt}_A^0(Q, \roman{Alt}_A^1(H,A))
@>>>
\roman{Alt}_A^3(Q, \roman{Alt}_A^0(H,A))
$$
entails the compatibility property (1.9.7). Indeed, given
$\xi,\eta,\vartheta \in Q$ and $\alpha \colon H \to A$,
$$
((d_1d_2+d_2d_1)\alpha)(\xi,\eta,\vartheta)
=
\sum_{(\xi,\eta,\vartheta)\quad\roman{cyclic}}
\alpha\left(\delta([\xi,\eta]_Q,\vartheta)
-\xi \cdot \delta(\eta,\vartheta)\right) .
$$

\smallskip\noindent
{\smc 2.9. The spectral sequence.\/} Let $(A,H,Q)$ be a Lie-Rinehart triple. 
The filtration of $\roman{Alt}_A(Q, \roman{Alt}_A(H,A))$ by $Q$-degree
leads to a spectral sequence
$$
(\roman E_r^{*,*},d_r)
\tag2.9.1
$$
having
$$
(\roman E_0,d_0) = (\roman{Alt}_A(Q, \roman{Alt}_A(H,A)), d_0)
\tag2.9.2
$$
whence $\roman E^{p,q}_1$ amounts to the Lie-Rinehart cohomology 
$\roman H^q(H,\roman{Alt}^p_A(Q,A))$ of $H$ with values in the left 
$(A,H)$-module $\roman{Alt}^p_A(Q,A)$. 
There is a slight conflict of notation here but it will always be 
clear from the context whether
$d_j$ ($j \geq 0$) refers to the differentials of a spectral sequence
or to a system of multicomplex operators.
{\sl The spectral sequence\/}
(2.9.1)
{\it is an 
invariant of the Lie-Rinehart triple structure.\/} In particular,
$\roman E_1^{0,0} = A^H$ and $\roman E_1^{1,0} = \roman{Hom}(Q,A)^H$, and 
$\roman H^*(H,A)$ inherits an $(A^H,Q^H)$-module structure, with reference to 
the Lie-Rinehart structure on $(A^H,Q^H)$, cf. Corollary 1.11. Thus the 
Rinehart complex $(\roman{Alt}_{A^H}(Q^H, \roman H^*(H,A)),d)$ is defined.

\smallskip\noindent
{\smc 2.10. Illustration.}
The spectral sequence (2.9.1) includes as special cases
that of a foliation and the Hodge-de Rham spectral sequence.
This provides a conceptually simple approach
to these spectral sequences and subsumes them 
under a single more general 
construction.
We will now make this precise. 
\newline\noindent
(i) Consider a 
foliated manifold $M$, the foliation being written as
$\Cal F$. Recall that a 
$p$-form $\omega$ on $M$ is called {\it horizontal\/} (with reference to the 
foliation $\Cal F$) provided $\omega(X_1,\dots,X_p) = 0$ if some $X_j$ is 
vertical, i.~e. tangent to the foliation, or, equivalently, $i_X\omega = 0$ 
whenever $X$ is vertical; a horizontal $p$-form $\omega$ is said to be 
{\it basic\/} provided it is constant on the leaves (i.~e. 
$\lambda_X\omega = 0$ whenever $X$ is vertical). The sheaf of germs of basic 
$p$-forms is in general not fine and hence gives rise to in general non-trivial 
cohomology in non-zero degrees, cf. \cite\breinhar. Thus, under the 
circumstances of the Example 1.4.1, and those of (1.12) as well, so that 
$(A,H)$ is the Lie-Rinehart algebra $(C^{\infty}(M),L_{\Cal F})$ arising from a 
foliation $\Cal F$ of a smooth manifold $M$, for every $p \geq 0$, the Rinehart 
complex $(\roman{Alt}^*_A(H,\roman{Alt}^p(Q,A)),d)$ for the Lie-Rinehart 
algebra $(A,H)=(C^{\infty}(M),L_{\Cal F})$ with coefficients in 
$\roman{Alt}^p(Q,A)$ which computes the cohomology 
$\roman H^*(L_{\Cal F},\roman{Alt}^p(Q,A))$, is the standard complex 
arising from a fine resolution of the sheaf of germs of basic $p$-forms on $M$.
Thus the cohomology 
$\roman E_1^{p,*}-\roman H^*(L_{\Cal F},\roman{Alt}^p(Q,A))$
is the cohomology of $M$ with values in the  sheaf of germs 
of basic $p$-forms on $M$. The corresponding spectral sequence (2.9.1)
comes down to the ordinary spectral sequence of a foliation, 
studied already in the literature, cf. 
\cite{\breinhar,\,\sarkhone,\,\sarkhtwo}; this spectral sequence is an
invariant of the foliation.
The cohomology $\roman E_2^{p,0}$ is sometimes called \lq\lq basic
cohomology\rq\rq, since it may be viewed as the cohomology of the 
\lq\lq space of leaves\rq\rq.

\smallskip
\noindent
(ii) Suppose that the foliation $\Cal F$ 
arises from a fiber bundle with fiber $F$, and 
write $\xi \colon P \to B$ for an associated principal bundle, the structure 
group being written as $G$. In this case, the spectral sequence (2.9.1) 
comes down to that of the fibration.
Furthermore, as a $C^{\infty}(B)$-module, the cohomology
$\roman H^*(L_{\Cal F},A)$ is the space of sections of the induced graded 
vector bundle $\zeta^* \colon P \times_G\roman H^*(F,\Bobb R) \to B$.
This vector bundle is flat and therefore inherits  a left 
$(C^{\infty}(B),\roman{Vect}(B))$-module structure, and $(E^{*,*}_1,d_1)$
coincides with the Rinehart complex 
$(\roman{Alt}^*_{C^{\infty}(B)}(\roman{Vect}(B), \Gamma(\zeta^*) ), d)$
which, in turn, is just the de Rham complex of $B$ with values in the 
flat vector bundle $\zeta^*$ and thus computes the cohomology
$E^{*,*}_2=\roman H^*(B,\zeta^*)$; equivalently, the flat connection on 
$\zeta^*$ turns $\roman H^*(F,\Bobb R)$ into a local system on $B$, and the 
de Rham complex of $B$ with values in the flat vector bundle $\zeta^*$
computes the cohomology $E^{*,*}_2=\roman H^*(B,\roman H^*(F,\Bobb R))$
of $B$ with coefficients in this local system.

\smallskip \noindent
(iii) Returning to (i) above, 
suppose in particular that the foliation is transversely complete
\cite \almolone. Then the closures of the leaves constitute a smooth fiber bundle
$M \to W$, the algebra $A^H$ is isomorphic to that of smooth functions
on $W$ in an obvious fashion, and the obvious map from
$Q^H$ to $\roman{Vect}(W) \cong \roman{Der}(A^H)$ which is part
of the Lie-Rinehart structure of $(A^H,Q^H)$ is surjective 
\cite\molinobo\ and hence
fits into an extension of $(\Bobb R, A^H)$-Lie algebras of the kind
$$
0 @>>> L' @>>> Q^H @>>> \roman{Vect}(W) @>>> 0.
\tag2.10.1
$$ 
Here $L'$ is the space of sections of a Lie algebra bundle on $W$, and the 
underlying extension of Lie algebroids on $W$ is referred to as the 
{\it Atiyah sequence\/} of the (transversely complete) foliation $\Cal F$ 
\cite\molinobo. Thus we see that the interpretation of $Q^H$ as the space 
of vector fields on the 
\lq\lq space of leaves\rq\rq\ 
requires, perhaps, some care, since
$L'$ will then consist of the \lq\lq vector fields on the 
\lq\lq space of leaves\rq\rq\ 
which act trivially on every function\rq\rq.
\smallskip
To get a concrete example, let $M = \roman{SU}(2) \times \roman{SU}(2)$,
and let $\Cal F$ be the foliation defined by a dense one-parameter subgroup
in a maximal torus $S^1 \times S^1$ in $\roman{SU}(2) \times \roman{SU}(2)$.
Then the space $W$ is $S^2 \times S^2$, and $L'$ is the space of sections of 
a real line bundle on $S^2 \times S^2$, necessarily trivial. One easily chooses 
a  vector bundle $\zeta$ on $\roman{SU}(2) \times \roman{SU}(2)$ which is 
complementary to $\tau_{\Cal F}$, and the Lie-Rinehart triple structure is 
defined on $(C^{\infty}(M), L_{\Cal F}, \Gamma(\zeta))$. In particular, the 
operation $\delta$ is non-zero. We note that the Chern-Weil construction 
in \cite\extensta\ yields a characteristic class in
$\roman H^2_{\roman{de Rham}}(S^2 \times S^2,\Bobb R)$ for the extension 
(2.10.1), and this class may be viewed as an {\it irrational Chern  class\/} 
\cite\extensta\ (Section 4). The non-triviality of this class entails that 
the differential $d_2$ of the spectral sequence (2.9.1) is non-trivial. 
We also note that, in view of a result of {\smc Almeida and Molino} 
\cite\almolone, the transitive Lie algebroid corresponding to (2.10.1) 
does not integrate to a principal bundle; in fact, {\smc Mackenzie's} 
integrability obstruction \cite\mackone\ is non-zero.

\smallskip \noindent
(iv) Under the circumstances of the Example 1.4.2, the cohomology 
$\roman H^*(H,\roman{Alt}^*(Q,A))$ is the Hodge cohomology of the smooth 
complex manifold $M$, i.~e. $\roman H^*(H,\roman{Alt}^p(Q,A))$ is the 
cohomology of $M$ with values in the sheaf of germs of holomorphic $p$-forms,
and the spectral sequence (2.9.1) is the {\it Hodge--de Rham spectral 
sequence\/}, sometimes 
referred to as the
{\it Fr\"olicher spectral sequence\/}
in the literature.

\smallskip \noindent
(v) Under the circumstances of Corollary 1.9.8,
so that $(A,Q,H)$ is a Lie-Rinehart triple with trivial
$(A,H)$-module structures on $A$ and $Q$,
the spectral sequence (2.9.1) is the ordinary spectral sequence
for the corresponding extension of Lie-Rinehart algebras.
If, furthermore, $A$ is the ground ring so that $Q$ and $H$ are ordinary
Lie algebras, this comes down to the Hochschild-Serre spectral sequence
of the Lie algebra extension.

\medskip\noindent {\bf 3. The additional structure on $Q$}
\smallskip\noindent
Let $(A,H,Q)$ be a Lie-Rinehart triple. Theorem 1.9  gives a possible answer
to Question 1.2 as well as to Question 1.1. What is missing is an intrinsic 
description of the structure induced on the constituent $(A,Q)$ which, in 
turn, should then in particular encapsulate the Lie-Rinehart triple structure 
on $(A,H,Q)$.We now proceed towards finding such an intrinsic description. 
To this end, we will introduce, on the constituent $Q$, certain operations 
similar to those introduced by Nomizu on the constituent $\fra q$ of a 
reductive decomposition $\fra g = \fra h \oplus \fra q$ of a Lie algebra 
\cite\nomiztwo; the operations in \cite\nomiztwo\ come from the curvature and 
torsion of an affine connection of the second kind. We note that the naive 
generalization to Lie-Rinehart algebras of the notion of reductive 
decomposition of a Lie algebra is not consistent with the Lie-Rinehart axioms. 
Given a Lie-Rinehart algebra $L$ and an $A$-module decomposition 
$L = H\oplus Q$ where $(A,H)$ inherits a Lie-Rinehart structure, since 
for $x \in H$, $\xi \in Q$, and $a \in A$, necessarily
$$
[x, a\xi] = a[x,\xi]-\xi(a) x,
$$
the defining property $[H,Q] \subset Q$ of a reductive decomposition cannot 
be satisfied {\it unless\/} the constituent $Q$ acts trivially on $A$.
\smallskip
Let $(A,H,Q)$ be an almost pre-Lie-Rinehart triple. We will now define 
triple@-, quadruple@-, and quintuple products of the kind
$$
\align
\{\cdot,\cdot;\cdot\}&\colon Q \otimes_R Q \otimes_R A @>>> A
\tag3.1
\\
\{\cdot;\cdot,\cdot;\cdot\}&\colon Q \otimes_R Q \otimes_R Q \otimes_R A @>>> A
\tag3.2
\\
\{\cdot;\cdot;\cdot,\cdot;\cdot\}
&\colon Q\otimes_R Q\otimes_R Q\otimes_R Q\otimes_R A 
@>>> A
\tag3.3
\\
\{\cdot,\cdot;\cdot\}&\colon Q \otimes_R Q \otimes_R Q @>>> Q
\tag3.4
\\
\{\cdot;\cdot,\cdot;\cdot\}&\colon Q \otimes_R Q \otimes_R Q \otimes_R Q @>>> Q
\tag3.5
\\
\{\cdot;\cdot;\cdot,\cdot;\cdot\}
&\colon Q\otimes_R Q\otimes_R Q\otimes_R Q\otimes_R Q 
@>>> Q .
\tag3.6
\endalign
$$
To this end, pick $\alpha,\beta,\gamma,\xi,\eta,\vartheta,\kappa \in Q$ and 
$a \in A$. For $1 \leq j \leq 6$, we will spell out an explicit description
of each of the operations $(3.j)$ and label it as $(3.j')$, as follows.
$$
\align
\{\xi,\eta;a\} &= (\delta(\xi,\eta))(a)
\tag3.1$'$
\\
\{\alpha;\xi,\eta;a\} &=(\alpha \cdot \delta(\xi,\eta))(a)
\tag3.2$'$
\\
\{\alpha;\beta;\xi,\eta;a\}
&=
(\alpha\cdot(\beta  \cdot \delta(\xi,\eta))(a)
\tag3.3$'$
\\
\{\xi,\eta;\vartheta\} &= (\delta(\xi,\eta))\cdot\vartheta
\tag3.4$'$
\\
\{\alpha;\xi,\eta;\kappa\}&=(\alpha \cdot \delta(\xi,\eta))\cdot \kappa
\tag3.5$'$
\\
\{\alpha;\beta;\xi,\eta;\gamma\}
&=
(\alpha\cdot(\beta  \cdot \delta(\xi,\eta))\cdot \gamma .
\tag3.6$'$
\endalign
$$

\proclaim{Proposition 3.7} Suppose that
$(A,H,Q)$ is a pre-Lie-Rinehart triple.
\newline\noindent
{\rm (i)} The operations
$\{\xi,\eta;\cdot\} \colon A \to A$,
$\{\alpha;\xi,\eta;\cdot\}\colon A \to A$,
$\{\alpha;\beta;\xi,\eta;\cdot\}\colon A \to A$
are derivations.
\newline\noindent
{\rm (ii)} The operations
$\{\xi,\eta;\cdot\} \colon A \to A$,
$\{\alpha;\xi,\eta;\cdot\}\colon A \to A$,
$\{\alpha;\beta;\xi,\eta;\cdot\}\colon A \to A$
are skew in the variables $\xi$ and $\eta$.
\newline\noindent
{\rm (iii)} The operations $\{\xi,\eta;\cdot\}$
{\rm (}on $A$ as well as on $Q${\rm )}
are $A$-linear in the variables $\xi$ and $\eta$, and the operations
{\rm (}on $A$ as well as on $Q${\rm )} $\{\alpha;\xi,\eta;\cdot\}$ and
$\{\alpha;\beta;\xi,\eta;\cdot\}$ are 
$A$-linear in the variable $\alpha$.
\newline\noindent
{\rm (iv)} The triple, quadruple, and quintuple products
$\{\xi,\eta;\vartheta\}$,
$\{\alpha;\xi,\eta;\kappa\}$,
$\{\alpha;\beta;\xi,\eta;\gamma\}$
are skew in the variables $\xi$ and $\eta$.
\newline\noindent
{\rm (v)}
Furthermore, these operations are related by the following identities.
$$
\alignat1
\{\xi,\eta;a\vartheta\} &= a\{\xi,\eta;\vartheta\} + \{\xi,\eta;a\} \vartheta
\\
\{\alpha;a\xi,\eta;b\} &
=\{\alpha;\xi,a\eta;b\} 
= a\{\alpha;\xi,\eta;b\} + \alpha(a) \{\xi,\eta;b\}
\\
\{\alpha;a\xi,\eta;\kappa\}&=
\{\alpha;\xi,a\eta;\kappa\}=
a \{\alpha;\xi,\eta;\kappa\}
+ \alpha(a) \{\xi,\eta;\kappa\}
\\
\{\alpha;\xi,\eta;a\kappa\}&=
a \{\alpha;\xi,\eta;\kappa\}
+\{\alpha;\xi,\eta;a\}\kappa
\\
\{\alpha;a\beta;\xi,\eta;b\}
&=
a\{\alpha;\beta;\xi,\eta;b\}
+
\alpha(a)\{\beta;\xi,\eta;b\}
\\
\{\alpha;\beta;a\xi,\eta;b\}
&=
\{\alpha;\beta;\xi,a\eta;b\}
\\
&=
a\{\alpha;\beta;\xi,\eta;b\}
+ \alpha(\beta(a))\{\xi,\eta;b\}
+ \beta(a) \{\alpha; \xi,\eta; b\}
+ \alpha(a) \{\beta; \xi,\eta; b\}
\\
\{\alpha;a\beta;\xi,\eta;\gamma\}
&=
a\{\alpha;\beta;\xi,\eta;\gamma\}
+
\alpha(a)\{\beta;\xi,\eta;\gamma\}
\\
\{\alpha;\beta;a\xi,\eta;\gamma\}
&=
\{\alpha;\beta;\xi,a\eta;\gamma\}
\\
&=
a\{\alpha;\beta;\xi,\eta;\gamma\}
+ \alpha(\beta(a))\{\xi,\eta;\gamma\}
+ \beta(a) \{\alpha; \xi,\eta; \gamma\}
+ \alpha(a) \{\beta; \xi,\eta; \gamma\}
\\
\{\alpha;\beta;\xi,\eta;a\gamma\}
&=
a\{\alpha;\beta;\xi,\eta;\gamma\}
+ \{\alpha;\beta;\xi,\eta;a\}\gamma
\endalignat
$$
\endproclaim

\demo{Proof} These assertions are immediate consequences of the
pre-Lie-Rinehart triple properties of $(A,H,Q)$.\qed
\enddemo

\proclaim{Proposition 3.8} 
Suppose that $(A,H,Q)$ is a pre-Lie-Rinehart triple, and let
$\alpha,\beta,\gamma,\zeta, \xi,\eta,\vartheta,\kappa \in Q$ and $a \in A$.
With the notation $x= \delta(\alpha,\beta)$ and $y=\delta(\gamma,\zeta)$, the 
compatibility properties {\rm (1.9.1)--(1.9.7)} take the following form.
$$
\alignat1
\xi \{\alpha,\beta; a\}
-
\{\alpha,\beta; \xi(a)\}
&=
\{\xi;\alpha,\beta; a\}
-
\{\alpha,\beta;\xi\}(a)
\tag3.8.1
\\
\{\alpha,\beta;[\xi,\eta]_Q\}
&=
[\{\alpha,\beta;\xi\},\eta]_Q
+
[\xi,\{\alpha,\beta;\eta\}]_Q
\\
&\quad
-\{\xi;\alpha,\beta; \eta\}
+\{\eta;\alpha,\beta; \xi\}
\tag3.8.2
\\
(\xi \cdot [\delta(\alpha,\beta),\delta(\gamma,\zeta)]_H)\cdot\kappa
&=
\{\xi; \alpha,\beta; \{\gamma,\zeta;\kappa\}\}
-
\{\gamma,\zeta,\{\xi; \alpha,\beta;\kappa\}\}
\\
& \quad
+ \{\alpha,\beta;\{\xi;\gamma,\zeta;\kappa\}\}
-\{\xi;\gamma,\zeta;\{\alpha,\beta;\kappa\}\}
\\
& \quad
-\{\{\alpha,\beta;\xi\};\gamma,\zeta;\kappa\}
+
\{\{\gamma,\zeta;\xi\};\alpha,\beta;\kappa\}
\tag3.8.3
\\
\xi (\eta(a))- \eta(\xi(a)) &= [\xi,\eta]_Q(a) + \{\xi,\eta;a\}
\tag3.8.4
\\
\{[\xi,\eta]_Q; \alpha,\beta; \gamma\}
&=\{\xi;\eta;\alpha,\beta;\gamma\}
-\{\eta;\xi;\alpha,\beta;\gamma\}
\\
&\quad
-\{\{\alpha,\beta;\xi\},\eta;\gamma\}
-\{\xi,\{\alpha,\beta; \eta\};\gamma\}
\\&
\quad
+
\{\alpha,\beta;\{\xi,\eta;\gamma\}\}
-
\{\xi,\eta;\{\alpha,\beta;\gamma\}\}
\tag3.8.5
\\
\sum_{(\xi,\eta,\vartheta)\ \roman{cyclic}}\left(
[[\xi,\eta]_Q,\vartheta]_Q \right.&+ 
\left .
\{\xi,\eta;\vartheta\} \right) = 0
\tag3.8.6
\\
\sum_{(\xi,\eta,\vartheta)\ \roman{cyclic}}
\{[\xi,\eta]_Q,\vartheta;\kappa\}
&=
\sum_{(\xi,\eta,\vartheta)\ \roman{cyclic}}
\{\xi;\eta,\vartheta;\kappa\}
\tag3.8.7
\endalignat
$$
Furthermore,
the compatibility property {\rm (1.5.12)}
takes the form
$$
[\delta(\alpha,\beta),\delta(\xi,\eta)]_H \cdot \xi 
= \{\alpha,\beta; \{\xi,\eta;\xi\}\}
- \{\xi,\eta; \{\alpha,\beta; \xi\}\}
\tag3.8.8
$$
\endproclaim

\demo{Proof} This is an immediate consequence
of Theorem 1.9.
We leave the details to the reader. \qed
\enddemo

We note that (3.8.5) is equivalent to
$$
\aligned
\{\alpha,\beta;\{\xi,\eta;\gamma\}\}
&-
\{\xi,\eta;\{\alpha,\beta;\gamma\}\}
\\
&
=
\{\{\alpha,\beta;\xi\},\eta;\gamma\}
+\{\xi,\{\alpha,\beta; \eta\};\gamma\}
\\
&\quad
+
\{[\xi,\eta]_Q; \alpha,\beta; \gamma\}
\\
&\quad
-\{\xi;\eta;\alpha,\beta;\gamma\}
+\{\eta;\xi;\alpha,\beta;\gamma\}
\endaligned
\tag3.8.5$'$
$$
Somewhat more explicitly, (3.8.7) reads
$$
\align
\{[\xi,\eta]_Q,\vartheta,\kappa\}
&+
\{[\eta,\vartheta]_Q,\xi,\kappa\}
+
\{[\vartheta,\xi]_Q,\eta,\kappa\}
\\
&=
(\xi \cdot \delta(\eta,\vartheta))\cdot \kappa
+
(\eta \cdot \delta(\vartheta,\xi))\cdot \kappa
+
(\vartheta \cdot \delta(\xi,\eta))\cdot \kappa
\endalign
$$
Moreover, with the notation $x = \delta(\xi,\eta)$,
(3.8.2) comes down to
$$
x\cdot [\vartheta,\kappa]_Q
=
[x\cdot\vartheta,\kappa]_Q
+
[\vartheta,x\cdot \kappa]_Q
-(\vartheta \cdot x)\cdot \kappa
+(\kappa \cdot x)\cdot \vartheta,
$$
which is just (1.9.2),
and
(3.8.5) reads
$$
\align
[x,\delta(\alpha,\beta)]_H 
&=\delta(x\cdot\alpha,\beta)
+\delta(\alpha,x\cdot \beta)
\\
&
\quad +
[\alpha,\beta]_Q \cdot x 
-\alpha\cdot(\beta  \cdot x)
+\beta\cdot (\alpha \cdot x),
\endalign
$$
which is (1.9.5).
\smallskip\noindent
{\smc Remark 3.9.} The description of the structure on $(A,Q)$ 
given in Propositions 3.7 and 3.8 is nearly intrinsic: Only the left-hand side
$(\xi \cdot [\delta(\alpha,\beta),\delta(\gamma,\zeta)]_H)(\kappa)$
of the equation (3.8.3) and the left-hand side
$[\delta(\alpha,\beta),\delta(\xi,\eta)]_H \cdot \xi $ of (3.8.8)
involve the Lie-Rinehart bracket $[\cdot,\cdot]_H$ on $H$ explicitly,
and this bracket is not covered by the structure on $(A,Q)$. The Lie-Rinehart 
structure of $(A,H)$ encapsulates a whole bunch of additional compatibility 
conditions which the triple@-, quadruple@-, quintuple products necessarily 
satisfy.
\smallskip\noindent
{\smc 3.10. Reconstruction of the Lie-Rinehart triple structure.\/}
Starting from $(A,Q)$, endowed with the pieces of structure
(1.5.1.$Q$) and (1.5.2.$Q$) which are supposed to satisfy
(1.5.6.$Q$) and (1.5.7.$Q$) and, furthermore,
with the triple@-, quadruple@-, quintuple products 
(3.1)--(3.6), to reconstruct an $(R,A)$-Lie algebra complement $H$ such that 
$E = H \oplus Q$ inherits an $(R,A)$-Lie algebra  structure which, in turn, 
then determines the given structure on $(A,Q)$, we might proceed as follows, 
where we pursue a reasoning similar to that in the proof of Theorem 18.1 in 
\cite\nomiztwo\ and that of Theorem 7.1 in \cite\kikkaone: 
Suppose that those compatibility properties spelled out in
(3.7) and (3.8) which are merely phrased in terms of $Q$ and, in particular,
do not involve the bracket $[\cdot,\cdot]_H$ on $H$ explicitly,
are satisfied. Given  $\xi,\eta \in Q$, define 
$\delta(\xi,\eta) \in\roman{End}_R(Q)$ by
$$
\delta(\xi,\eta)(\vartheta) = \{\xi,\eta;\vartheta\}
$$
and let $H \subseteq \roman{End}_R(Q)$ be the $A$-linear span of the
$\delta(\xi,\eta)$'s in $\roman{End}_R(Q)$ ($\xi,\eta \in Q$); notice that, 
by assumption, $Q$ comes with an $A$-module structure whence it makes sense to 
take the $A$-linear span of the $\delta(\xi,\eta)$'s in $\roman{End}_R(Q)$
($\xi,\eta \in Q$). The restriction of the evaluation pairing 
$\roman{End}_R(Q) \otimes_RQ \to Q$ to $H$ yields the pairing (1.5.3), 
to be written as the association
$$
(\delta(\xi,\eta), \vartheta) \longmapsto
\delta(\xi,\eta) \cdot \vartheta,\quad
\xi,\eta,\vartheta \in Q,
$$
and the requisite bilinear pairing (1.5.5) is just $\delta$, viewed as a 
function from $Q\otimes_AQ$ to $H$. Since the triple product (3.4) is 
$A$-bilinear, the pairing (1.5.3) will then satisfy (1.5.9), and $\delta$ is 
well defined on  $Q\otimes_AQ$. Next, define a pairing
$$
H \otimes_RA @>>> A,\quad (x,a) \mapsto x(a),
$$
by means of 
$$
\delta(\xi,\eta)(a) = \{\xi,\eta;a\},\quad  \xi,\eta \in Q,\ a \in A.
$$
This yields the requisite pairing (1.5.2.$H$). Since the triple product (3.1)
is $A$-bilinear, (1.5.6.$H$) will hold. Thereafter, define a pairing
$$
\cdot\,\colon Q \otimes_RH @>>> H
$$
by setting
$$
(\alpha \cdot \delta(\xi,\eta))(\kappa) = \{\alpha;\xi,\eta;\kappa\},
\quad  \alpha,\xi,\eta,\kappa \in Q.
$$
This yields the requisite pairing (1.5.4). Since the quadruple product is
$A$-linear in $\alpha$, (1.5.11) will hold. The compatibility properties
in (3.7) and (3.8) imply that the pairings (1.5.3) and (1.5.4) will satisfy
(1.5.8) and (1.5.10).
\smallskip
To complete the construction, we must require that {\sl the ordinary 
commutator bracket on $\roman{End}_R(Q)$ descend to a bracket 
$[\cdot,\cdot]_H$ on $H$ in such a way that $(A,H)$, with this bracket and the 
pairing {\rm (1.5.2.$H$)}} (which we reconstructed from the triple product 
(3.4)), {\sl be a Lie-Rinehart algebra in such a way that\/} (3.8.3) {\sl and\/}
(3.8.8) {\sl are satisfied\/}. The remaining compatibility properties in order 
for $(A,H,Q)$ to be a Lie-Rinehart triple will then be implied by the
structure isolated in (3.7) and (3.8).

\beginsection 4. Quasi-Lie-Rinehart algebras

Let $(A,H,Q)$ be a Lie-Rinehart triple. Thus $(A,H)$ is a Lie-Rinehart algebra
whence the Rinehart complex $\Cal A = (\roman{Alt}_A(H,A),d)$ inherits a 
differential graded $R$-algebra structure and $Q$ is, in particular,
an $(A,H)$-module whence the Rinehart complex $\Cal Q = (\roman{Alt}_A(H,Q),d)$ 
is a differential graded $\Cal A$-module in an obvious fashion. 
For the special case where $(A,H,Q)$ is a twilled Lie-Rinehart algebra
(i.~e. the operation $\delta \colon Q \otimes_AQ \to H$, cf. (1.5.5), is zero),
we have shown in \cite\twilled\ (3.2) that the pair $(\Cal A,\Cal Q)$ 
acquires a differential graded Lie-Rinehart structure and that the  
twilled Lie-Rinehart algebra compatibility conditions can be characterized in 
terms of this differential graded Lie-Rinehart structure. We will now show 
that, for a general Lie-Rinehart triple $(A,H,Q)$
(i.~e. with in general non-zero $\delta$), the pair $(\Cal A,\Cal Q)$
inherits a higher homotopy version of a differential graded Lie-Rinehart
algebra structure; abstracting from the structure which thus emerges, we 
isolate the notion of {\it quasi-Lie-Rinehart algebra\/}.
This structure  provides a complete solution of the problem of 
describing the 
structure on the constituent of a Lie-Rinehart triple written as $Q$
and hence yields a complete answer to Question 1.1. 
\smallskip
We begin by describing the requisite pieces of structure, independently of any 
given (pre-)Lie-Rinehart triple, in the following fashion:
Let $\Cal A$ be a 
graded commutative algebra
concentrated in non-negative degrees
($\Cal A^q = 0$ for $q<0$),
at this stage {\it not\/} a differential graded commutative algebra, and let 
$\Cal Q$ be a 
graded (left) $\Cal A$-module which we suppose to 
be an {\it induced\/} graded 
$\Cal A$-module of the kind $\Cal Q = \Cal A \otimes_AQ$ where $A=\Cal A^0$ and
where $Q$ is concentrated in degree zero; the notation
$(\Cal A,\Cal Q)$ will refer to this kind of structure throughout,
perhaps endowed with additional structure. 
A homogeneous $A$-multilinear function $\phi$
on $\Cal Q$ in $\ell$ variables with values in a graded 
$\Cal A$-module $\Cal M$ is said to be
$\Cal A$-{\it graded multilinear\/} if, for every
$\alpha_1, \dots,\alpha_{\ell} \in \Cal A$ and
every $\xi_1, \dots,\xi_{\ell} \in \Cal Q$,
$$
\align
\phi&(\xi_1,\dots,\xi_{j-1},\alpha_j \xi_j,\xi_{j+1},\dots, \xi_{\ell})
\\
&=
(-1)^{(|\phi|+|\xi_1|+ \dots + |\xi_{j-1}|)|\alpha_j|}
\alpha_j \phi(\xi_1,\dots,\xi_{j-1},\xi_j,\xi_{j+1},\dots, \xi_{\ell});
\endalign
$$
it is called 
{\it graded alternating\/} if, for every
$\xi_1, \dots,\xi_{\ell} \in \Cal Q$,
$$
\phi(\xi_1,\dots,\xi_j,\xi_{j+1},\dots, \xi_{\ell})
=-
(-1)^{|\xi_j||\xi_{j+1}|}
\phi(\xi_1,\dots,\xi_{j+1},\xi_j,\dots,\xi_{\ell}).
$$
A pairing is {\it graded skew-symmetric\/} provided it is
graded alternating as a graded bilinear function.
\smallskip
With these preparations out of the way
suppose  that, in addition, 
$(\Cal A,\Cal Q)$ carries
\newline\noindent
---  a graded skew-symmetric
$R$-bilinear pairing of degree zero
$$
[\,\cdot\,,\cdot\,]_{\Cal Q}
\colon
\Cal Q
\otimes_R
\Cal Q
@>>>
\Cal Q,
\tag4.1
$$
\newline\noindent
---  an $R$-bilinear pairing of degree zero
$$
\Cal Q \otimes_R \Cal A
@>>>
\Cal A,
\quad (\xi,\alpha) \mapsto \xi(\alpha),
\tag4.2
$$
\newline\noindent
---  an $A$-trilinear operation of degree $-1$
$$
\langle\cdot,\cdot\,;\cdot\rangle_{\Cal Q}
\colon
Q
\otimes_{A}
Q
\otimes_{A}
\Cal A
@>>>
\Cal A
\tag4.3.$Q$
$$
which is graded skew-symmetric in the first two variables
(i.~e. in the $Q$-variables).
\smallskip
We will say that the pair $(\Cal A,\Cal Q)$ constitutes a 
{\it pre-quasi-Lie-Rinehart algebra\/}
provided it satisfies (i) and (ii) below.
\newline\noindent
(i) The values of the adjoints $\Cal Q @>>> \roman{End}_R(\Cal A)$ and 
$Q \otimes_{A}Q @>>> \roman{End}_R(\Cal A)$ of (4.2) and 
(4.3.$Q$)
respectively, lie in $\roman{Der}_R(\Cal A)$ so that, in particular,
given 
$\xi,\eta \in Q$ and homogeneous $\alpha,\beta \in \Cal A$,
$$
\langle \xi,\eta;\beta\alpha \rangle_{\Cal Q}
=
\langle \xi,  \eta;\beta \rangle_{\Cal Q}\alpha
+(-1)^{|\beta|} \beta \langle \xi,  \eta;\alpha \rangle_{\Cal Q};
$$
\newline\noindent
(ii) the bracket (4.1), the operation (4.2), and the graded $\Cal A$-module 
structure on $\Cal Q$ satisfy the following graded Lie-Rinehart axioms (4.4) 
and (4.5):
$$
\alignat1
(a \xi) (b) &= a (\xi (b)),
\quad a, b \in \Cal A,\  \xi\in \Cal Q,
\tag4.4
\\
[\xi, a \eta]_{\Cal Q} &= \xi (a) \eta + a [\xi, \eta]_{\Cal Q},
\quad a \in \Cal A,\ \xi,\eta \in \Cal Q.
\tag4.5
\endalignat
$$
The graded Lie-Rinehart algebra axioms (4.4) and (4.5)
imply that (4.1) and (4.2)
are determined by their 
restrictions
$$
\gather
[\,\cdot\,,\cdot\,]_Q
\colon
Q
\otimes_R
Q
@>>>
Q
\tag4.1.$Q$
\\
Q \otimes_R \Cal A
@>>>
\Cal A,
\quad (\xi,\alpha) \mapsto \xi(\alpha)
\tag4.2.$Q$
\endgather
$$
Here the values of (4.1.$Q$) necessarily lie in $Q$ since 
$[\,\cdot\,,\cdot\,]_{\Cal Q}$ is supposed to be of degree zero; in particular, 
(4.1.$Q$) 
is skew-symmetric in the usual sense.
We note that,
when $\Cal A$ is concentrated in degree zero, the operation (4.3.$Q$)
is necessarily zero.
\smallskip
Given a  pre-quasi-Lie-Rinehart algebra $(\Cal A,\Cal Q)$, consider the 
bigraded algebra
$$
\roman{Alt}_A(Q,\Cal A)
\cong
\roman{Alt}_{\Cal A}(\Cal Q,\Cal A),
\tag4.6
$$
of $\Cal A$-valued $A$-multilinear alternating functions on $Q$
and define the operators
$$
d_1\colon
\roman{Alt}_A^p(Q,\Cal A^q)
@>>>
\roman{Alt}_A^{p+1}(Q,\Cal A^q)
\quad (p,q \geq 0)
\tag4.7.1
$$
and
$$
d_2\colon
\roman{Alt}_A^p(Q,\Cal A^q)
@>>>
\roman{Alt}_A^{p+2}(Q,\Cal A^{q-1})
\quad (p,q \geq 0)
\tag4.8.1
$$
by 
$$
\aligned
(-1)^{|f|+1}(d_1f)&(\xi_1,\dots,\xi_{p+1})
= 
\sum_{j=1}^{p+1} 
(-1)^{j-1}
\xi_j(f(\xi_1,\dots\widehat{\xi_j}\dots,\xi_{p+1}))
\\
&\phantom{=}+\quad
\sum_{1 \leq j<k \leq p+1} (-1)^{j+k} f(\lbrack \xi_j,\xi_k \rbrack_Q,
\xi_1,\dots\widehat{\xi_j}\dots\widehat{\xi_k}\dots,\xi_{p+1})
\endaligned
\tag4.7.2
$$
(the graded CCE formula)
$$
\aligned
(-1)^{|f|+1}&(d_2f)(\xi_1,\dots,\xi_{p+2})
\\
&= \phantom{=+}\quad
(-1)^p\sum_{1 \leq j<k \leq p+2}  
(-1)^{j+k}
\langle\xi_j,\xi_k;
f(\xi_1, \dots\widehat{\xi_j}\dots\widehat{\xi_k}\dots,\xi_{p+2})\rangle_Q 
\endaligned
\tag4.8.2
$$
where $\xi_1,\dots,\xi_{p+2} \in Q$. The graded Lie-Rinehart axioms (4.4) and 
(4.5) imply that the operator $d_1$ is well defined on 
$\roman{Alt}_A(Q,\Cal A)$
as an $R$-linear (beware, {\it not\/} $A$-linear)
operator.
The usual argument shows that
$d_1$ is a derivation on the
bigraded $A$-algebra 
$\roman{Alt}_A(Q,\Cal A)$.
Since the operation
$\langle\cdot,\cdot\,;\cdot\rangle_{\Cal Q}$,
cf. (4.3.$Q$), is  $A$-trilinear, 
the operator $d_2$ is well defined 
on $\Cal A$-valued $A$-multilinear functions on $Q$. 
Since
(4.3.$Q$) is skew-symmetric in the first two variables, the operator $d_2$ 
automatically has square zero, i.~e. is a differential.

\proclaim{Lemma 4.8.3}
The operator $d_2$ is an $A$-linear derivation 
on the bigraded $A$-algebra $\roman{Alt}_A(Q,\Cal A)$.
\endproclaim

\demo{Proof}
Since, as a graded $\Cal A$-module,
$\Cal Q$  is an induced graded
$\Cal A$-module,
the bigraded algebra $\roman{Alt}_A(Q,\Cal A)$
may be written as the bigraded tensor product
$\roman{Alt}_A(Q,\Cal A) \cong \roman{Alt}_A(Q,A) \otimes \Cal A$,
and it suffices to consider forms which
may be written as
$\beta \alpha$ where 
$\beta \in \roman{Alt}_A(Q,A)$ 
and $\alpha \in \Cal A$;
the formula (4.8.2) yields
$$
d_2(\beta) = 0,\quad
d_2(\beta \alpha) = (-1)^{|\beta|} \beta d_2(\alpha)
$$
and, 
since for $\xi,\eta \in Q$, the operation
$\langle\xi,\eta;\cdot \rangle_Q$ is a derivation of $\Cal A$,
we conclude that the operator $d_2$ is an $R$-linear derivation 
on $\roman{Alt}_A(Q,\Cal A)$.
Furthermore, since
for $a \in A = \Cal A^0$,
for degree reasons, 
$d_2(a)$ is necessarily zero
the operator $d_2$ is plainly well defined 
on $\Cal A$-valued $A$-multilinear functions on $Q$ 
and in fact an $A$-linear derivation 
on $\roman{Alt}_A(Q,\Cal A)$ as asserted.  \qed
\enddemo

\noindent
{\smc Remark 4.8.4.\/} 
On the formal level, the notion of quasi-Lie-Rinehart algebra isolated above
is somewhat unsatisfactory since the definition
involves the structure of $\Cal Q$ as an {\it induced\/} 
$\Cal A$-module. The operator $d_1$ may be written out as
an operator on the bigraded $A$-module $\roman{Alt}_{\Cal A}(\Cal Q,\Cal A)$
of $\Cal A$-graded multilinear alternating forms
on $\Cal Q$ directly in terms of the operations (4.1) and (4.2),
that is, in terms of the arguments of these operations, 
without explicit reference to the induced $\Cal A$-module structure. Indeed, 
given an $n$-tuple $\eta = (\eta_1,\dots,\eta_n)$
of homogeneous elements of $\Cal Q$, write
$|\eta| = |\eta_1| + \ldots +|\eta_n|$ and 
$|\eta|^{(j)} = |\eta_1| + \ldots +|\eta_j|$,
for $1 \leq j \leq n$, and define the operators
$$
d_{(\cdot,\cdot)}
\colon
\roman{Alt}_R^p(\Cal Q,\Cal A^q)
\to
\roman{Alt}_R^{p+1}(\Cal Q,\Cal A^q),
\quad
d_{[\cdot,\cdot]}
\colon
\roman{Alt}_R^p(\Cal Q,\Cal A^q)
\to
\roman{Alt}_R^{p+1}(\Cal Q,\Cal A^q),
$$
by means of
$$
\alignat1
{}
&(-1)^{|f|+1+ |\eta|}
(d_{(\cdot,\cdot)}(f))(\eta_1,\ldots,\eta_{p+1})
\\
&=
\sum_{j=1}^{p+1}
(-1)^{j-1+(|\eta|^{(j-1)}+|f|)|\eta_j|} 
\eta_j f(\eta_1,\ldots \widehat{\eta_j} \ldots,\eta_{p+1})
\\
{}
&(-1)^{|f|+1+|\eta|}(d_{[\cdot,\cdot]}(f))(\eta_1,\ldots,\eta_{p+1})
\\
&=
\sum_{1\leq j <k \leq p+1}
(-1)^{j+k+|\eta|^{(j-1)}|\eta_j|
+(|\eta|^{(k-1)} -|\eta_j|)|\eta_k|}
f([\eta_j,\eta_k],
\eta_1,\ldots\widehat \eta_j \ldots \widehat{\eta_k} \ldots ,\eta_{p+1})
\endalignat
$$
where $\eta_1,\dots,\eta_{p+1}$ are homogeneous elements of $\Cal Q$.
Then the sum $d_{(\cdot,\cdot)}+d_{[\cdot,\cdot]}$ descends to an operator 
on $\roman{Alt}_{\Cal A}(\Cal Q,\Cal A)$ which, in turn, coincides with $d_1$.
In this fashion, $d_1$ appears as being given by
the CCE formula (2.2.8) with respect to (4.1) and (4.2).
We were so far unable to give a similar description of the operator
$d_2$, though, in terms of a suitable extension of (4.3.$Q$) to an operation
of the kind 
$\Cal Q
\otimes_{A}
\Cal Q
\otimes_{A}
\Cal A
@>>>
\Cal A
$.

\smallskip\noindent
{\smc 4.9. Definition\/.} Let $(\Cal A,\Cal Q)$ be a pre-quasi-Lie-Rinehart 
algebra so that, in particular,
$\Cal A$ is a differential graded commutative algebra and
$\Cal Q$ a differential graded $\Cal A$-module.
Consider the bigraded $A$-algebra
$$
\roman{Alt}_A(Q,\Cal A)
\cong
\roman{Alt}_{\Cal A}(\Cal Q,\Cal A)
\subseteq
\roman{Mult}_R(\Cal Q,\Cal A),
$$
cf. (4.6) above, where $\roman{Mult}_R(\Cal Q,\Cal A)$ refers to the
bigraded algebra of $\Cal A$-valued $R$-multilinear forms on $\Cal Q$. The 
differentials on $\Cal Q$ and $\Cal A$ (both written as $d$, with an abuse 
of notation,) induce a differential $D$ on $\roman{Mult}_R(\Cal Q,\Cal A)$
in the usual way, that is, given an $R$-multilinear  $\Cal A$-valued form $f$ 
on $\Cal Q$, 
$$
Df = df + (-1)^{|f|+1} fd
$$
where, with a further abuse of notation, the \lq\lq $d$\rq\rq\ 
in the constituent $fd$ signifies the
induced operator on any of the tensor powers
$\Cal Q^{\otimes_R \ell}$ ($\ell \geq 1$).
We will say that the pre-quasi-Lie-Rinehart algebra $(\Cal A,\Cal Q)$
is a {\it quasi-Lie-Rinehart algebra} provided it satisfies the requirements
(4.9.1)--(4.9.6) below where $d_1$ and $d_2$ are the operators
(4.7.1) and (4.8.1), respectively.
\newline\noindent
(4.9.1)
The differential $D$ descends to an operator 
on $\roman{Alt}_{\Cal A}(\Cal Q,\Cal A)$,
necessarily a differential,
which we then write as $d_0$.
\newline\noindent
(4.9.2)
The differential on $\Cal Q$ is a derivation for the bracket (4.1).
\newline\noindent
(4.9.3)
The pairing (4.2) is compatible with the differentials on $\Cal A$ and $\Cal Q$.
\newline\noindent
(4.9.4)
For every $\xi,\eta \in Q$ and 
$\alpha \in \Cal A$,
$$
\xi(\eta(\alpha)) - \eta(\xi(\alpha)) -[\xi,\eta]_Q(\alpha) 
= \left((d_0 d_2+d_2 d_0)(\alpha)\right)(\xi,\eta) 
$$
\newline\noindent
(4.9.5)
For every $\xi,\eta, \vartheta\in Q$ and 
$\alpha \in \roman{Alt}_A^1(Q,\Cal A^0) = \roman{Hom}_A(Q,A)$,
$$
\sum_{(\xi,\eta,\vartheta)\ \roman{cyclic}}
\alpha([[\xi,\eta]_Q,\vartheta]_Q)
=
(d_2 d_0 \alpha)(\xi,\eta,\vartheta)
+
\sum_{(\xi,\eta,\vartheta)\ \roman{cyclic}}
(d_2 d_0 \alpha(\vartheta))(\xi,\eta)
$$
\newline\noindent
(4.9.6) The operators $d_1$ and $d_2$ satisfy the commutation relation
$$
d_1d_2 + d_2 d_1 = 0.
$$
In (4.9.6), it suffices to require 
the vanishing of the operator $d_1d_2 + d_2 d_1$
on $\roman{Alt}_A^0(Q,\Cal A^1)$.
We leave it to the reader to spell out a description
of this requirement directly in terms of the structure
(4.1)--(4.3); this description would be less concise than
the requirement given as (4.9.6).

\proclaim{Theorem 4.10}
Let $(\Cal A,\Cal Q)$ be a pre-quasi-Lie-Rinehart 
algebra, consider the bigraded $A$-algebra
$$
\roman{Alt}_A(Q,\Cal A)
\cong
\roman{Alt}_{\Cal A}(\Cal Q,\Cal A)
\subseteq
\roman{Mult}_R(\Cal Q,\Cal A),
$$
suppose that the operator $D$ 
on $\roman{Mult}_R(\Cal Q,\Cal A)$ descends to an operator $d_0$
on $\roman{Alt}_A(Q,\Cal A)$,
and let $d_1$ and $d_2$ be the operators on $\roman{Alt}_A(Q,\Cal A)$
given by {\rm (4.7.1)} and {\rm (4.8.1)}, respectively.
Then $(\Cal A,\Cal Q)$ is a quasi-Lie-Rinehart algebra
if and only if $(\roman{Alt}_{\Cal A}(\Cal Q,\Cal A),d_0, d_1, d_2)$ is a 
multialgebra.
\endproclaim

\demo{Proof}
(i) The identity
$0=d_0 d_1 + d_1 d_0$ on 
$\roman{Alt}_A^1(Q,\Cal A^0)$
is equivalent to (4.9.2), that is,
to the differential on $\Cal Q$ being a derivation for the bracket
$[\,\cdot\,,\cdot\,]_{\Cal Q}$, cf (4.1). See also (2.8.5(iii)).

\noindent
(ii) The identity
$0=d_0 d_1 + d_1 d_0$ on
$\roman{Alt}_A^0(Q,\Cal A^*)$
is equivalent to (4.9.3), that is,
to the differentials on $\Cal A$ and $\Cal Q$ being
compatible with the pairing (4.2).

\noindent
(iii) The identity
$0=d_0 d_2 + d_1 d_1+ d_2 d_0$ on
$\roman{Alt}_A^0(Q,\Cal A^0)$
is equivalent to the special case
of (4.9.4) where $\alpha \in A = \Cal A^0$. Cf. (2.8.5(iv)). 

\noindent
(iv) Once (4.9.4) holds, the identity
$0=d_0 d_2 + d_1 d_1+ d_2 d_0$ on
$\roman{Alt}_A^1(Q,\Cal A^0)$
is equivalent to (4.9.5). Cf. (2.8.5(v)).

\noindent
(v) The identity
$0=d_0 d_2 + d_1 d_1+ d_2 d_0$
on
$\roman{Alt}_A^0(Q,\Cal A^1)$
is equivalent to the special case
of (4.9.4) where $\alpha \in \Cal A^1$. 
Cf. (2.8.5(vi)). \quad \qed 
\enddemo

Under the circumstances of Theorem 4.10, we will refer to the multialgebra 
$$
(\roman{Alt}_{\Cal A}(\Cal Q,\Cal A),d_0,d_1,d_2)
$$
as the {\it Maurer-Cartan algebra\/} for the 
quasi-Lie-Rinehart algebra structure on $(\Cal A,\Cal Q)$.

\smallskip
\noindent
{\smc 4.11. Relationhip with almost pre-Lie-Rinehart triples.}
Our goal is to show how a Lie-Rinehart triple determines a
quasi-Lie-Rinehart algebra.
Here we explain the first step, that is,
how a structure of the kind (4.1.$Q$)--(4.3.$Q$) 
that underlies a pre-quasi-Lie-Rinehart algebra arises: 
Let $(A,Q,H)$ be an almost pre-Lie-Rinehart triple, and let 
$\Cal A = \roman{Alt}_A(H,A)$ and $\Cal Q = \roman{Alt}_A(H,Q)$. Then
$\Cal A = \roman{Alt}_A(H,A)$ is a graded commutative algebra
(beware, not necessarily a {\it differential\/} graded commutative algebra)
and $\Cal Q = \roman{Alt}_A(H,Q)$ is a graded $\Cal A$-module 
(not necessarily a {\it differential\/} graded module).
The pairings (1.5.2.$Q$) and (1.5.4) induce a pairing
$Q \otimes_R \Cal A \to\Cal A$ of the kind
(4.2.$Q$) by means of the association
$$
\xi \otimes \alpha \mapsto \xi(\alpha),\quad
\xi \in Q,\ \alpha \in \Cal A = \roman{Alt}_A(H,A)
\tag4.11.1
$$
where
$$
(\xi(\alpha))(x_1,\dots, x_n) =
\xi(\alpha(x_1,\dots, x_n))
-\sum_{j=1}^n \alpha (x_1,\dots, \xi \cdot x_j, \dots x_n).
\tag4.11.2
$$
The corresponding
induced pairing of the kind (4.2) has the form
$$
\Cal Q \otimes_R \Cal A
@>>>
\Cal A,
\quad (\xi,\alpha) \mapsto \xi(\alpha),\ \xi \in \Cal Q,\ 
\alpha \in \Cal A.
\tag4.11.3
$$
Furthermore, the bracket 
$[\,\cdot\,,\cdot\,]_Q$ is exactly of the kind (4.1.$Q$).
It extends to a graded skew-symmetric bracket
$$
[\,\cdot\,,\cdot\,]_{\Cal Q}
\colon
\Cal Q
\otimes_R
\Cal Q
@>>>
\Cal Q
\tag4.11.4
$$
of the kind (4.1). To get an explicit formula for this bracket we suppose, for 
simplicity, that the canonical map from $\Cal A \otimes_A Q$ to
$\Cal Q = \roman{Alt}_A(H,Q)$ is an isomorphism of graded $A$-modules so that
$\Cal Q$ is indeed an induced graded $\Cal A$-module of the  kind considered 
above. This will be the case, for example, when $H$ is finitely generated
and projective as an $A$-module or when $Q$ is projective as an $A$-module.
Under these circumstances, given homogeneous elements $\alpha,\beta \in \Cal A$
and $\xi,\eta \in Q$, the value 
$[\alpha\otimes \xi,\beta \otimes \eta ]_{\Cal Q}$ of the bracket (4.11.4)
is given  by
$$
[\alpha\otimes \xi,\beta \otimes \eta ]_{\Cal Q}
=
(\alpha\xi(\beta)) \otimes \eta
-(\beta\eta(\alpha))\otimes\xi +
(\alpha\beta) \otimes[\xi,\eta ]_Q .
\tag4.11.5
$$
Furthermore,
setting
$$
\langle\xi,\eta;\alpha\rangle_Q
=
i_{\delta(\xi,\eta)}\alpha,\quad
\alpha \in \Cal A,\, \xi,\eta \in Q,
\tag4.11.6
$$
where, for $x \in H$, $i_x$ refers to the operation of
contraction, that is,
$$
\langle\xi,\eta;\alpha\rangle_Q(x_1,\dots,x_{q-1})
=
\alpha(\delta(\xi,\eta),x_1,\dots,x_{q-1}),\quad
x_1,\dots,x_{q-1} \in H,
\tag4.11.7
$$
we obtain a pairing of the kind (4.3.$Q$).
Thus, summing up, we conclude that, on $(\Cal A,\Cal Q)$,
the operations (4.11.1), (4.11.4), and (4.11.6) which, in turn, come from
the almost pre-Lie-Rinehart triple structure on $(A,Q,H)$
determine a structure of the kind (4.1.$Q$)--(4.3.$Q$) 
which underlies that of a
pre-quasi-Lie-Rinehart algebra.
Indeed, the structure on $(\Cal A,\Cal Q)$ 
given by (4.11.1), (4.11.4), and (4.11.6)
is essentially a rewrite of the 
almost pre-Lie-Rinehart triple structure on $(A,Q,H)$; the two structures are 
equivalent when $H$ is finitely 
generated  projective as an $A$-module and when $Q$ has property P.
At this stage we do not make any claim as to whether or not
the structure given by (4.11.1), (4.11.4), and (4.11.6)
turns  $(\Cal A,\Cal Q)$ into a pre-quasi-Lie-Rinehart algebra.

\smallskip
\noindent
{\smc 4.12. Lie-Rinehart triples and quasi-Lie-Rinehart algebras.}
Suppose now that $(A,Q,H)$ is a pre-Lie-Rinehart triple; with reference to the 
Lie-Rinehart structure on $(A,H)$ and the left $(A,H)$-module structure on $Q$, 
the Lie-Rinehart differentials then turn 
$\Cal A= \roman{Alt}_A(H,A)$
into a {\it differential\/} graded commutative algebra and 
$\Cal Q= \roman{Alt}_A(H,Q)$
into a {\it differential\/} graded (left) $\Cal A$-module;
cf. (1.5.13). Furthermore,
the bigraded algebra $\roman{Alt}_A(Q,\Cal A)$ of alternating $A$-multilinear 
$\Cal A$-valued forms on $Q$ may be rewritten in the form 
$\roman{Alt}_A(H,\roman{Alt}_A(Q,A))$; equivalently, the algebra 
$\roman{Alt}_A(Q,\Cal A)$ may be viewed as the bigraded algebra 
$\roman{Alt}_{\Cal A}(\Cal Q,\Cal A)$ of alternating $\Cal A$-multilinear 
$\Cal A$-valued forms on $\Cal Q$. 
We write the resulting operator $d_0$, cf. (2.4.5) and (2.5.1),
as
$$
d_0 \colon \roman{Alt}^p_A(Q,\Cal A^q)
@>>>
\roman{Alt}^p_A(Q,\Cal A^{q+1}) \quad (p,q \geq 0).
\tag4.12.1
$$
\smallskip
Consider the operators $d_1$ and $d_2$  on $\roman{Alt}_A(Q,\Cal A)$
given as (4.7.1) and (4.8.1) above, respectively.
These operators  now come down to the operators
(2.4.6) and (2.4.7), respectively.
By Theorem 2.7, when $(A,Q,H)$ is a genuine Lie-Rinehart triple, 
$$
(\roman{Alt}_A(Q,\Cal A),d_0,d_1,d_2)
=
(\roman{Alt}_A(Q,\roman{Alt}_A(H,A)),d_0,d_1,d_2)
\tag4.12.2
$$
is a Maurer-Cartan algebra, that is, 
$d = d_0 + d_1 + d_2$ turns $\roman{Alt}_{\Cal A}(\Cal Q,\Cal A)$ into a 
differential graded algebra.
Furthermore, still by Theorem 2.7,
under the assumption that $H$ and $Q$ both have property P,
the converse holds, i.~e. when (4.12.2) is a Maurer-Cartan algebra,
$(A,Q,H)$ is a genuine Lie-Rinehart triple. In view of Theorem 4.10
we conclude the following.

\proclaim{Theorem 4.13}
Let $(A,H,Q)$ be a pre-Lie-Rinehart triple and suppose that
both $H$ and $Q$ have property {\rm P}, (e.~g. $H$ and $Q$ are
both projective as $A$-modules).
Then $(A,H,Q)$ is a genuine Lie-Rinehart triple
if and only if
$$
(\Cal A,\Cal Q)= (\roman{Alt}_A(H,A),\roman{Alt}_A(H,Q)),
$$ 
endowed with the pairing {\rm (4.11.1)}, the bracket
$[\,\cdot\,,\cdot\,]_{\Cal Q}$, cf. {\rm (4.11.4)}, and the 
operation 
$\langle\xi,\eta;\alpha\rangle_Q$, cf. {\rm (4.11.6)}, is a quasi-Lie-Rinehart 
algebra. \qed
\endproclaim

The proof of the following is straightforward and left to the reader:

\proclaim{Proposition 4.14}
The homology $(\roman H^*(\Cal A),\roman H^*(\Cal Q))$ of a quasi-Lie-Rinehart 
algebra $(\Cal A,\Cal Q)$ inherits a graded Lie-Rinehart algebra structure. \qed
\endproclaim

Given a Lie-Rinehart triple $(A,H,Q)$,
the graded Lie-Rinehart algebra
\linebreak
$(\roman H^*(\Cal A),\roman H^*(\Cal Q))$
of the corresponding quasi-Lie-Rinehart algebra
$$
(\Cal A,\Cal Q)= (\roman{Alt}_A(H,A),\roman{Alt}_A(H,Q))
$$
contains more information than the Lie-Rinehart algebra
$(A^H,Q^H) =(\roman H^0(\Cal A),\roman H^0(\Cal Q))$
spelled out in Corollary 1.11.
\smallskip
\noindent
{\smc Illustration 4.15.} Let $(M,\Cal F)$ be a foliated manifold,
maintain the notation established earlier
in (1.4.1), (1.12), and (2.11),
let $(A,H,Q) = (C^{\infty}(M), L_{\Cal F}, Q)$,
the corresponding Lie-Rinehart triple, and consider
the resulting quasi-Lie-Rinehart algebra
$(\Cal A,\Cal Q)= (\roman{Alt}_A(H,A),\roman{Alt}_A(H,Q))$.
We may view
$\Cal A$ as the {\it algebra of generalized functions\/}
and 
$\Cal Q$
as the {\it generalized Lie algebra of vector fields\/}
for the foliation.
Thus $\Cal A$ is the standard complex arising from a fine resolution
of the sheaf of germs of functions on $M$ which are constant on the leaves.
Likewise, the constituent $Q^H$ of the Lie-Rinehart algebra $(A^H,Q^H)$ 
(discussed earlier), cf. (1.12) and (2.10) (iii),
amounts to the space of global sections
of the sheaf $\Cal V_Q$ of germs of vector fields on $M$ which are 
horizontal (with respect to the decomposition 
$\Gamma(\tau_M) = L_{\Cal F} \oplus Q$) and constant on the leaves, and
$\Cal Q$ is the standard complex arising from a fine resolution
of this sheaf.
Thus
$\roman H^*(\Cal A)$
is the cohomology
of $M$ with values in the sheaf of germs of functions
which are constant on the leaves,
and
$\roman H^*(\Cal Q)$
is the cohomology
of $M$ with values in the sheaf
$\Cal V_Q$.
\smallskip
Under the circumstances of (2.10(ii)),
so that the foliation $\Cal F$ comes from a fiber bundle
and the space of leaves coincides with the base $B$ of the
corresponding fibration, from the graded commutative $\Bobb R$-algebra 
structure of $\roman H^*(F,\Bobb R)$,
the space $\Gamma(\zeta^*)$ of sections of 
the induced graded vector bundle 
$$
\zeta^*\colon P \times_G \roman H^*(F,\Bobb R) \to B
$$
inherits a graded 
$C^{\infty}(B)$-algebra structure
and, as a graded $C^{\infty}(B)$-algebra,
$\roman H^*(\Cal A)$ coincides with
the graded commutative algebra
$\Gamma(\zeta^*)$ of sections of $\zeta^*$;
in particular,
$\roman H^0(\Cal A) =C^{\infty}(B)$.
Furthermore, 
$\roman H^0(\Cal Q)$
is the $(\Bobb R,C^{\infty}(B))$-Lie algebra $\roman{Vect}(B)$ of smooth vector
fields on the base $B$ and, as a graded
$(\Bobb R,\roman H^*(\Cal A))$-Lie algebra,
$\roman H^*(\Cal Q)$
is the graded crossed product
$$
\roman H^*(\Cal Q) = 
\roman H^*(\Cal A) \otimes_{C^{\infty}(B)} \roman{Vect}(B)
\tag4.15.1
$$
(cf. \cite\twilled\ for the notion of graded crossed product Lie-Rinehart
algebra).
\smallskip
Under the circumstances of (2.10(i)), when the foliation 
does {\it not\/} come from a fiber bundle,
the structure of the graded Lie-Rinehart algebra
$(\roman H^*(\Cal A),\roman H^*(\Cal Q))$
will in general be more complicated than that for the case when the foliation
comes from a fiber bundle. The significance of this more complicated
structure has been commented on already in the introduction.

\smallskip\noindent
{\smc Remark 4.16.} We are indebted to P. Michor for having pointed out to us
a possible connection of the notion of quasi-Lie-Rinehart bracket with that of
Fr\"olicher-Nijenhuis bracket \cite\fronione,\, \cite\nijenhui.
Given a smooth manifold $M$,
the Fr\"olicher-Nijenhuis bracket is defined on the graded vector space
of forms on $M$ with values in the tangent bundle $\tau_M$ of $M$ and endows 
this graded vector space with a graded Lie algebra structure
which in degree zero amounts to the ordinary Lie bracket of
vector fields on $M$.
Given a Lie-Rinehart algebra $(A,L)$, an obvious generalization 
of the Fr\"olicher-Nijenhuis bracket  
endows the graded $A$-module $\roman{Alt}_A(L,L)$
with a graded $R$-Lie algebra structure.
Given a Lie-Rinehart triple $(A,H,Q)$, with correponding
Lie-Rinehart algebra $(A,L)$ where $L = H \oplus Q$,
the induced quasi-Lie-Rinehart bracket (4.11.4)
is defined on $\roman{Alt}_A(H,Q)$, and the obvious question
arises how this quasi-Lie-Rinehart bracket is related with the
Fr\"olicher-Nijenhuis bracket on 
$\roman{Alt}_A(L,L)$.

\medskip\noindent{\bf 5. Quasi-Gerstenhaber algebras}
\smallskip\noindent
The notion of Gerstenhaber algebra has recently been isolated in the literature
but implicitly occurs already in Gerstenhaber's paper \cite\gersthtw;
see \cite\bv\ for details and more references. In this section we will introduce
a notion of 
quasi-Gerstenhaber algebra 
which
generalizes that of
{\it strict differential 
bigraded\/} Gerstenhaber algebra isolated in \cite{\twilled,\,\banach}
(where the attribute \lq\lq strict\rq\rq\ refers to the requirement that
the differential be a derivation for the Gerstenhaber bracket). 
The generalization
consists 
in admitting a bracket which does {\sl not necessarily satisfy the graded 
Jacobi identity and incorporating an additional piece of structure which 
measures the deviation from the graded Jacobi identity\/}. 
\smallskip
For intelligibility, we recall the notion of graded Lie algebra, tailored to
our purposes. As before, $R$ denotes a commutative ring with 1. A graded 
$R$-module $\fra g$, endowed with a graded skew-symmetric degree zero bracket
${
[\,\cdot\,,\cdot\,]\colon \fra g \otimes \fra g @>>> \fra g,
}$
is called a {\it graded Lie algebra\/} provided the bracket satisfies
the {\it graded Jacobi identity\/}
$$
\sum_{(a,b,c)\ \roman{cyclic}}(-1)^{|a||c|}[a,[b,c]] = 0,
$$
for every triple $(a,b,c)$ of homogeneous elements of $\fra g$.
\smallskip
Given a graded commutative algebra $\Cal A$, 
an (ordered) $m$-tuple $\bold a =(a_1,\dots,a_m)$ of homogeneous
elements thereof, and a permutation $\sigma$ of $m$ objects, we denote by 
$\varepsilon(\bold a,\sigma)$ the sign
defined by
$$
a_1\cdot \ldots \cdot a_m = \varepsilon(\bold a,\sigma)
a_{\sigma 1}\cdot \ldots \cdot a_{\sigma m} 
$$
according to the {\it Eilenberg-Koszul\/} convention.
\smallskip
We will consider bigraded $R$-algebras.
Such a bigraded algebra is said to be {\it bigraded commutative\/} provided 
it is commutative in the bigraded sense, that is, graded commutative with 
respect to the total degree. Given such a bigraded commutative 
algebra $\Cal G$, for bookkeeping purposes, we will write its homogeneous 
components in the form $\Cal G^q_p$, the superscript being viewed as a 
cohomology degree and the subscript as a homology degree; the {\it total 
degree\/} 
$|\alpha|$ of an  element $\alpha$ of $\Cal G^q_p$ is, then, $|\alpha|=p-q$.
\smallskip
We will explore differential operators in the bigraded context.
We recall the requisite notions from \cite\koszulon\ (Section 1),
cf. also \cite\akmantwo.
Let $\Cal G$ be a bigraded commutative $R$-algebra with 1, and
let $r\geq 1$.
A (homogeneous) {\it differential operator\/} on $\Cal G$ {\it of order\/} 
$\leq r$ is a homogeneous $R$-endomorphism $D$ of $\Cal G$ such that a certain
$\Cal G$-valued $(r+1)$-form $\Phi^{r+1}_D$ on $\Cal G$
(the definiton of which for general $r$ we do not reproduce here)
vanishes. For our purposes, it suffices to recall explicit descriptions
of these forms in low degrees. Thus, given
the homogeneous $R$-endomorphism $D$ of $\Cal G$, for
homogeneous $\xi,\eta,\vartheta$,
$$
\align
\Phi^1_D(\xi) &= D(\xi) - D(1) \xi
\\
\Phi^2_D(\xi,\eta) &= D(\xi \eta) - D(\xi) \eta 
- (-1)^{|\xi||\eta|} D(\eta) \xi + D(1) \xi \eta
\\
\Phi^3_D(\xi,\eta,\vartheta) &= 
D(\xi \eta \vartheta) 
\\
& \quad
-                                  D(\xi \eta) \vartheta 
- (-1)^{|\xi|(|\eta|+|\vartheta|)}   D(\eta \vartheta) \xi 
- (-1)^{|\vartheta|(|\xi|+|\eta|)} D(\vartheta \xi) \eta
\\
& \quad
+                                  D(\xi) \eta \vartheta 
+ (-1)^{|\xi|(|\eta|+|\vartheta|)}   D(\eta) \vartheta \xi
+ (-1)^{|\vartheta|(|\xi|+|\eta|)} D(\vartheta) \xi \eta
\\
& \quad
- D(1) \xi \eta \vartheta .
\endalign
$$
In the literature, a (homogeneous) differential operator $D$
of order $\leq r$ with $D(1)=0$ is also referred to
as a (homogeneous) {\it derivation\/} of order $\leq r$.
In particular, a homogeneous derivation
$d$ of (total) degree 1 and  order 1 is precisely a differential turning
$\Cal G$ into a differential graded $R$-algebra.
\smallskip
With these preparations out of the way, consider a bigraded commutative 
$R$-algebra $\Cal G$ with 1, with $\Cal G^q_p$ zero when $q<0$ or $p<0$,
together with \newline\noindent
---
a homogeneous bracket
$[\,\cdot\,,\cdot\,]\colon\Cal G \otimes_R \Cal G \to \Cal G$
of bidegree $(0,-1)$, where \lq\lq bidegree $(0,-1)$\rq\rq\ 
means that, in given bidegrees $(q_1,p_1)$ and
$(q_2,p_2)$,
the bracket takes the form
$$
[\,\cdot\,,\cdot\,]\colon
\Cal G^{q_1}_{p_1}
\otimes
\Cal G^{q_2}_{p_2}
@>>>
\Cal G^{q_1+q_2}_{p_1+p_2-1} ; 
$$
\newline\noindent
---
a 
differential $d\colon \Cal G^*_* \to \Cal G^{*+1}_*$
of bidegree $(1,0)$
which endows $\Cal G$ 
(with respect to the total degree) with
a differential graded
$R$-algebra structure, and
\newline\noindent
---
a homogeneous differential operator 
$\Psi \colon \Cal G \to \Cal G$
of order $\leq 3$ 
with $\Psi(1) =0$ which is $\Cal G^0_0$-linear
and of bidegree $(-1,-2)$, i.~e. in bidegree $(q,p)$, $\Psi$ may
be depicted as
$$
\Psi \colon \Cal G^q_p @>>> \Cal G^{q-1}_{p-2}\quad (q \geq 1, p \geq 2).
\tag5.1
$$
In particular, $\Psi$ is zero on 
$\Cal G^0_*$, $\Cal G^*_0$, $\Cal G^*_1$.
Notice that $d$ and $\Psi$ both lower total degree by 1,
that is, are homogeneous operators on $\Cal G$ of degree $-1$.
\smallskip
We will refer to the bracket $[\,\cdot\,,\cdot\,]$ as a 
{\it quasi-Gerstenhaber bracket\/} 
and to $\Psi$
as an $h$-{\it Jacobiator\/} for the bracket $[\,\cdot\,,\cdot\,]$
provided 
$[\,\cdot\,,\cdot\,]$ and $\Psi$
satisfy (5.i)--(5.vi) below.
\newline\noindent
(5.i) The bracket
$[\,\cdot\,,\cdot\,]$ is 
graded skew-symmetric
when the total degree of
$\Cal G$ is regraded down by one,
that is, for homogeneous
$\alpha,\beta \in \Cal G$,
$$
[\alpha,\beta] =-(-1)^{(|\alpha|-1)(|\beta|-1)} [\beta,\alpha];
\tag5.2
$$
(5.ii) for each homogeneous element $\alpha$ of $\Cal G$
of bidegree $(q,p)$,
the operation
$[\alpha,\cdot\,]$ is a derivation of $\Cal G$ of bidegree $(q-1,p)$
for the multiplicative structure on $\Cal G$;
that is to say,
$[\alpha,\cdot\,]$
may be depicted as
$$
[\alpha,\cdot\,] \colon \Cal G^*_* @>>> \Cal G^{*+q-1}_{*+p} 
$$
and, for homogeneous $\beta,\gamma \in \Cal G$,
$$
[\alpha, \beta\gamma] 
=
[\alpha, \beta]\gamma +{(-1)}^{(|\alpha|-1)|\beta|} \beta[\alpha,\gamma].
\tag5.3
$$
(5.iii) The differential
$d$ behaves as a derivation for
the bracket $[\,\cdot\,,\cdot\,]$, that is,
for homogeneous $x,y \in \Cal G$,
$$
d[x,y] = [dx,y] -(-1)^{|x|} [x,dy]. 
\tag5.4
$$
(5.iv)
Given homogeneous elements $\xi,\eta,\vartheta$ of  $\Cal G$,
$$
\sum_{(\xi,\eta,\vartheta)\ \roman{cyclic}}
(-1)^{(|\xi|-1)(|\vartheta|-1)}
[\xi,[\eta,\vartheta]]
=
(-1)^{(|\xi| + |\eta|+|\vartheta|)}
\Phi^3_{d\Psi + \Psi d}(\xi, \eta,\vartheta).
\tag5.5
$$
\noindent
(5.v) the differential operator $\Psi$ has square zero and
\newline\noindent
(5.vi) the bracket
$[\,\cdot\,,\cdot\,]$ and $\Psi$ 
are related by the following requirement:
{\sl For every ordered quadruple 
$\bold a=(a_1,a_2,a_3,a_4)$ of homogeneous 
elements  of $\Cal G$,
$$
\sum_{\sigma} 
\varepsilon(\sigma)\varepsilon(\bold a,\sigma)
[\Phi^3_{\Psi}(a_{\sigma 1}, a_{\sigma 2}, a_{\sigma 3}),a_{\sigma 4}]
=
\sum_{\tau} 
\varepsilon(\tau)\varepsilon(\bold a,\tau)
\Phi^3_{\Psi}([a_{\tau 1},a_{\tau 2}], a_{\tau 3},a_{\tau 4})
\tag5.6
$$
where $\sigma$ runs through {\rm (3,1)-}shuffles and $\tau$ through 
{\rm (2,2)-}shuffles and where $\varepsilon(\sigma)$ and $\varepsilon(\tau)$ 
are the signs of the permutations $\sigma$ and $\tau$.\/}
The data $(\Cal G;d,[\,\cdot\,,\cdot\,],\Psi)$  will then be referred to as a
{\it quasi-Gerstenhaber algebra\/}.
Notice that (5.3) implies that
$[\alpha,1]=0$ for every homogeneous element $\alpha$ of $\Cal G$.
\smallskip
We note that, given an $L_{\infty}$-algebra
$\fra h$ with only two-variable and three-variable
brackets $[\cdot,\cdot]$ and $[\cdot,\cdot,\cdot]$, respectively
(and no non-zero higher order bracket operation), 
the compatibility condition which 
relates $[\cdot,\cdot]$ and $[\cdot,\cdot,\cdot]$
is exactly an identity of the kind (5.6),
when $[\cdot,\cdot,\cdot]$ is substituted for $\Phi^3_{\Psi}$.
\smallskip
A quasi-Gerstenhaber algebra having $\Psi$ zero is just an ordinary strict 
differential bigraded Gerstenhaber algebra. Indeed, in the general (quasi-) 
case, in view of the requirement (5.iv), the operation $\Psi$ {\sl measures the 
failure of the quasi-Gerstenhaber bracket $[\,\cdot\,,\cdot\,]$ 
to satisfy the graded 
Jacobi identity in a coherent fashion\/}. A strict differential bigraded
Gerstenhaber algebra having zero differential is called a {\it bigraded
Gerstenhaber algebra\/} \cite{\twilled,\,\banach}. Given a 
quasi-Gerstenhaber algebra $(\Cal G;d,[\,\cdot\,,\cdot\,],\Psi)$, we denote its 
$d$-homology by $\roman H^*_*(\Cal G)_d$. The following is straightforward.

\proclaim{Proposition 5.7}
Given a quasi-Gerstenhaber algebra $(\Cal G;d,[\,\cdot\,,\cdot\,],\Psi)$, 
the quasi-Gerstenhaber bracket $[\,\cdot\,,\cdot\,]$ induces a bracket
$$
[\,\cdot\,,\cdot\,]\colon
\roman H^{q_1}_{p_1}(\Cal G)_d
\otimes
\roman H^{q_2}_{p_2}(\Cal G)_d
@>>>
\roman H^{q_1+q_2}_{p_1+p_2-1}(\Cal G)_d
\tag5.7.1
$$
on the $d$-homology $\roman H^*_*(\Cal G)_d$ which turns 
$\roman H^*_*(\Cal G)_d$ into an ordinary bigraded Gerstenhaber algebra. \qed
\endproclaim

\smallskip
\noindent
{\smc 5.8. Relationship with Lie-Rinehart triples.\/}
We will now explain how 
quasi-Gerstenhaber algebras arise from
Lie-Rinehart triples. To this 
end, we recall that, given an ordinary  Lie-Rinehart algebra $(A,L)$, the Lie 
bracket on $L$ 
and the $L$-action on $A$
determine a Gerstenhaber bracket on the exterior $A$-algebra
$\Lambda_AL$ on $L$; for $\alpha_1,\ldots,\alpha_n \in L$, the bracket $[u,v]$ 
in $\Lambda_AL$ of $u=\alpha_1\wedge \ldots \wedge\alpha_{\ell}$ and
$v=\alpha_{\ell+1}\wedge \ldots \wedge \alpha_n$ is given by the 
expression
$$
[u,v]
=
(-1)^{\ell} 
\sum_{1 \leq j\leq \ell <k \leq n} (-1)^{j+k}
\lbrack \alpha_j,\alpha_k \rbrack \wedge 
\alpha_1\wedge \ldots \widehat{\alpha_j} \ldots \widehat{\alpha_k}
\ldots \wedge \alpha_n,
\tag5.8.1
$$
where $\ell = |u|$ is the degree of $u$, cf. \cite\bv\ (1.1). In fact, given 
the $R$-algebra $A$ and the $A$-module $L$, a bracket of the kind (5.8.1) 
yields 
a bijective correspondence between Lie-Rinehart structures on $(A,L)$ and 
Gerstenhaber algebra structures on $\Lambda_AL$. 
Our goal, which will be achieved in the next section,
is now to extend this 
observation to a relationship between Lie-Rinehart triples, 
quasi-Lie-Rinehart algebras, and 
quasi-Gerstenhaber algebras. 
\smallskip
Thus, let $(A,Q,H)$ be a pre-Lie-Rinehart triple. 
Consider the graded exterior $A$-algebra $\Lambda_A Q$, and let 
$\Cal G = \roman{Alt}_A(H,\Lambda_AQ)$, with the bigrading 
$\Cal G^q_p = \roman{Alt}^q_A(H,\Lambda^p_AQ)$ ($p,q \geq 0$). 
Suppose for the moment that $(A,Q,H)$ is merely an almost pre-Lie-Rinehart 
triple. Recall that the almost pre-Lie-Rinehart triple structure
induces operations of the kind 
(4.11.3), (4.11.4), and (4.11.6)
on the pair
$$
(\Cal A,\Cal Q) = (\roman{Alt}_A(H,A), \roman{Alt}_A(H,Q))
$$
but, at the present stage, this pair is not necessarily
a quasi-Lie-Rinehart algebra.
Consider the bigraded algebra $\roman{Alt}_A(H,\Lambda_AQ)$; 
at times we will view it
as the exterior $\Cal A$-algebra on $\Cal Q$, and we will 
accordingly write
$$
\Lambda_{\Cal A} \Cal Q = \roman{Alt}_A(H,\Lambda_AQ).
\tag5.8.2
$$
The graded skew-symmetric bracket (4.11.4) on $\Cal Q\  (=\roman{Alt}_A(H,Q))$
extends to a (bigraded) bracket 
$$
[\,\cdot\,,\cdot\,]
\colon
\Lambda_{\Cal A}\Cal Q
\otimes_R
\Lambda_{\Cal A}\Cal Q
@>>>
\Lambda_{\Cal A}\Cal Q
\tag5.8.3
$$
on $\Lambda_{\Cal A}\Cal Q=\roman{Alt}_A(H,\Lambda_AQ)$. Indeed,
with reference to the graded bracket $[\,\cdot\,,\cdot\,]$ on $\Cal Q$ spelled 
out as (4.11.4) (and written there as $[\,\cdot\,,\cdot\,]_{\Cal Q}$)
and the pairing (4.11.1), the bigraded bracket (5.8.3)
on $\Lambda_{\Cal A}\Cal Q = \roman{Alt}_A(H,\Lambda_AQ)$ is determined
by the formulas
$$
\aligned
[\alpha \beta,\gamma] 
&=
\alpha [\beta,\gamma] +{(-1)}^{|\alpha||\beta|} \beta[\alpha,\gamma],
\\
[\xi,a] &= \xi(a), 
\\ 
[\alpha,\beta] 
&= -{(-1)}^{(|\alpha|-1)(|\beta|-1)}[\beta,\alpha]
\endaligned
\tag5.8.4
$$
where $\alpha,\beta,\gamma$ are homogeneous elements of 
$\Lambda_{\Cal A}\Cal Q=\roman{Alt}_A(H,\Lambda_AQ)$,
and where $\xi \in Q$ and $a \in\Cal A$. 
\smallskip
We now construct an operation
$\Psi$
of the kind (5.1)
from the operation
$\langle \cdot,\cdot;\cdot\rangle_Q$, that is, one which formally
looks like an $h$-Jacobiator for (5.8.3).
To this end we suppose that, as an $A$-module,
at least one of $H$ or $Q$ is finitely generated and projective;
then the canonical $A$-linear morphism
from $\roman{Alt}_A(H,A) \otimes \Lambda_AQ$ to 
$\roman{Alt}_A(H,\Lambda_AQ)$
is an isomorphism of bigraded $A$-algebras.
Let $\xi_1,\dots,\xi_p \in Q$.
Now, given  a homogeneous
element $\beta$ of $\roman{Alt}_A(H,A)$, 
with reference to the operation $\langle \cdot,\cdot;\cdot\rangle_Q$
induced by $\delta$,
cf. (4.11.6),
let
$$
\Psi (\beta \xi_1\wedge \ldots\wedge \xi_p)
=
\sum_{1 \leq j <k \leq p} 
(-1)^{j+k}\langle \xi_j,\xi_k;\beta \rangle_Q \xi_1 \wedge 
\ldots \widehat \xi_j \ldots \widehat \xi_k \ldots \wedge \xi_p;
\tag5.8.5
$$
we will write $\Psi_{\delta}$ rather than just $\Psi$ whenever appropriate.
As an operator
on the graded $A$-algebra $\roman{Alt}_A(H,\Lambda_AQ)$,
$\Psi$ may be written as a finite sum of operators
which are three consecutive contractions each;
since 
an operator which consists of three consecutive contractions
is a differential operator of order $\leq 3$,
the operator 
$\Psi$ is a differential operator of order $\leq 3$.
Furthermore, 
since for $\xi,\eta \in Q$, the operation $\langle \xi,\eta;\cdot\rangle_Q$
is a derivation of 
the graded $A$-algebra $\roman{Alt}_A(H,A)$,
given homogeneous elements
$\beta_1$ and $\beta_2$ of
$\roman{Alt}_A(H,A)$,
$$
\aligned
\Psi&(\beta_1\beta_2 \xi_1 \wedge\ldots \wedge\xi_p)
\\
&=(-1)^{|\beta_1|}
\beta_1 \Psi(\beta_2 \xi_1 \wedge\ldots \wedge\xi_p)
+(-1)^{(|\beta_1|+1)|\beta_2|}
\beta_2 \Psi(\beta_1 \xi_1 \wedge\ldots \wedge\xi_p)
\endaligned
\tag5.8.6
$$
A somewhat more intrinsic description of $\Psi$ 
results from the observation that the operation
$$
\Psi \colon \roman{Alt}_A^1(H,\Lambda_A^2Q) =
\roman{Hom}_A(H,\Lambda_A^2Q) @>>> A
\cong\roman{Alt}_A^0(H,\Lambda_A^0Q)
$$
is simply given by the assignment to $\chi\colon H \to \Lambda_A^2Q$
of the {\it trace\/} of the $A$-module endomorphism
$\delta \circ \chi$ of $H$
when $H$ is finitely generated and projective
as an $A$-module,
and of
the {\it trace\/} of the $A$-module endomorphism
$\chi\circ \delta$ of $\Lambda^2_AQ$
when $Q$ is finitely generated and projective as an $A$-module.
\smallskip
We now give another description of $\Psi$,
cf. (5.8.11) below, under an additional hypothesis:
Suppose that, as an $A$-module, $Q$
is finitely generated and projective of constant rank $n$.
Then the canonical $A$-module isomorphism
$$
\phi\colon\Lambda^*_A Q \to \roman{Alt}^{n-*}_A(Q,\Lambda^n_A Q)
$$
extends to an isomorphism
$$
\phi\colon
\roman{Alt}^*_A(H,\Lambda^*_A Q)
@>>>
\roman{Alt}^*_A(H, \roman{Alt}^{n-*}_A(Q,\Lambda^n_A Q))
\tag5.8.7
$$
of graded $A$-modules. In this fashion,
$\roman{Alt}^*_A(H,\Lambda^*_A Q)$
acquires a
bigraded
\linebreak
$\roman{Alt}^*_A(H, \roman{Alt}^{*}_A(Q,A))$-module structure,
induced from the graded $A$-module
$\Lambda^n_A Q$.
Further, the skew-symmetric $A$-bilinear pairing
(1.5.5) induces an operator
$$
d_2
\colon
\roman{Alt}_A^*(H, \roman{Alt}_A^{n-*}(Q,\Lambda^n_A Q))
@>>>
\roman{Alt}_A^{*-1}(H, \roman{Alt}_A^{n-(*-2)}(Q,\Lambda^n_A Q)).
\tag5.8.8
$$
This is just the operator (2.4.7$'$), suitably rewritten,
with $M=\Lambda^n_A Q$, where the degree of the latter $A$-module
forces the correct sign:
The $A$-module $\Lambda^n_A Q$
is concentrated in degree $n$,                                     
and  a form in $\roman{Alt}_A^{n-p}(Q,\Lambda^n_A Q)$
has degree $p$.
In bidegree $(q,p)$, given 
$$
\psi \in \roman{Alt}_A^q(H, \roman{Alt}_A^{n-p}(Q,\Lambda^n_A Q)),
$$
the value 
$$
d_2(\psi)\in \roman{Alt}_A^{q-1}(H, \roman{Alt}_A^{n-p+2}(Q,\Lambda^n_A Q))
$$ 
of the operator (5.8.8) is given by the formula
$$
\align
&(-1)^{|\psi|+1}\left((d_2 \psi)
(x_1,\dots,x_{q-1})\right)(\xi_{p-1},\dots ,\xi_n)
\\
&=
\sum_{p-1\leq j<k\leq n} (-1)^{j+k} \left(\psi 
(\delta(\xi_j,\xi_k),x_1,\dots,x_{q-1})\right)
(\xi_{p-1},\dots \widehat \xi_j\dots \widehat \xi_k\dots ,\xi_n),
\endalign
$$
where
$x_1,\dots,x_{q-1}\in H$ and $\xi_{p-1},\dots, \xi_n \in Q$
and, with $|\psi| = q+p$ (the correct degree would be
$|\psi| = p-q$ but modulo 2 this makes no difference), this simplifies to
$$
\aligned
&(-1)^{p}\left((d_2 \psi)
(x_1,\dots,x_{q-1})\right)(\xi_{p-1},\dots ,\xi_n)
\\
&=
\sum_{p-1\leq j<k\leq n} (-1)^{j+k} \left(\psi (x_1,\dots,x_{q-1},
\delta(\xi_j,\xi_k))\right)
(\xi_{p-1},\dots \widehat \xi_j\dots \widehat \xi_k\dots ,\xi_n);
\endaligned
\tag5.8.9
$$
cf. (2.5.4).

\proclaim{Lemma 5.8.10}
The operator $\Psi$ makes the diagram
$$
\CD
\roman{Alt}_A^*(H,\Lambda^*_AQ)
@>{\Psi}>>
\roman{Alt}_A^{*-1}(H,\Lambda^{*-2}_AQ)
\\
@V{\phi}VV
@VV{\phi}V
\\
\roman{Alt}_A^*(H, \roman{Alt}_A^{n-*}(Q,\Lambda^n_A Q))
@>>{d_2}>
\roman{Alt}_A^{*-1}(H, \roman{Alt}_A^{n-(*-2)}(Q,\Lambda^n_A Q))
\endCD
$$
commutative. 
\endproclaim

Thus, under the 
isomorphism {\rm (5.8.7)}, 
the operator $\Psi$ is induced by
the operator $d_2$ (on the right-hand side 
of (5.8.7)).

\demo{Proof}
In a given bidegree $(q,p)$, the isomorphism (5.8.7) sends 
$\alpha \in \roman{Alt}_A^q(H, \Lambda_A^pQ)$ to
$$
\phi_{\alpha} \in \roman{Alt}_A^q(H, \roman{Alt}_A^{n-p}(Q,\Lambda^n_A Q))
$$
determined by the identity
$$
\left(\phi_{\alpha}(x_1,\ldots,x_q)\right)(\xi_{p+1},\ldots,\xi_n)
=
\left(\alpha(x_1,\ldots,x_q)\right)\wedge \xi_{p+1}\wedge \ldots \wedge\xi_n,
$$
for arbitrary $x_1,\ldots,x_q \in H$ and 
$\xi_{p+1},\ldots,\xi_n \in Q$.
Under the isomorphism (5.8.7), the operator $d_2$ (on the right-hand side 
of (5.8.7)) induces an operator 
$$
\Theta =\Theta_{\delta}\colon 
\roman{Alt}_A^q(H,\Lambda^p_AQ)
@>>>
\roman{Alt}_A^{q-1}(H,\Lambda^{p-2}_AQ)
\tag5.8.11
$$
of the kind (5.1) on the left-hand side of {\rm (5.8.7)}; by construction,
for $x_1,\dots,x_{q-1}\in H$ and $\xi_{p-1},\dots, \xi_n \in Q$, 
in view of (5.8.9),
$$
\align
(-1)^p&\left((\Theta \alpha)(x_1,\dots,x_{q-1})\right)
\wedge \xi_{p-1}\wedge\dots \wedge\xi_n
\\
&=
\sum_{p-1\leq j<k\leq n} (-1)^{j+k} \left(\alpha (x_1,\dots,x_{q-1},
\delta(\xi_j,\xi_k))\right)
\wedge \xi_{p-1} \wedge \ldots \widehat \xi_j\ldots \widehat \xi_k\ldots 
\wedge \xi_n.
\endalign
$$
Let $\beta \in \roman{Alt}_A^q(H,A)$, 
$\eta_1,\dots,\eta_p\in Q$, and $\alpha = 
(\eta_1 \wedge \ldots \wedge \eta_p)\beta$;
then, for $p-1\leq j <k \leq n$,
$$
\align
\alpha (x_1,\dots,x_{q-1},\delta(\xi_j,\xi_k))
&=
(\eta_1 \wedge \ldots \wedge \eta_p)
\beta(x_1,\dots,x_{q-1},\delta(\xi_j,\xi_k))
\\
&=
\beta(x_1,\dots,x_{q-1},\delta(\xi_j,\xi_k))
\eta_1 \wedge \ldots \wedge \eta_p
\\
&=
(-1)^{q-1}\beta(\delta(\xi_j,\xi_k),x_1,\dots,x_{q-1})
\eta_1 \wedge \ldots \wedge \eta_p
\\
&=
(-1)^{q-1}\langle\xi_j,\xi_k;\beta\rangle_Q(x_1,\dots,x_{q-1})
\eta_1 \wedge \ldots \wedge \eta_p
\\
&=
(-1)^{q-1+p(q-1)}
\left(\langle\xi_j,\xi_k;\beta\rangle_Q
\eta_1 \wedge \ldots \wedge \eta_p\right)(x_1,\dots,x_{q-1})
\endalign
$$
whence
$$
\align
{}&(-1)^{p+q-1}\left((\Theta \alpha)(x_1,\dots,x_{q-1})\right)
\wedge \xi_{p-1}\wedge\dots \wedge\xi_n
\\
&=
\sum_{p-1\leq j<k\leq n} (-1)^{j+k} 
\langle\xi_j,\xi_k;\beta\rangle_Q(x_1,\dots,x_{q-1})
\eta_1 \wedge \dots \wedge \eta_p
\wedge \xi_{p-1} \wedge \ldots \widehat \xi_j\ldots \widehat \xi_k\ldots 
\wedge \xi_n.
\endalign
$$
Let $(\eta_1,\dots,\eta_p) = (\xi_1,\dots,\xi_p)$.
With
$j=p-1$ and $k=p$, this yields
$$
\align
(-1)^{p+q-1}&\left((\Theta \alpha)(x_1,\dots,x_{q-1})\right)
\wedge \xi_{p-1}\wedge\dots \wedge\xi_n
\\
&=
-\langle\xi_{p-1},\xi_p;\beta\rangle_Q(x_1,\dots,x_{q-1})
\xi_1 \wedge \ldots 
\wedge \xi_p \wedge \xi_{p+1} \wedge \ldots 
\wedge \xi_n
\endalign
$$
or, equivalently, since
$|x_1| +\dots +|x_{q-1}| = q-1$ and $|\xi_1 \wedge \ldots \wedge \xi_p|=p$,
$$
\align
{}&(-1)^{pq}\left((\Theta \alpha)(x_1,\dots,x_{q-1})\right)
\wedge \xi_{p-1}\wedge\dots \wedge\xi_n
\\
&=
- \left(\langle\xi_{p-1},\xi_p;\beta\rangle_Q
\eta_1 \wedge \ldots \wedge \eta_p\right)(x_1,\dots,x_{q-1})
\wedge \xi_{p+1} \wedge \ldots \wedge \xi_n,
\\
&=
- \left(\langle\xi_{p-1},\xi_p;\beta\rangle_Q
\xi_1 \wedge \ldots \wedge \xi_{p-2}\right)(x_1,\dots,x_{q-1})
\wedge \xi_{p-1} \wedge \ldots \wedge \xi_n.
\endalign
$$
Hence
$$
(-1)^{pq}(\Theta \alpha)
=
- \langle\xi_{p-1},\xi_p;\beta\rangle_Q \xi_1 \wedge \ldots \wedge \xi_{p-2}
\pm \ldots
$$
or, equivalently,
$$
\Theta (\beta \xi_1 \wedge \ldots \wedge \xi_p)
=
- \langle\xi_{p-1},\xi_p;\beta\rangle_Q \xi_1 \wedge \ldots \wedge \xi_{p-2}
\pm \ldots
$$
where $\ldots$ stands for terms involving
$\xi_1 \wedge \ldots \widehat \xi_j\ldots \widehat \xi_k\ldots 
\wedge \xi_p$
with $(j,k) \ne (p-1,p)$.
Consequently
$$
\Theta (\beta \xi_1 \wedge \ldots \wedge \xi_p)
=
\sum_{1 \leq j <k \leq p}
(-1)^{j+k}
 \langle\xi_j,\xi_k;\beta\rangle_Q
\xi_1 \wedge \ldots \widehat \xi_j\ldots \widehat \xi_k\ldots \wedge \xi_p.
$$
However, this is exactly the definition (5.8.5) of $\Psi$. \qed
\enddemo

In view of Remark 2.6, the operator $\Psi$
thus calculates essentially the 
Lie algebra cohomology $\roman H^*(L_{\roman{nil}},\Lambda^n_A Q)$ of the 
(nilpotent) $A$-Lie algebra $L_{\roman{nil}}$ ($=H \oplus Q$ as an $A$-module)
with values in the $A$-module $\Lambda^n_A Q$, viewed as a trivial 
$L_{\roman{nil}}$-module.
In particular, $\Psi$ is $A$-linear.

\smallskip
Suppose finally that $(A,Q,H)$ is a genuine Lie-Rinehart triple, not just 
an almost 
pre-Lie-Rinehart triple.  By Proposition 4.13, $(\Cal A,\Cal Q)$ then acquires 
a quasi-Lie-Rinehart structure. Our ultimate goal is now to prove that, 
likewise, $\Lambda_{\Cal A} \Cal Q$ endowed with the bigraded bracket (5.8.3) 
and the operation $\Psi$, cf. (5.8.5), 
(which formally looks like an $h$-Jacobiator)
acquires a 
quasi-Gerstenhaber structure. The verification of the 
requirements (5.2)--(5.4) does not present any difficulty at this stage,
and the vanishing of $\Psi \Psi$ is immediate. 
However we were so far unable to establish 
(5.5) and (5.6)
without an additional piece of structure, that of a 
{\it generator\/} of a (quasi-Gerstenhaber) bracket.
The next section is devoted to the notion of generator
and the consequences it entails.
A precise statement is given as Corollary 6.10.4 below.

\beginsection 6. Quasi-Batalin-Vilkovisky algebras
and quasi-Gerstenhaber algebras

Let $\Cal G = \Cal G^*_*$ be a bigraded commutative $R$-algebra, endowed with
a bigraded bracket
$[\,\cdot\,,\cdot\,] \colon \Cal G \otimes_R \Cal G \to \Cal G$
of bidegree $(0,-1)$ which is graded skew-symmetric when the total
degree is regraded down by 1.
Extending terminology due to Koszul, cf. the definition 
of $[\,\cdot\,,\cdot\,]_D$ 
on p.~260 of \cite\koszulon, we will say that an $R$-linear operator $\Delta$
on $\Cal G$  of bidegree $(0,-1)$ {\it generates\/} the 
bracket $[\,\cdot\,,\cdot\,]$ provided, for every homogeneous $a, b \in \Cal G$,
$$
[a,b] = (-1)^{|a|}\left(
\Delta(ab) -(\Delta a) b - (-1)^{|a|} a (\Delta b)\right)
\quad \left(= (-1)^{|a|} \Phi^2_{\Delta}(a,b)\right);
\tag6.1
$$
we then refer to the operator $\Delta$ as a {\it generator\/}. 
\smallskip
In particular,
let $(\Cal G;d,[\,\cdot\,,\cdot\,],\Psi)$ be a quasi-Gerstenhaber 
algebra over $R$. In view of the identity (1.4) on p. 260 of \cite\koszulon, 
a generator $\Delta$ is then necessarily a differential operator on 
$\Cal G$ of order $\leq 2$.
Indeed, given a differential operator $D$, this identity reads
$$
\Phi_D^3(a,b,c) = 
\Phi_D^2(a,bc) 
-
\Phi_D^2(a,b)c
-
(-1)^{|b||c|} \Phi_D^2(a,c)b.
$$
Hence, when a differential operator $\Delta$ generates a  quasi-Gerstenhaber 
bracket $[\,\cdot\,,\cdot\,]$, 
$$
\Phi_{\Delta}^3(a,b,c) = (-1)^{|a|}
\left([a,bc]-[a,b]c-(-1)^{|b||c|} [a,c]b\right).
$$
However, by virtue of (5.3), the right-hand side of this identity is zero, 
whence $\Delta$ is necessarily of order $\leq 2$.
\smallskip
A generator $\Delta$ of a quasi-Gerstenhaber bracket
satisfies $\Delta(1) a=0$ for every $a \in \Cal G$ since, with respect 
to the multiplication map on $\Cal G$, the  
quasi-Gerstenhaber bracket behaves 
as a derivation of the appropriate degree in each variable of the bracket,
cf. (5.3). We will say that a generator $\Delta$  is {\it strict\/} provided 
$\Delta(1) = 0$ and
$$
\align
d\Delta + \Delta d &= 0,
\tag6.2
\\
d \Psi +\Delta\Delta + \Psi d&= 0;
\tag6.3
\\
\Delta \Psi + \Psi \Delta&= 0;
\tag6.4
\endalign
$$
a strict generator will henceforth often be written as $\partial$. 

\smallskip
Let $\Cal G$ be a bigraded commutative algebra, with differential operators
$$
d\colon \Cal G^*_* \to \Cal G^{*+1}_*,
\quad
\Delta\colon \Cal G^*_* \to \Cal G^*_{*-1},
\quad
\Psi\colon \Cal G^*_* \to \Cal G^{*-1}_{*-2},
$$
having orders, respectively, $\leq 1$, $\leq 2$, $\leq 3$, 
and having the properties $d(1) =0$, $\Delta(1)=0$, $\Psi(1)=0$. 
We will say that $(\Cal G;d,\Delta,\Psi)$ is a quasi-Batalin-Vilkovisky algebra
provided $dd=0$ and the operators $d,\Delta,\Psi$ satisfy
the identities {\rm (6.2)--(6.4)} as well as the identity
$$
\Psi \Psi = 0.
\tag6.5
$$
Thus $(\Cal G;d,\Delta,\Psi)$ is a quasi-Batalin-Vilkovisky algebra
if and only if, on the totalization,
the operator $\Cal D = d + \Delta + \Psi$ has square zero.

\smallskip\noindent
{\smc 6.6. From quasi-Batalin-Vilkovisky algebras to
quasi-Gerstenhaber algebras.\/}

\proclaim{Theorem 6.6.1} Given a quasi-Batalin-Vilkovisky algebra
$(\Cal G;d,\partial,\Psi)$ with $\Cal G^q_p = 0$ for $q<0$ and $p<0$,
let $[\cdot,\cdot]$ be the bracket on $\Cal G$ generated by $\partial$.
Then $(\Cal G;d,[\,\cdot\,,\cdot\,],\Psi)$ is a quasi-Gerstenhaber algebra
provided, as a bigraded $\Cal G^0_0$-algebra, $\Cal G$ is generated by its 
homogeneous constituents $\Cal G^1_0$ and $\Cal G^0_1$. 
\endproclaim

To prepare for the proof, we need the following.

\proclaim{Lemma 6.6.2}
Let $\Cal G= \{\Cal G^*_*\}$ be a bigraded commutative $R$-algebra
with $\Cal G^q_p = 0$ for $q<0$ and $p<0$, let 
$\Delta \colon \Cal G^*_* \to \Cal G^*_{*-1}$
and $\Psi \colon \Cal G^*_* \to \Cal G^{*-1}_{*-2}$
be differential operators of orders $\leq 2$ and $\leq 3$,
respectively, and let
$[\cdot,\cdot]\colon \Cal G^*_* \to \Cal G^*_{*-1}$
be the bracket {\rm (6.1)} generated by $\Delta$.
Let $A= \Cal G^0_0$, and suppose that, as a bigraded $A$-algebra,
$\Cal G$ is generated by its homogeneous constituents
$\Cal G^1_0$ and $\Cal G^0_1$ and that $\Psi$ is $A$-linear. Then
$$
\Delta \Psi +\Psi\Delta = 0
\tag6.6.3
$$
if and only if, for every ordered quadruple 
$\bold a=(a_1,a_2,a_3,a_4)$ of homogeneous 
elements  of $\Cal G$, 
$$
\sum_{\sigma} 
\varepsilon(\sigma)\varepsilon(\bold a,\sigma)
[\Phi^3_{\Psi}(a_{\sigma 1}, a_{\sigma 2}, a_{\sigma 3}),a_{\sigma 4}]
=
\sum_{\tau} 
\varepsilon(\tau)\varepsilon(\bold a,\tau)
\Phi^3_{\Psi}([a_{\tau 1},a_{\tau 2}], a_{\tau 3},a_{\tau 4})
\tag6.6.4
$$
where $\sigma$ runs through {\rm (3,1)-}shuffles and $\tau$ through 
{\rm (2,2)-}shuffles and where $\varepsilon(\sigma)$ and $\varepsilon(\tau)$ 
are the signs of the permutations $\sigma$ and $\tau$.
\endproclaim

We note the identity (6.6.4) is formally the same
as (5.6), but the circumstances are now more general.
We also note that the hypotheses of the Lemma imply that
$\Delta(1) = 0$ and $\Psi(1) = 0$.

\noindent
{\smc Remark 6.6.5.\/} For a graded commutative (not bigraded) algebra,
multiplicatively generated by its homogeneous degree 1 constituent
and endowed with a suitable Batalin-Vilkovisky structure, 
formally the same identity as (5.6) has been derived in 
Theorem 3.2 of \cite\bangotwo. Our totalization $\roman {Tot}$ 
yields a notion of Batalin-Vilkovisky algebra not equivalent
to that explored in \cite\bangotwo, though; see Remark 6.17 below for details.
The distinction between the ground ring $R$ and the $R$-algebra $A$,
crucial for our approach (involving in particular
Lie-Rinehart algebras and variants thereof), complicates the situation further.
We therefore give a complete proof of the Lemma.

\demo{Proof of Lemma 6.6.2}
We start by exploring the operator
$$
\Delta \Psi +\Psi\Delta\colon \Cal G^1_3 @>>> \Cal G^0_0 =A.
$$
Let $\alpha \in \Cal G^1_0$
and $\xi_1,\xi_2,\xi_3 \in \Cal G^0_1$;
then 
$\alpha \xi_1\xi_2\xi_3 \in \Cal G^1_3$.
Since $\Psi$ is of order $\leq 3$,
$$
\Phi^4_{\Psi}(\alpha,\xi_1,\xi_2,\xi_3) = 0
$$
and, for degree reasons, this identity boils down to
$$
\Psi(\alpha \xi_1\xi_2\xi_3) =
\Psi(\alpha \xi_1\xi_2)\xi_3
+
\Psi(\alpha \xi_2\xi_3)\xi_1
+
\Psi(\alpha \xi_3\xi_1)\xi_2 .
$$
In view of the definition (6.1) of the bracket $[\cdot,\cdot]$,
$$
[\Psi(\alpha \xi_1\xi_2),\xi_3]
=
\Delta(\Psi(\alpha \xi_1\xi_2)\xi_3) -\Psi(\alpha \xi_1\xi_2)\Delta(\xi_3)
$$
whence
$$
\align
\Delta\Psi(\alpha \xi_1\xi_2\xi_3)
&=
 [\Psi(\alpha \xi_1\xi_2),\xi_3]
+[\Psi(\alpha \xi_2\xi_3),\xi_1]
+[\Psi(\alpha \xi_3\xi_1),\xi_2]
\\
&
\quad
+\Psi(\alpha \xi_1\xi_2)\Delta(\xi_3)
+\Psi(\alpha \xi_2\xi_3)\Delta(\xi_1)
+\Psi(\alpha \xi_3\xi_1)\Delta(\xi_2) .
\endalign
$$
On the other hand,
$$
\align
\Delta(\xi_1\xi_2\xi_3)
&=
 \Delta(\xi_1)\xi_2\xi_3
+\Delta(\xi_2)\xi_3\xi_1
+\Delta(\xi_3)\xi_1\xi_2
\\
&\quad
-[\xi_1,\xi_2]\xi_3
-[\xi_2,\xi_3]\xi_1
-[\xi_3,\xi_1]\xi_2
\\
[\alpha,\xi_1\xi_2\xi_3] &
=-\left(\Delta(\alpha\xi_1\xi_2\xi_3) + \alpha \Delta(\xi_1\xi_2\xi_3)\right) .
\endalign
$$
Hence
$$
\align
\Psi \Delta(\alpha \xi_1\xi_2\xi_3)
&= - \Psi\left([\alpha,\xi_1\xi_2\xi_3] +
  \alpha \Delta(\xi_1\xi_2\xi_3)\right)
\\
&=- \Psi[\alpha,\xi_1\xi_2\xi_3] 
\\
&\quad
-\Psi (\alpha \Delta(\xi_1)\xi_2\xi_3)
-\Psi (\alpha \Delta(\xi_2)\xi_3\xi_1)
-\Psi (\alpha \Delta(\xi_3)\xi_1\xi_2)
\\
&\quad
+\Psi (\alpha[\xi_1,\xi_2]\xi_3)
+\Psi (\alpha[\xi_2,\xi_3]\xi_1)
+\Psi (\alpha[\xi_3,\xi_1]\xi_2)
\endalign
$$
that is,
$$
\align
\Psi \Delta(\alpha \xi_1\xi_2\xi_3)
&=
- \Psi[\alpha,\xi_1\xi_2\xi_3] 
\\
&\quad
-\Psi (\alpha \xi_2\xi_3)\Delta(\xi_1)
-\Psi (\alpha \xi_3\xi_1)\Delta(\xi_2)
-\Psi (\alpha \xi_1\xi_2)\Delta(\xi_3)
\\
&\quad
+\Psi (\alpha[\xi_1,\xi_2]\xi_3)
+\Psi (\alpha[\xi_2,\xi_3]\xi_1)
+\Psi (\alpha[\xi_3,\xi_1]\xi_2)
\endalign
$$
since $\Psi$ is $A$-linear.
Exploiting the identity
$$
\Psi[\alpha,\xi_1\xi_2\xi_3] 
=
\Psi[\alpha,\xi_1]\xi_2\xi_3
+\Psi[\alpha,\xi_2]\xi_3\xi_1 
+\Psi[\alpha,\xi_3]\xi_1\xi_2 
$$
we conclude
$$
\align
(\Delta\Psi+\Psi \Delta)(\alpha \xi_1\xi_2\xi_3)
&=
 [\Psi(\alpha \xi_1\xi_2),\xi_3]
+[\Psi(\alpha \xi_2\xi_3),\xi_1]
+[\Psi(\alpha \xi_3\xi_1),\xi_2]
\\
&\quad
+\Psi (\alpha[\xi_1,\xi_2]\xi_3)
+\Psi (\alpha[\xi_2,\xi_3]\xi_1)
+\Psi (\alpha[\xi_3,\xi_1]\xi_2)
\\
&\quad
-\Psi([\alpha,\xi_1]\xi_2\xi_3)
-\Psi([\alpha,\xi_2]\xi_3\xi_1)
-\Psi([\alpha,\xi_3]\xi_1\xi_2) .
\endalign
$$
Thus the graded commutator $\Delta \Psi +\Psi\Delta$
vanishes on $\Cal G^1_3$ if and only if, for every 
$\alpha \in \Cal G^1_0$ and 
$\xi_1,\xi_2,\xi_3 \in \Cal G^0_1$,
$$
\align
 [\Psi(\alpha \xi_1\xi_2),\xi_3]&
+[\Psi(\alpha \xi_2\xi_3),\xi_1]
+[\Psi(\alpha \xi_3\xi_1),\xi_2] =
\\
&
\Psi([\alpha,\xi_1]\xi_2\xi_3)
+\Psi([\alpha,\xi_2]\xi_3\xi_1)
+\Psi([\alpha,\xi_3]\xi_1\xi_2)
\\
&
-\Psi (\alpha[\xi_1,\xi_2]\xi_3)
-\Psi (\alpha[\xi_2,\xi_3]\xi_1)
-\Psi (\alpha[\xi_3,\xi_1]\xi_2) ;
\endalign
$$
since $\Psi(\xi_1\xi_2\xi_3)=0$, the latter identity is equivalent to
$$
\align
 [\Psi(\alpha \xi_1\xi_2),\xi_3]&
-[\Psi(\xi_1\xi_2\xi_3),\alpha]
+[\Psi(\xi_2\xi_3\alpha),\xi_1]
-[\Psi(\xi_3\alpha\xi_1),\xi_2] 
=
\\
&
\Psi([\alpha,\xi_1]\xi_2\xi_3)
+\Psi([\alpha,\xi_2]\xi_3\xi_1)
+\Psi([\alpha,\xi_3]\xi_1\xi_2)
\\
&
+\Psi ([\xi_1,\xi_2]\alpha\xi_3)
+\Psi ([\xi_2,\xi_3]\alpha\xi_1)
+\Psi ([\xi_3,\xi_1]\alpha\xi_2) .
\endalign
$$
With the more neutral notation $(a_1,a_2,a_3,a_4) = (\alpha,\xi_1,\xi_2,\xi_3)$
since, for degree reasons,
$$
\Phi^3_{\Psi}(a_1,a_2,a_3)=-\Psi(a_1 a_2a_3), 
\quad
\Phi^3_{\Psi}([a_1,a_2],a_3,a_4)=-\Psi([a_1,a_2]a_3a_4)
$$
etc.,
the identity takes the form
$$
\align
 [\Phi^3_{\Psi}(a_1,a_2,a_3),a_4]&
-[\Phi^3_{\Psi}(a_2,a_3,a_4),a_1]
+[\Phi^3_{\Psi}(a_3,a_4,a_1),a_2]
-[\Phi^3_{\Psi}(a_4,a_1,a_2),a_3] 
=
\\
&
\Phi^3_{\Psi}([a_1,a_2],a_3,a_4)
-\Phi^3_{\Psi}([a_1,a_3],a_2,a_4)
+\Phi^3_{\Psi}([a_1,a_4],a_2,a_3)
\\
&
+\Phi^3_{\Psi} ([a_2,a_3],a_1,a_4)
-\Phi^3_{\Psi} ([a_2,a_4],a_1,a_3)
+\Phi^3_{\Psi} ([a_3,a_4],a_1,a_2) .
\endalign
$$
This is the identity (6.6.4) for the special case where the elements
$a_1,a_2,a_3,a_4$ are from $\Cal G^1_0 \cup \Cal G^0_1$.
The operator $\Delta$ being of order $\leq 2$ means precisely that
the bracket $[\cdot,\cdot]= \pm \Phi^2_{\Delta}$ 
(generated by it) behaves as a derivation in each argument and, accordingly,
the operator $\Psi$ being of order $\leq 3$ means that the operation
$\Phi^3_{\Psi}$ is a derivation in each of its three arguments.
The equivalence between the identities (6.6.3) and (5.7.1) 
for arbitrary arguments is now etablished by induction on the degrees
of the arguments. \qed
\enddemo

\demo{Proof of Theorem {\rm 6.6.1}} The quasi-Gerstenhaber bracket 
$[\,\cdot\,,\cdot\,]$ on
$\Cal G$  is that generated by $\Delta=\partial$
via (6.1). This bracket is plainly
graded skew-symmetric in the correct sense, and the reasoning in Section 1 of
\cite\koszulon\ shows that this bracket satisfies  
the identities (5.3)--(5.5).
In particular, the identity (5.5) is a consequence of 
the identity (6.3):
This identity  may be rewritten as
$$
d \Psi  + \Psi d= -\Delta \Delta.
\tag6.3$'$
$$
Hence, given homogeneous elements $\xi,\eta,\vartheta$ of $\Cal G$,
the identity (5.5) takes the form
$$
\sum_{(\xi,\eta,\vartheta)\ \roman{cyclic}}
(-1)^{(|\xi|-1)(|\vartheta|-1)}
[\xi,[\eta,\vartheta]]
=
-(-1)^{(|\xi| + |\eta|+|\vartheta|)}
\Phi^3_{\Delta^2}(\xi, \eta,\vartheta).
\tag5.5$'$
$$
This is exactly the identity in line -5 on p. 260 of \cite\koszulon, which 
measures the failure of the bracket $[\,\cdot\,,\cdot\,]$ 
to satisfy the graded Jacobi identity
in terms of the square $\Delta^2$ of the generating 
operator $\Delta$. 
The identity (5.6) holds by virtue of Lemma 6.6.2. \qed
\enddemo
An observation due to Koszul \cite\koszulon\ (p.~261) extends to the present 
case in the following fashion: For any  quasi-Batalin-Vilkovisky algebra
$(\Cal G; d,\partial,\Psi)$, the operator $\partial$ (which is strict by 
assumption) behaves as a derivation for the quasi-Gerstenhaber bracket 
$[\,\cdot\,,\cdot\,]$, up to a suitable {\it correction term\/}
which we now determine:
The identity in line 6 on p. 261 of \cite\koszulon\ implies that,
for homogeneous $a,b \in \Cal G$,
$$
\partial [a,b] - ([\partial a,b] -(-1)^{|a|} [a,\partial b]), 
=
(-1)^{|a|}\Phi^2_{\partial^2}(a,b).
$$
Since, by virtue of (6.3),
$\partial \partial + d \Psi + \Psi d = 0$,
we conclude
$$
\partial [a,b] - ([\partial a,b] -(-1)^{|a|} [a,\partial b]) 
=(-1)^{|a|-1}\Phi^2_{d\Psi + \Psi d}(a,b) .
$$
The correction term
$\Phi^2_{d\Psi + \Psi d}(a,b)$
is {\sl plainly an instance of the occurrence
of a homotopy\/}.
We also note that, in view of (6.1), a generator, even if strict,
behaves as a derivation for the multiplication
of $\Cal G$ only if the quasi-Gerstenhaber bracket $[\,\cdot\,,\cdot\,]$ 
is zero.
\smallskip
A quasi-Batalin-Vilkovisky algebra having $\Psi$ zero is just an ordinary
differential bigraded Batalin-Vilkovisky algebra, and a differential bigraded
Batalin-Vilkovisky algebra having zero differential is called a
{\it bigraded Batalin-Vilkovisky algebra\/} \cite{\twilled,\,\banach}.
Maintaining notation introduced in the previous section, given a 
quasi-Batalin-Vilkovisky algebra $(\Cal G;d,\partial,\Psi)$, we denote its 
$d$-homology by $\roman H^*_*(\Cal G)_d$; Proposition 5.7 above says that
$\roman H^*_*(\Cal G)_d$ inherits a bigraded Gerstenhaber bracket.
Plainly, under the present circumstances this homology inherits more structure;
indeed, the proof of the following is straightforward and left to the reader.

\proclaim{Proposition 6.7}
Given a quasi-Batalin-Vilkovisky algebra $(\Cal G;d,\partial,\Psi)$,
the strict operator $\partial$ induces a generator
$$
\partial\colon
\roman H^*_*(\Cal G)_d
@>>>
\roman H^*_{*-1}(\Cal G)_d
\tag6.7.1
$$
for the bigraded Gerstenhaber bracket on its $d$-homology 
$\roman H^*_*(\Cal G)_d$ and hence turns the latter
into a bigraded Batalin-Vilkovisky algebra. \qed
\endproclaim

A quasi-Batalin-Vilkovisky algebra has an invariant which is finer than just 
ordinary homology, though: Let $(\Cal G;d,\partial,\Psi)$ be a 
quasi-Batalin-Vilkovisky algebra, and consider the following
$\roman{Tot} \Cal G$ of $\Cal G$ given by
$$
(\roman{Tot} \Cal G)_n = \sum_{q-p = n} \Cal G^q_p =
\Cal G^n_0
\oplus
\Cal G^{n+1}_1
\oplus
\ldots
\oplus
\Cal G^{n+k}_k
\oplus
\ldots
\tag6.7.2
$$
This totalization is forced by the isomorphism (6.8) and by Theorem 6.10
below.
In a given bidegree $(q,p)$, the operators $d,\partial,\Psi$ may be depicted as
$$
d \colon \Cal G^q_p \to \Cal G^{q+1}_p,
\quad
\partial \colon \Cal G^q_p \to \Cal G^q_{p-1},
\quad
\Psi \colon \Cal G^q_p \to \Cal G^{q-1}_{p-2},
\tag6.7.3
$$
and the defining properties (6.2)--(6.5)  say that the sum
$$
\Cal D = d + \partial + \Psi
\tag6.7.4
$$
is a square zero operator on 
$\roman{Tot} \Cal G$, i.~e. a differential. Consider 
the ascending filtration $\{F_r\}_{r\geq 0}$ of $\roman{Tot} \Cal G$ given by
$$
F_r (\roman{Tot} \Cal G)_n = \sum_{q-p = n,\,p \leq r} \Cal G^q_p
=\Cal G^n_0
\oplus
\Cal G^{n+1}_1
\oplus
\ldots
\oplus
\Cal G^{n+r}_r
\tag6.7.5
$$
This filtration gives rise to a spectral sequence
$$
({\roman E}_*^*(r),d(r)),
\quad
d(r) \colon {\roman E}_p^q(r)
@>>> {\roman E}_{p-r}^{q-r+1}(r)
\tag6.7.6
$$
having
$$
({\roman E}(0),d(0)) = 
(\Cal G,d)
\tag6.7.7
$$
whence
$$
({\roman E}(1),d(1))
= 
(\roman H^*_*(\Cal G)_d,\partial),
\tag6.7.8
$$
which is the bigraded homology Batalin-Vilkovisky algebra spelled out in 
Proposition 6.7 above. This spectral sequence is an invariant for the
quasi-Batalin-Vilkovisky algebra $\Cal G$ which is {\it finer\/} than just the 
bigraded homology Batalin-Vilkovisky algebra 
$(\roman H^*_*(\Cal G)_d,\partial)$.
\smallskip
We will now take up and extend the discussion in (5.8) and
describe how quasi-Gerstenhaber and 
quasi-Batalin-Vilkovisky algebras arise from Lie-Rinehart triples.
To this end, let $(A,H,Q)$ be a pre-Lie-Rinehart triple and suppose that, as 
an $A$-module, $Q$ is finitely generated and projective, of constant rank $n$.
Consider the graded exterior $A$-algebra $\Lambda_A Q$, and let
$\Cal G = \roman{Alt}_A(H,\Lambda_AQ)$, with
$\Cal G^q_p = \roman{Alt}^q_A(H,\Lambda^p_AQ)$;
this is a bigraded commutative $A$-algebra. The Lie-Rinehart differential $d$,
with respect to the canonical graded $(A,H)$-module structure on 
$\Lambda_A Q$, turns $\Cal G$ into a differential graded $R$-algebra.
Our aim is to 
determine when $(A,Q,H)$ is a genuine Lie-Rinehart triple 
in terms of conditions on $\Cal G$.
\smallskip
The graded $A$-module $\roman{Alt}^{*}_A(Q,\Lambda^n_A Q)$
acquires a canonical graded $(A,H)$-module structure.
Further, 
since $(A,H,Q)$ is a pre-Lie-Rinehart triple (not just an almost 
pre-Lie-Rinehart triple), the 
canonical bigraded $A$-module isomorphism (5.8.7)
is now an isomorphism
$$
\phi\colon (\roman{Alt}^*_A(H,\Lambda^*_A Q),d)
@>>>
(\roman{Alt}^*_A(H, \roman{Alt}^{n-*}_A(Q,\Lambda^n_A Q)),d)
\tag6.8
$$
of Rinehart complexes, with reference to the graded $(A,H)$-module structures 
on $\Lambda^*_A Q$ and $\roman{Alt}^{n-*}_A(Q,\Lambda^n_A Q)$.
We will say that $(A,H,Q)$ is {\it weakly orientable\/} if 
$\Lambda^n_A Q$ is a free $A$-module, that is, if there is an $A$-module 
isomorphism $\omega \colon \Lambda^n_A Q \to A$, and $\omega$ will then be 
referred to as a {\it weak orientation form\/}. Under the circumstances of
Example 1.4.1, this notion of weak orientability means that the foliation 
$\Cal F$ is transversely orientable, with transverse volume form $\omega$.
For a general pre-Lie-Rinehart triple $(A,H,Q)$, we will say that a weak 
orientation form $\omega$ is {\it invariant\/} provided it is invariant under 
the $H$-action; we will then refer to $\omega$ as an {\it orientation\/} form,
and we will say that $(A,H,Q)$ is {\it orientable\/}. In the situation of 
Example 1.4.1, 
with a grain of salt, an orientation form in this sense amounts 
to an orientation for the 
\lq\lq space of leaves\rq\rq, that is, with reference to the 
spectral sequence (2.9.1), the class in the top basic
cohomology group $\roman E_2^{n,0}$ (cf. 2.10(i)) of such a form
is non-zero and generates this cohomology group. Likewise, in the situation 
of Example 1.4.2, an orientation form is a {\it holomorphic volume form\/},
and the requirement that an (invariant) orientation form exist is 
precisely the Calabi-Yau condition.
\smallskip
Let $(A,H,Q)$ be a general orientable pre-Lie-Rinehart triple, and let 
$\omega$ be an invariant orientation form. Then $\omega$ induces an isomorphism
$\roman{Alt}^{*}_A(Q,\Lambda^n_A Q)\to\roman{Alt}^{*}_A(Q,A)$
of graded $(A,H)$-modules and hence an isomorphism
$$
\phi^{\omega}\colon
(\roman{Alt}^*_A(H,\Lambda^*_A Q),d)
@>>>
(\roman{Alt}^*_A(H, \roman{Alt}^{n-*}_A(Q,A)),d_0)
\tag6.9
$$
of Rinehart complexes. Here, on the right-hand side of (6.9), 
the operator $d_0$ is that given earlier as 
(2.4.5), with the orders of $H$ and $Q$ interchanged.
On the right-hand side of (6.9), 
we have as well the operator 
$d_1$ given as (2.4.6) and the operator $d_2$ given as (2.4.7) 
(the order of $H$ and $Q$ being interchanged), cf. also (5.8.8).
The operator $d_1$ induces an operator
$$
\Delt_{\omega}
\colon
\roman{Alt}^*_A(H,\Lambda^*_A Q)
@>>>
\roman{Alt}^*_A(H,\Lambda^{*-1}_A Q)
\tag6.10.1
$$
on the left-hand side of (6.9) by means of the
the relationship
$$
\phi^{\omega}_{\Delt_{\omega}(\alpha)} = 
(-1)^{n+1}d_1 (\phi^{\omega}_\alpha),\quad \alpha \in \Lambda^*_A Q.
$$
By Lemma 5.8.10, the operator $d_2$ on the right-hand side 
of (6.9) corresponds to the operator $\Psi_{\delta}$ on the left-hand side 
of (6.9) given 
as (5.8.5) above. Notice that $\Delt_{\omega}$ is an $R$-linear 
operator on $\Cal G^*_*=\roman{Alt}^*_A(H,\Lambda^*_A Q)$  
of bidegree $(0,-1)$ which looks like a generator for the 
corresponding bracket (5.8.3). We will now describe the circumstances
where $\Delt_{\omega}$ is a generator.

\proclaim{Theorem 6.10}
Let $(A,H,Q)$ be an  orientable pre-Lie-Rinehart triple, with invariant 
orientation form $\omega$. 
If $(A,H,Q)$ is a genuine Lie-Rinehart triple,
then $(\roman{Alt}^*_A(H,\Lambda^*_A Q);d, \Delta_{\omega},\Psi_{\delta})$
is a quasi-Batalin-Vilkovisky algebra, and 
$\Delta_{\omega}$ is a strict generator for the bracket 
$[\,\cdot\,,\cdot\,]$ given by {\rm (5.8.3)}.
Conversely, under the additional hypothesis that $H$ satisfy the property 
{\rm P}, if
$(\roman{Alt}^*_A(H,\Lambda^*_A Q),d, \Delta_{\omega},\Psi_{\delta})$ is a 
quasi-Batalin-Vilkovisky algebra,
then $(A,H,Q)$ is a genuine Lie-Rinehart triple.
\endproclaim

\demo{Proof} 
We note first that, when $(\Cal A,d_0,d_1,d_2)$ is a multialgebra,
so is $(\Cal A,d_0,-d_1,d_2)$.
Furthermore, when 
$\left(\roman{Alt}^*_A(H,\roman{Alt}^*_A(Q,A)),d_0,-d_1,d_2\right)$
is a multialgebra, 
$$
\left(\roman{Alt}^*_A(H,\roman{Alt}^*_A(Q,\Lambda_A^nQ)),d_0,-d_1,d_2\right)
$$
is a multicomplex,
the operators $d_j$ ($0 \leq j \leq 2$) 
(where the notation $d_j$ is abused somewhat)
being the induced ones, with the 
correct sign, that is,
$$
\omega_*(d_j (\cdot)) = (-1)^n d_j \omega_*( (\cdot))
$$
where $\omega_*$ is the induced bigraded morphism of degree $n$.
Hence the equivalence between the Lie-Rinehart triple
and quasi-Batalin-Vilkovisky properties is straightforward, 
in view of Theorem 2.7 and
Theorem 6.6.1.
In particular, the identities (2.1.4.2)--(2.1.4.5) correspond to
the identities (6.2)--(6.5) which characterize
$(\roman{Alt}^*_A(H,\Lambda^*_A Q);d,\Delta_{\omega},\Psi_{\delta})$
being a quasi-Batalin-Vilkovisky algebra. 
It remains to show that, when $(A,H,Q)$ is a genuine 
Lie-Rinehart triple, the operator
$\Delta_{\omega}$ (given by (6.10.1))
is indeed a strict generator for the bigraded bracket (5.8.3)
to which the rest of the proof is devoted.
\smallskip
\noindent
{\smc 6.10.2. Verification of the generating property.\/}
We note first that,
in view of the derivation properties of a quasi-Gerstenhaber bracket,
it suffices to establish the generating property (6.1) on $\Lambda_AQ$, 
viewed as the bidegree $(0,*)$-constituent of 
$\roman{Alt}_A(H,\Lambda_AQ) =\Lambda_{\Cal A}\Cal Q$. 
To make the operator  $\Delta_{\omega}$
somewhat more explicit,
we note that
the pairing (1.5.2.Q) and the choice of $\omega$
determine a {\it generalized\/} $Q$-connection 
$$
\nabla \colon Q \otimes \Lambda_A^nQ @>>> \Lambda_A^nQ
$$
on $\Lambda_A^nQ$ determined by requiring that the diagram
$$
\CD
Q \otimes_R \Lambda_A^nQ @>{\nabla}>> \Lambda_A^nQ
\\
@V{\roman{Id}\otimes \omega}VV
@VV{\omega}V
\\
Q \otimes_R A @>>{(1.5.2.Q)}> A
\endCD
$$
be commutative, and the 
multialgebra compatibility property $d_0 d_1 + d_1 d_0 = 0$
(cf. (2.1.4.2)) is 
equivalent to this generalized $Q$-connection being compatible
with the $H$-module structures. 
In the situation of Example 1.4.2, such a generalized $Q$-connection on 
$\Lambda_A^nQ$ amounts to a flat holomorphic connection on
the highest exterior power of the holomorphic tangent bundle.
In the present general case,
the operator $d_1$ on $\roman{Alt}_A(H, \roman{Alt}_A(Q,A))$ (given by (2.5.3)
then corresponds to an operator
$$
\do \colon \roman{Alt}_A(H,\roman{Alt}_A^p(Q,\Lambda_A^nQ)) 
@>>> \roman{Alt}_A(H,\roman{Alt}_A^{p+1}(Q,\Lambda_A^nQ))
\quad (p \geq 0)
$$
determined by the commutativity of the diagram
$$
\CD
\roman{Alt}_A(H,\roman{Alt}_A^p(Q,\Lambda_A^nQ))
@>{\do}>> \roman{Alt}_A(H,\roman{Alt}_A^{p+1}(Q,\Lambda_A^nQ))
\\
@V{}VV
@VV{}V
\\
\roman{Alt}_A(H,\roman{Alt}_A^p(Q,A)) 
@>>{(-1)^nd_1}> \roman{Alt}_A(H,\roman{Alt}_A^{p+1}(Q,A))
\endCD
$$
whose vertical arrows are induced by $\omega$.
Consider, then, the operator $\Del$
determined by the requirement that
the diagram
$$
\CD
\roman{Alt}_A(H,\Lambda^p_AQ)
@>{\phi}>>
\roman{Alt}_A(H,\roman{Alt}^{n-p}_A(Q,\Lambda_A^nQ))
\\
@V{\Del}VV
@VV{-\do }V
\\
\roman{Alt}_A(H,\Lambda^{p-1}_AQ)
@>>{\phi}>
\roman{Alt}_A(H,\roman{Alt}^{n-(p-1)}_A(Q,\Lambda_A^nQ))
\endCD
$$
be commutative. This operator coincides with
the operator
$\Delta_{\omega}$
but we prefer to use a neutral notation.
In view of the derivation properties of a quasi-Gerstenhaber bracket,
to establish the generating property,
it will suffice to study the restriction
$$
D \colon \Lambda^p_AQ @>>> \Lambda^{p-1}_AQ\quad (1 \leq p \leq n)
$$
of this operator.
\smallskip
Given $\alpha \in \Lambda^p_AQ$, we will write
$\phi_{\alpha} \in \roman{Alt}^{n-p}_A(Q,\Lambda^n_AQ)$
for the image under $\phi$ so that, for $\xi_{p+1},\ldots,\xi_n$,
$$
\phi_{\alpha}(\xi_{p+1},\ldots,\xi_n) = \alpha \wedge 
\xi_{p+1} \wedge \ldots \wedge \xi_n.
$$
Let $\alpha_1$ and $\alpha_2$ be homogeneous elements of $\Lambda^*_AQ$.
We will now establish the generating property
$$
(-1)^{|\alpha_1|} [\alpha_1,\alpha_2] = 
\Del(\alpha_1 \alpha_2) -(\Del\alpha_1) \alpha_2 
- (-1)^{|\alpha_1|} \alpha_1 (\Del\alpha_2) .
\tag6.10.3
$$
Let $\beta \in \Lambda_A^{n+1-|\alpha_1| -|\alpha_2|}Q$;
it will suffice to study the expression
$$
(\Del(\alpha_1 \alpha_2)) \wedge \beta
-((\Del\alpha_1) \alpha_2) \wedge \beta
- (-1)^{|\alpha_1|} (\alpha_1 (\Del\alpha_2))\wedge \beta
-(-1)^{|\alpha_1|} [\alpha_1,\alpha_2] \wedge \beta
\in \Lambda_A^nQ
$$
or, equivalently, the expression
$$
\phi_{\Del(\alpha_1 \alpha_2)}(\beta) 
-\phi_{(\Del\alpha_1) \alpha_2} (\beta)
- (-1)^{|\alpha_1|} \phi_{\alpha_1 \Del\alpha_2}(\beta)
-(-1)^{|\alpha_1|}
\phi_{[\alpha_1, \alpha_2]}(\beta) \in \Lambda_A^nQ .
$$
To this end, we note first that
$$
\align
-\phi_{\Del(\alpha_1 \alpha_2)}(\beta) 
&= (\do \phi_{\alpha_1 \alpha_2}) (\beta)
\\
-\phi_{(\Del \alpha_1)\alpha_2}(\beta) 
&= 
-(\Del \alpha_1)\wedge \alpha_2 \wedge\beta 
=(\do \phi_{\Del\alpha_1})(\alpha_2 \wedge \beta)
\\
-\phi_{\alpha_1 (\Del \alpha_2)}(\beta) 
&= 
-\alpha_1\wedge (\Del \alpha_2) \wedge\beta 
=- (-1)^{|\alpha_1|(|\alpha_2|-1)}
(\Del \alpha_2)\wedge \alpha_1 \wedge\beta 
\\
&
=(-1)^{|\alpha_1|(|\alpha_2|-1)}
(\do \phi_{\Del\alpha_2})(\alpha_1 \wedge \beta) .
\endalign
$$
Let $\vartheta_1,\vartheta_2 \in Q$ and
$\xi_2,\dots, \xi_n \in Q$.
Letting $\xi_1 = \vartheta_2$ we obtain
$$
\align
(\do \phi_{\vartheta_1})
&
(\vartheta_2,\xi_2,\dots, \xi_n)
=(\do \phi_{\vartheta_1})
(\xi_1,\xi_2,\dots, \xi_n)
\\
&=
\sum_{1\leq j\leq n} 
(-1)^{j-1}
\nabla_{\xi_j}
(\vartheta_1 \wedge \xi_1 \wedge\dots \widehat \xi_j\ldots\wedge\xi_n)
\\
&
\quad
+
\sum_{1\leq j<k\leq n} 
(-1)^{j+k}
\vartheta_1 \wedge [\xi_j,\xi_k]
\wedge\xi_1\wedge\dots\widehat\xi_j\dots\widehat\xi_k\ldots\wedge\xi_n
\endalign
$$ 
and a straightforward calculation gives
$$
\align
(\do \phi_{\vartheta_1})
(\xi_1,\xi_2,\dots, \xi_n)
&=
\nabla_{\vartheta_2}
(\vartheta_1 \wedge \xi_2 \wedge \ldots\wedge\xi_n)
\\
&\quad
+
\sum_{2\leq j\leq n} 
(-1)^{j-1}
\nabla_{\xi_j}
(\vartheta_1 \wedge \vartheta_2 \wedge \xi_2
\wedge\dots \widehat \xi_j\ldots\wedge\xi_n)
\\
&\quad
+
\sum_{1<k\leq n} 
(-1)^{1+k}
\vartheta_1 \wedge [\vartheta_2,\xi_k]
\wedge\xi_2\wedge\dots\widehat\xi_k\ldots\wedge\xi_n
\\
&\quad
+
\sum_{2\leq j<k\leq n} 
(-1)^{j+k}
\vartheta_1 \wedge [\xi_j,\xi_k]
\wedge\vartheta_2\wedge \xi_2\wedge\dots\widehat\xi_j\dots
\widehat\xi_k\ldots\wedge\xi_n .
\endalign
$$
Likewise letting $\xi_1 = \vartheta_1$ we obtain
$$
\align
(\do \phi_{\vartheta_2})
&
(\vartheta_1,\xi_2,\dots, \xi_n)
=(\do \phi_{\vartheta_2})
(\xi_1,\xi_2,\dots, \xi_n)
\\
&=
\sum_{1\leq j\leq n} 
(-1)^{j-1}
\nabla_{\xi_j}
(\vartheta_2 \wedge \xi_1 \wedge\dots \widehat \xi_j\ldots\wedge\xi_n)
\\
&\quad
+
\sum_{1\leq j<k\leq n} 
(-1)^{j+k}
\vartheta_2 \wedge [\xi_j,\xi_k]
\wedge\xi_1\wedge\dots\widehat\xi_j\dots\widehat\xi_k\ldots\wedge\xi_n
\endalign
$$
and again a calculation yields
$$
\align
(\do \phi_{\vartheta_2})
(\xi_1,\xi_2,\dots, \xi_n)
&=
\nabla_{\vartheta_1}
(\vartheta_2 \wedge \xi_2 \wedge \ldots\wedge\xi_n)
\\
&\quad
+
\sum_{2\leq j\leq n} 
(-1)^{j-1}
\nabla_{\xi_j}
(\vartheta_2 \wedge \vartheta_1 \wedge \xi_2
\wedge\dots \widehat \xi_j\ldots\wedge\xi_n)
\\
&\quad
+
\sum_{1<k\leq n} 
(-1)^{1+k}
\vartheta_2 \wedge [\vartheta_1,\xi_k]
\wedge\xi_2\wedge\dots\widehat\xi_k\ldots\wedge\xi_n
\\
&\quad
+
\sum_{2\leq j<k\leq n} 
(-1)^{j+k}
\vartheta_2 \wedge [\xi_j,\xi_k]
\wedge\vartheta_1\wedge \xi_2\wedge\dots\widehat\xi_j\dots
\widehat\xi_k\ldots\wedge\xi_n .
\endalign
$$
Next, let $\alpha =\vartheta_1 \wedge\vartheta_2$;
then the corresponding $(n-2)$-form 
$\phi_{\alpha} \in \roman{Alt}_A^{n-2}(Q,\Lambda_A^nQ)$
has degree 2 whence
$$
\align
(-1)^3(\do \phi_{\alpha})
&
(\xi_2,\dots, \xi_n)
=
-\sum_{2\leq j\leq n} 
(-1)^{j-1}
\nabla_{\xi_j}
(\vartheta_1 \wedge\vartheta_2 
\wedge 
\xi_2 \wedge\dots \widehat \xi_j\ldots\wedge\xi_n)
\\
&\quad
+
\sum_{2\leq j<k\leq n} 
(-1)^{j+k}
\vartheta_1 \wedge\vartheta_2 \wedge [\xi_j,\xi_k]
\wedge\xi_2\wedge\dots\widehat\xi_j\dots\widehat\xi_k\ldots\wedge\xi_n
\endalign
$$
that is
$$
\align
(\do \phi_{\alpha})
&
(\xi_2,\dots, \xi_n)
=
\sum_{2\leq j\leq n} 
(-1)^{j-1}
\nabla_{\xi_j}
(\vartheta_1 \wedge\vartheta_2 
\wedge 
\xi_2 \wedge\dots \widehat \xi_j\ldots\wedge\xi_n)
\\
&\quad
-
\sum_{2\leq j<k\leq n} 
(-1)^{j+k}
\vartheta_1 \wedge\vartheta_2 \wedge [\xi_j,\xi_k]
\wedge\xi_2\wedge\dots\widehat\xi_j\dots\widehat\xi_k\ldots\wedge\xi_n .
\endalign
$$
We now take $(\vartheta_1,\vartheta_2) = (\xi_1,\xi_2)$. Then
$$
\align
(\do \phi_{\alpha})
&
(\xi_2,\dots, \xi_n)
=
\sum_{2\leq j\leq n} 
(-1)^{j-1}
\nabla_{\xi_j}
(\vartheta_1 \wedge\vartheta_2 
\wedge 
\xi_2 \wedge\dots \widehat \xi_j\ldots\wedge\xi_n)
\\
&\quad
-\sum_{2\leq j<k\leq n} 
(-1)^{j+k}
\vartheta_1 \wedge\vartheta_2 \wedge [\xi_j,\xi_k]
\wedge\xi_2\wedge\dots\widehat\xi_j\dots\widehat\xi_k\ldots\wedge\xi_n
\\
&=
- \nabla_{\xi_2}
(\xi_1 \wedge\xi_2 
\wedge 
\xi_3 \wedge\dots \wedge\xi_n)
\\
&\quad
+\sum_{2<k\leq n} 
(-1)^{k+1}
\xi_1 \wedge\xi_2 \wedge [\xi_2,\xi_k]
\wedge\xi_3\wedge\dots\widehat\xi_k\ldots\wedge\xi_n
\endalign
$$
and
$$
(\do \phi_{\vartheta_1})
(\vartheta_2,\xi_2,\dots, \xi_n)
=0 
$$
whereas, by a calculation the details of which are not given here,
$$
\align
(\do \phi_{\vartheta_2})
(\vartheta_1,\xi_2,\dots, \xi_n)
&=
\nabla_{\xi_2}
(\xi_1 \wedge \xi_2 \wedge \xi_3\wedge\ldots\wedge\xi_n)
\\
&\quad
+[\xi_1,\xi_2] \wedge\xi_2 \wedge \xi_3\wedge\ldots\wedge\xi_n
\\
&\quad
+
\sum_{2<k\leq n} 
(-1)^k
\xi_1\wedge\xi_2 \wedge [\xi_2,\xi_k]
\wedge \xi_3\wedge\dots\widehat\xi_k\ldots\wedge\xi_n .
\endalign
$$
Consequently
$$
\align
 (\do \phi_{\vartheta_1 \wedge \vartheta_2}) (\xi_2,\ldots ,\xi_n)
&-(\do \phi_{\vartheta_1}) (\vartheta_2,\xi_2,\ldots ,\xi_n)
+(\do \phi_{\vartheta_2}) (\vartheta_1,\xi_2,\ldots ,\xi_n)
\\
&=
\sum_{2<k\leq n} 
(-1)^{k+1}
\xi_1 \wedge\xi_2 \wedge [\xi_2,\xi_k]
\wedge\xi_3\wedge\dots\widehat\xi_k\ldots\wedge\xi_n
\\
&\quad
+[\xi_1,\xi_2] \wedge\xi_2 \wedge \xi_3\wedge\ldots\wedge\xi_n
\\
&\quad
+
\sum_{2<k\leq n} 
(-1)^{k}
\xi_1\wedge\xi_2 \wedge [\xi_2,\xi_k]
\wedge \xi_3\wedge\dots\widehat\xi_k\ldots\wedge\xi_n
\\
&=
[\xi_1,\xi_2]\wedge \xi_2  \wedge \xi_3 \ldots \wedge \xi_n
\endalign
$$
that is, with $\beta = \xi_2  \wedge \xi_3\ldots \wedge \xi_n$,
$$
(\Del(\vartheta_1 \wedge \vartheta_2)) \wedge \beta
-((\Del\vartheta_1) \vartheta_2) \wedge \beta
- (\vartheta_1 (\Del\vartheta_2))\wedge \beta
=-[\vartheta_1,\vartheta_2] \wedge \beta
\in \Lambda_A^nQ.
$$
This etablishes the generating property (6.10.3)
for $\alpha_1$ and $\alpha_2$ homogeneous of degree 1
since, as an $A$-module, $Q$ is 
finitely generated and projective of constant rank $n$.
Since, as an $A$-algebra, $\Lambda_AQ$ is generated by its elements of degree 1,
a straightforward induction completes the proof
of Theorem 6.10.  \qed
\enddemo

For the special case where
$\delta$ and hence $\Psi_{\delta}$ is zero, the statement of the theorem
is a consequence of Theorem 5.4.4 in \cite\twilled.

\proclaim{Corollary 6.10.4}
Let $(A,H,Q)$ be an orientable Lie-Rinehart triple,
and let $\Cal G = \roman{Alt}_A(H,\Lambda_AQ)$
be endowed with the Lie-Rinehart differential $d$,
the bigraded bracket {\rm (5.8.3)}, 
and Jacobiator {\rm (5.8.5)}.
Then $(\Cal G,d,[\cdot,\cdot],\Psi)$
is a quasi-Gerstenhaber algebra.
\endproclaim

Indeed, the identity (5.6) then corresponds to (1.9.7);
cf. also (4.9.6) and (6.11)(vii) below.

\smallskip
\noindent
{\smc Remark 6.11.} It is instructive to spell out the relationship between 
the quasi-Batalin-Vilkovisky compatibility conditions (6.2)--(6.4) and the 
Lie-Rinehart triple axioms (1.9.1)--(1.9.7); cf. (2.8.5) above.
As before, write $\Cal G = \roman{Alt}_A(H,\Lambda_A Q)$, and recall that
$n$ is the rank of $Q$ as a projective $A$-module.
\newline\noindent
(i) The vanishing of
$d\Delta + \Delta d \colon \Cal G^0_n \to\Cal G^1_{n-1}$
(a special case of (6.2)) corresponds to (1.9.1).
\newline\noindent
(ii) The vanishing of the operator
$d\Delta + \Delta d \colon \Cal G^1_n \to\Cal G^2_{n-1}$
(a special case of (6.2), too) corresponds to (1.9.2).
\newline\noindent
(iii) The vanishing of
$d\Delta + \Delta d \colon \Cal G^0_{n-1} \to\Cal G^1_{n-2}$
(still a special case of (6.2)) corresponds to (1.9.3).
\newline\noindent
(iv) The vanishing of
$
\Delta\Delta+\Psi d=
d\Psi + \Delta\Delta+\Psi d \colon \Cal G^0_{n} @>>>\Cal G^0_{n-2}
$
(a special case of (6.3)) corresponds to (1.9.4).
\newline\noindent
(v) The vanishing of
$
\Delta\Delta+\Psi d=
d\Psi + \Delta\Delta+\Psi d \colon \Cal G^0_{n-1} @>>>\Cal G^0_{n-3}
$
(a special case of (6.3), too) corresponds to (1.9.5).
\newline\noindent
(vi) The vanishing of
$
d\Psi + \Delta\Delta+\Psi d \colon \Cal G^1_n @>>>\Cal G^1_{n-2}
$
(still a special case of (6.3)) corresponds to (1.9.6).
\newline\noindent
(vii) The vanishing of
$\Delta\Psi + \Psi\Delta\colon \Cal G^1_n \to\Cal G^0_{n-3}$
(a special case of (6.4)) corresponds to (1.9.7).
Cf. also (4.9.6) above.
\smallskip
When $(A,H,Q)$ is an orientable Lie-Rinehart triple, with orientation form 
$\omega$, pursuing the philosophy developed in Section 7 of \cite\twilled\ 
(cf. in particular (7.14)), we may view
$(\roman{Alt}^*_A(H,\Lambda^*_A Q),d,\Delta_{\omega},\Psi_{\delta})$
as an object the category of 
$\Cal A$-modules calculating the \lq\lq quasi-Lie-Rinehart {\it homology\/} 
$\roman H^*_*(\Cal Q,\Cal A_{\Delta})$ of the quasi-Lie-Rinehart algebra
$(\Cal A,\Cal Q)$, with values in the right $(\Cal A,\Cal Q)$-module
$\Cal A_{\Delta}$\rq\rq, the right $(\Cal A,\Cal Q)$-module structure
being induced by $\Delta$. The isomorphism
$$
(\roman{Alt}^*_A(H,\Lambda^*_A Q);d,\Delta_{\omega},\Psi_{\delta})
@>>>
(\roman{Alt}^*_A(H, \roman{Alt}^{n-*}_A(Q,A));d_0,-d_1,d_2)
\tag6.12
$$
is then
a kind of \lq\lq duality isomorphism\rq\rq\  of chain complexes inducing
a \lq\lq duality isomorphism\rq\rq\ 
which, in bidegree $(q,p)$, is
of the kind
$$
\roman H^q_p(\Cal Q,\Cal A_{\Delta})
@>>>
\roman H^{q,n-p}(\Cal Q,\Cal A) \cong \roman H^{q,n-p}(L,A)
\tag6.13
$$
where
$L=H \oplus Q$ is 
the $(R,A)$-Lie algebra
which corresponds to the given Lie-Rinehart triple
$(A,H,Q)$.
Proposition 7.14 in \cite\twilled\ 
makes this precise for the special case
where
$(A,Q,H)$
is a twilled Lie-Rinehart algebra.
In our case, pushing further,
consider the filtrations of
$\roman{Alt}^*_A(H, \Lambda^*_AQ)$
and
$\roman{Alt}^*_A(H, \roman{Alt}^{n-*}_A(Q,\Lambda^n_A Q))$
by $Q$-degree.
In view of what was said above,
the corresponding spectral sequence
(6.7.6), which we now write in the form
$$
(\widehat{\roman E}_*^*(r),d(r)),
\quad
d(r) \colon \widehat{\roman E}_p^q(r)
@>>> \widehat{\roman E}_{p-r}^{q-r+1}(r),
\tag6.14.1
$$
has
$$
(\widehat {\roman E}(0),d(0)) = 
(\roman{Alt}^*_A(H,\Lambda^*_A Q),d)
\tag6.14.2
$$
whence
$$
(\widehat {\roman E}(1),d(1))
= 
(\roman H^*_*(\Cal G)_d,\partial);
\tag6.14.3
$$
this is the bigraded homology Batalin-Vilkovisky algebra
spelled out in Proposition 6.7 above,
for the quasi-Batalin-Vilkovisky algebra
$\Cal G^*_* =(\roman{Alt}^*_A(H,\Lambda^*_A Q);
d,\Delta_{\omega},\Psi_{\delta})$.
The isomorphism (6.8) is compatible with these filtrations.
Hence it identifies the corresponding spectral sequence (2.9.1) 
with (6.14.1). 

\smallskip\noindent
{\smc Illustration 6.15.}
Return to the situation of (1.4.1), and maintain the notation etablished there 
as well as in (2.10), cf. also (4.15). Thus $(M,\Cal F)$ is a 
foliated 
manifold, $(A,H,Q) = (C^{\infty}(M), L_{\Cal F}, Q)$ is the corresponding 
Lie-Rinehart triple,
and $(\Cal A, \Cal Q) = \roman{Alt}_A(H,A),\roman{Alt}_A(H,Q)$
is the corresponding  quasi-Lie-Rinehart algebra. We
now push further the interpretation, advertised already in (4.15) above,
of $\Cal A$ as the {\sl algebra of generalized 
functions\/} and of $\Cal Q$ as the {\sl generalized Lie algebra of vector 
fields for the foliation. This interpretation relies
crucially on the totalization spelled out as\/} (6.7.2) {\sl above;
with the more familiar totalization\/}
$\roman{Tot}' \Cal G$ {\sl given by\/}
$$
(\roman{Tot}' \Cal G)^n = \sum_{p+q = n} \Cal G^q_p,
$$
{\sl such an interpretation is not visible.\/}
\smallskip
Thus, consider the bigraded algebra
$\Cal G^*_*=\roman{Alt}_A^*(H,\Lambda^*_AQ) = \Lambda_{\Cal A} \Cal Q$,
where as before
$\Cal A= \roman{Alt}_A(H,A)$ and  $\Cal Q=\roman{Alt}_A^*(H,Q)$.
Suppose that the foliation is transversely orientable with a 
{\it basic\/} transverse volume form $\omega$, and consider the resulting
quasi-Batalin-Vilkovisky algebra
$(\roman{Alt}^*_A(H,\Lambda^*_A Q);d, \Delta_{\omega},\Psi_{\delta})$,
cf. Theorem 6.10.
In particular,  $\Cal G^*_*$ is then a quasi-Gerstenhaber algebra.
This quasi-Gerstenhaber algebra yields a kind of 
generalized Schouten algebra (algebra of multivector fields) for the 
foliation; the cohomology 
$\roman H_*^0(\Cal G)$ may be viewed as the Schouten algebra
for the \lq\lq space of leaves\rq\rq.
However the entire cohomology contains more information about the
foliation than just $\roman H_*^0(\Cal G)$.
\smallskip
Under the circumstances of (2.10(ii)), where the foliation comes from
a fiber bundle, cf. also (4.15),
let $B$ denote the \lq\lq space of leaves\rq\rq\  or, equivalently,
the base of the corresponding bundle;
an orientation $\omega$ in our sense is now essentially equivalent
to a volume form $\omega_B$ for the base $B$.
Let $L_B=\roman{Vect}(B)$.
The volume form $\omega_B$ induces an exact generator
$\partial_{\omega_B}$
for the ordinary Gerstenhaber algebra
$G_* = \Lambda_{C^{\infty}(B)}L_B$,
and the corresponding bigraded homology Batalin-Vilkovisky algebra
$(\roman H^*_*(\roman{Alt}_A(H,\Lambda^*_A Q))_d,\partial_{\omega})$ 
coming into play
in Theorem 6.10
may then be written
as the bigraded crossed product 
$$
(\roman H^*_*(\roman{Alt}_A(H,\Lambda^*_A Q))_d,\partial_{\omega}) = 
\roman H^*(\Cal A) \otimes_{C^{\infty}(B)} (G_*,\partial_{\omega_B})
\tag6.15.1
$$
of 
$\roman H^*(\Cal A)$
with the ordinary Batalin-Vilkovisky algebra
$(G_*,\partial)= (\Lambda^*_{C^{\infty}(B)} L_B,\partial_{\omega_B})$
(cf. \cite\twilled\ for the notion of bigraded crossed product
Batalin-Vilkovisky algebra);
here $\Cal A = (\roman{Alt}_A(H,A),d)$
which, cf. (2.10(ii)), computes the cohomology of $M$ with values in
the sheaf of germs of functions which are constant on the leaves,
i.~e. fibers.
\smallskip
Under the circumstances of
(2.10(i)), when the foliation does {\it not\/} come from a fiber bundle,
the structure
of
the bigraded homology Batalin-Vilkovisky algebra
$\roman H^*_*(\roman{Alt}_A(H,\Lambda^*_A Q))_d$ 
may be more intricate.

\smallskip\noindent
{\smc Illustration 6.16.}
For a (finite dimensional) quasi-Lie bialgebra $(\fra h,\fra h^*)$ 
\cite\kosmatwe, with Manin pair $(\fra g,\fra h)$, where 
$\fra g =\fra h \oplus \fra h^*$, the resulting quasi-Batalin-Vilkovisky 
algebra has the form
$$
\roman{Alt}(\fra h,\Lambda \fra h^*)
\cong
\Lambda \fra h^* \otimes \Lambda \fra h^*
\cong
\Lambda (\fra h^* \oplus\fra h^*) .
$$

\noindent
{\smc Remark 6.17.}
Given the bigraded commutative algebra $\Cal G^*_*$, consider the totalization
$\roman{Tot}' \Cal G$ spelled out above.
Suppose there be given operators $d,\partial,\Psi$
which endow $\Cal G$ with a quasi-Batalin-Vilkovisky algebra structure
in our sense. These operators induce operators
$$
d\colon
(\roman{Tot}' \Cal G)^* @>>> (\roman{Tot}' \Cal G)^{*+1},
\
\partial\colon
(\roman{Tot}' \Cal G)^* @>>> (\roman{Tot}' \Cal G)^{*-1},
\
\Psi\colon
(\roman{Tot}' \Cal G)^* @>>> (\roman{Tot}' \Cal G)^{*-3}
$$
such that $\Cal L=d \partial + \partial d = 0$,
$d \Psi +\partial\partial + \Psi d= 0$,
$\partial \Psi + \Psi \partial= 0$,
$\Psi \Psi = 0$,
whence, endowed with these operators,
$\roman{Tot}' \Cal G$ is precisely a quasi-Batalin-Vilkovisky
algebra in the sense of \cite\getzlthr\ with zero Laplacian $\Cal L$.
This notion of quasi-Batalin-Vilkovisky algebra
extends that of differential GBV-algebra
in \cite\maninbtw\ (III.9.5)
(which corresponds to the structure under discussion
with $\Psi = 0$, with reference to the totalization $\roman{Tot}' \Cal G$,)
and is a special
case of a more general notion of generalized $\roman{BV}$-algebra
explored in \cite\kravcthr. 
In \cite\bangotwo\ (Definition 3.2), 
a corresponding notion of quasi-Gerstenhaber algebra has been isolated.
When $(\Cal G,d,[\cdot,\cdot],\Psi)$
is a quasi-Gerstenhaber algebra in our sense,
the operations $d,[\cdot,\cdot],\Psi$
induce as well corresponding
pieces of structure $d,[\cdot,\cdot],\Psi$ on $\roman{Tot}' \Cal G$
and, in view of Lemma 2.2 in \cite\bangotwo,
the requirement (5.5) above
(which makes precise how
under our circumstances the $h$-Jacobiator $\Psi$
controls the failure of the strict Jacobi identity)
entails the requirement (3.7) 
in \cite\bangotwo\ 
which, in turn,
describes the failure of the strict Jacobi identity
under the circumstances of \cite\bangotwo. 
Moreover, our requirements (5.i)--(5.iii) in Section 5
above now amount to the corresponding requirements (3.6)--(3.8) in
\cite\bangotwo. Likewise the requirement (5.6) 
corresponds to the requirement
(3.9) in \cite\bangotwo. These observations make precise the relationship
between our notions of quasi-Gerstenhaber 
and of quasi-Batalin-Vilkovisky algebra and that of  
quasi-Gerstenhaber algebra in \cite\bangotwo\ 
and those of quasi-Batalin-Vilkovisky algebra
(with zero Laplacian) explored
in \cite\bangotwo\ and \cite\getzlthr. 
However the notion of Laplacian
does not seem to have a meaning for the totalization
Tot which we use in this paper,
in particular, does not have an interpretation
(at least not an obvious one)
in terms of foliations.

\bigskip 

\widestnumber\key{999}
\centerline{References}
\smallskip\noindent

\ref \no \akmantwo
\by F. Akman
\paper On some generalizations of Batalin-Vilkovisky algebras 
\jour J. of Pure and Applied Algebra
\vol 120
\yr 1997
\pages 105--141
\endref

\ref \no \almolone
\by R. Almeida and P. Molino
\paper Suites d'Atiyah et feuilletages transversalement complets
\jour C. R. Acad. Sci. Paris I 
\vol 300
\yr 1985
\pages 13--15
\endref

\ref \no \atiyathr
\by M. F. Atiyah
\paper Complex analytic connections in fibre bundles
\jour Trans. Amer. Math. Soc.
\vol 85
\yr 1957
\pages 181--207
\endref

\ref \no \bangoone 
\by M. Bangoura 
\paper Quasi-big\`ebres de Lie et alg\`ebres
quasi-Batalin-Vilkovisky 
\linebreak
diff\'erentielles
\jour Comm. in Algebra
\vol 31 (1)
\yr 2003
\pages 29--44
\endref

\ref \no \bangotwo
\by M. Bangoura 
\paper Alg\`ebres quasi-Gerstenhaber diff\'erentielles
\paperinfo preprint
\endref

\ref \no \batviltw
\by I. A. Batalin and G. S. Vilkovisky
\paper Quantization of gauge theories
with linearly dependent generators
\jour  Phys. Rev. 
\vol D 28
\yr 1983
\pages  2567--2582
\endref

\ref \no \batvilfo
\by I. A. Batalin and G. S. Vilkovisky
\paper Closure of the gauge algebra, generalized Lie equations
and Feynman rules
\jour  Nucl. Phys. B
\vol 234
\yr 1984
\pages  106-124
\endref

\ref \no \batavilk
\by I. A. Batalin and G. S. Vilkovisky
\paper Existence theorem for gauge algebra
\jour Jour. Math. Phys.
\vol 26
\yr 1985
\pages  172--184
\endref

\ref \no \ecartone
\by \'E. Cartan
\paper La g\'eom\'etrie des groupes de transformations
\jour 
Journal de Math\'e-
\linebreak
matiques pures et appliqu\'ees
\vol  6 
\yr 1927
\pages 1--119
\endref

\ref \no \ecarttwo
\by \'E. Cartan
\paper Sur les invariants int\'egraux de certains espaces homog\`enes clos
et les propri\'et\'es topologiques de ces espaces
\jour Ann. Soc. Pol. Math.
\vol  8 
\yr 1929
\pages 181--225
\endref

\ref \no \drinftwo
\by V. G. Drinfeld
\paper Quasi-Hopf algebras
\jour Leningrad Math. J.
\vol 1
\yr 1990
\pages 1419--1457
\endref

\ref \no \fronione
\by A. Fr\"olicher and A. Nijenhuis
\paper 
Theory of vector-valued differential forms, Part I.
Derivations in the graded ring of differential forms
\jour Indagationes Math.
\vol 18
\yr 1956 
\pages 338--359
\endref

\ref \no \gersthtw
\by M. Gerstenhaber
\paper The cohomology structure of an associative ring
\jour Ann. of Math.
\vol 78
\yr 1963
\pages  267-288
\endref

\ref \no \getzlthr
\by E. Getzler
\paper Manin pairs and topological conformal field theory
\jour Annals of Phys.
\vol 237
\yr 1995
\pages 161--201
\endref

\ref \no \helleron
\by A. Heller
\paper Homological resolutions of complexes with operators
\jour  Ann. of Math.
\vol 60
\yr 1954
\pages  283--303
\endref

\ref \no \poiscoho
\by J. Huebschmann
\paper Poisson cohomology and quantization
\jour 
J. f\"ur die reine und angewandte Mathematik
\vol  408 
\yr 1990
\pages 57--113
\endref

\ref \no  \souriau
\by J. Huebschmann
\paper On the quantization of Poisson algebras
\paperinfo Symplectic Geometry and Mathematical Physics,
Actes du colloque en l'honneur de Jean-Marie Souriau,
P. Donato, C. Duval, J. Elhadad, G.M. Tuynman, eds.;
Progress in Mathematics, Vol. 99
\publ Birkh\"auser Verlag
\publaddr Boston $\cdot$ Basel $\cdot$ Berlin
\yr 1991
\pages 204--233
\endref

\ref \no \extensta
\by J. Huebschmann
\paper 
Extensions of Lie-Rinehart algebras and the Chern-Weil construction
\paperinfo in: Festschrift in honor of J. Stasheff's 60th birthday
\jour Cont. Math. 
\vol 227
\yr 1999
\pages 145--176
\publ Amer. Math. Soc.
\publaddr Providence R. I.
\endref

\ref \no \bv
\by J. Huebschmann
\paper Lie-Rinehart algebras, Gerstenhaber algebras, and Batalin-
Vilkovisky algebras
\jour Annales de l'Institut Fourier
\vol 48
\yr 1998
\pages 425--440
\endref

\ref \no \duality
\by J. Huebschmann
\paper 
Duality for Lie-Rinehart algebras and the modular class
\jour Journal f\"ur die Reine und Angew. Math. 
\vol 510
\yr 1999
\pages 103--159
\endref

\ref \no \twilled
\by J. Huebschmann
\paper Twilled Lie-Rinehart algebras and differential Batalin-Vilkovisky 
algebras
\paperinfo {\tt math.DG/9811069}
\endref

\ref \no \banach
\by J. Huebschmann
\paper Differential Batalin-Vilkovisky algebras arising from
twilled Lie-Rinehart algebras
\paperinfo Poisson Geometry 
\jour Banach Center publications 
\vol 51
\yr 2000
\pages 87--102
\endref

\ref \no \lradq
\by J. Huebschmann
\paper Lie-Rinehart algebras, descent, and quantization
\jour Fields Institute Communications (to appear)
\finalinfo {\tt math.SG/0303016}
\endref

\ref \no \berikas
\by J. Huebschmann
\paper Berikashvili's functor $\Cal D$ and the deformation equation
\paperinfo Festschrift in honor of N. Berikashvili's 70th birthday;
{\tt math.AT/9906032}
\jour Proceedings of the A. Razmadze Mathematical Institute
\vol 119
\yr 1999
\pages 59--72
\endref

\ref \no \kaehler
\by J. Huebschmann
\paper K\"ahler spaces, nilpotent orbits, and singular reduction
\linebreak
\paperinfo {\tt math.dg/0104213}
\jour Memoirs of the AMS (to appear)
\endref

\ref \no \qr
\by J. Huebschmann
\paper K\"ahler quantization and reduction
\jour  J. f\"ur die reine und angewandte Mathematik (to appear)
\finalinfo {\tt math.SG/0207166}
\endref

\ref \no \perturba
\by J. Huebschmann
\paper Perturbation theory and free resolutions for nilpotent
groups of class 2
\jour J. of Algebra
\yr 1989
\vol 126
\pages 348--399
\endref

\ref \no \cohomolo
\by J. Huebschmann
\paper Cohomology of nilpotent groups of class 2
\jour J. of Algebra
\yr 1989
\vol 126
\pages 400--450
\endref

\ref \no \modpcoho
\by J. Huebschmann
\paper The mod $p$ cohomology rings of metacyclic groups
\jour J. of Pure and Applied Algebra
\vol 60
\yr 1989
\pages 53--105
\endref

\ref \no \intecoho
\by J. Huebschmann
\paper Cohomology of metacyclic groups
\jour Trans. Amer. Math. Soc.
\vol 328
\yr 1991
\pages 1-72
\endref

\ref \no \huebkade
\by J. Huebschmann and T. Kadeishvili
\paper Small models for chain algebras
\jour Math. Z.
\vol 207
\yr 1991
\pages 245--280
\endref

\ref \no \huebstas
\by J. Huebschmann and J. D. Stasheff
\paper Formal solution of the master equation via HPT and
deformation theory
\paperinfo {\tt math.AG/9906036}
\jour Forum mathematicum 14
\yr 2002
\pages 847--868
\endref

\ref \no \bkellone
\by B. Keller
\paper Introduction to $A_{\infty}$-algebras and modules
\jour 
Homology, Homotopy and its Applications
\vol  3 
\yr 2001
\pages 1--35
\moreref
Addendum
\jour 
Homology, Homotopy and its Applications
\vol  4 
\yr 2002
\pages 25--28
\endref

\ref \no \kikkaone
\by M. Kikkawa
\paper Geometry of homogeneous Lie loops
\jour Hiroshima Math. J.
\vol  5 
\yr 1975
\pages 141--179
\endref

\ref \no \kinywein
\by M. K. Kinyon and A. Weinstein
\paper Leibniz algebras, Courant algebroids, and multiplications on reductive
homogeneous spaces
\jour Amer. J. of Math.
\vol 123
\yr 2001
\pages 525--550
\finalinfo
\endref

\ref \no \kjesetwo
\by L. J. Kjeseth
\paper Homotopy Rinehart cohomology of homotopy Lie-Rinehart pairs
\jour Homology, Homotopy and its Applications
\vol 3
\yr 2001
\pages 139--163
\endref

\ref \no \kosmatwe
\by Y. Kosmann-Schwarzbach 
\paper Jacobian quasi-bialgebras and quasi-Poisson Lie groups
\jour  Cont. Math.
\vol 132
\yr 1992
\pages 459--489
\endref

\ref \no \kosmathr
\by Y. Kosmann-Schwarzbach 
\paper Exact Gerstenhaber algebras and Lie bialgebroids
\jour  Acta Applicandae Mathematicae
\vol 41
\yr 1995
\pages 153--165
\endref

\ref \no \kosmafou
\by Y. Kosmann-Schwarzbach 
\paper From Poisson algebras to Gerstenhaber algebras
\jour Annales de l'Institut Fourier
\vol 46
\yr 1996
\pages 1243--1274
\endref

\ref \no \kosmafte
\by Y. Kosmann-Schwarzbach 
\paper Quasi, twisted, and all that ..., in Poisson geometry
and Lie algebroid theory
\paperinfo Preprint, 2003
\endref

\ref \no \koszulon
\by J. L. Koszul
\paper Crochet de Schouten-Nijenhuis et cohomologie
\jour Ast\'erisque,
\vol hors-s\'erie,
\yr 1985
\pages 251--271
\paperinfo in E. Cartan et les Math\'ematiques d'aujourd'hui, 
Lyon, 25--29 Juin, 1984
\endref

\ref \no \kravcthr
\by O. Kravchenko
\paper Deformations of Batalin-Vilkovisky algebras
\paperinfo Poisson Geometry 
\jour Banach Center publications
\vol 51
\yr 2000
\pages 131--139
\endref

\ref \no \liuleone
\by A. Liulevicius
\paper Multicomplexes and a general change of rings theorem
\linebreak
\paperinfo mimeographed notes, University of Chicago
\endref

\ref \no \liuletwo
\by A. Liulevicius
\paper A theorem in homological algebra and stable homotopy of projective
spaces
\jour Trans. Amer. Math. Soc.
\vol 109
\yr 1963
\pages  540--552
\endref

\ref \no \mackone
\by K. Mackenzie
\book Lie groupoids and Lie algebroids in differential geometry
\bookinfo London Math. Soc. Lecture Note Series, vol. 124
\publ Cambridge University Press
\publaddr Cambridge, England
\yr 1987
\endref

\ref \no \maninbtw
\by Yu. I. Manin
\book Frobenius manifolds, quantum cohomology, and moduli spaces
\bookinfo Colloquium Publications, vol. 47
\publ Amer. Math. Soc.
\publaddr Providence, Rhode Island
\yr 1999
\endref

\ref \no \molinobo
\by P. Molino
\book Riemannian foliations
\bookinfo Progress in Mathematics, No. 73
\publ Birkh\"auser Verlag
\publaddr Boston $\cdot$ Basel $\cdot$ Berlin
\yr 1988
\endref

\ref \no \nijenhui
\by A. Nijenhuis
\paper Jacobi-type identities for bilinear differential concomitants
of tensor fields.~I
\jour Indag. Math.
\vol 17
\yr 1955
\pages 390--403
\endref

\ref \no \nomiztwo
\by K. Nomizu
\paper Invariant affine connections on homogeneous spaces
\jour Amer. J. of Math.
\vol  76 
\yr 1954
\pages 33--65
\endref

\ref \no \palaione
\by R. S. Palais
\paper The cohomology of Lie rings
\jour  Proc. Symp. Pure Math.
\vol III
\yr 1961
\pages 130--137
\paperinfo Amer. Math. Soc., Providence, R. I.
\endref

\ref \no \breinhar
\by B. L. Reinhart
\paper Foliated manifolds with bundle-like metrics
\jour Ann. of Math.
\vol  69 
\yr 1959
\pages 119--132
\endref

\ref \no \rinehone
\by G. Rinehart
\paper Differential forms for general commutative algebras
\jour  Trans. Amer. Math. Soc.
\vol 108
\yr 1963
\pages 195--222
\endref

\ref \no \roytbtwo
\by D. Roytenberg 
\paper Quasi-Lie bialgebroids and twisted Poisson manifolds
\jour Lett. in Math. Phys. 
\vol 61
\yr 2002
\pages 123--137
\endref

\ref \no \roytwein
\by D. Roytenberg and A. Weinstein
\paper Courant algebroids and strongly homotopy Lie algebras
\jour Lett. in Math. Phys. 
\vol 46
\yr 1998
\pages 81--93
\endref

\ref \no \sabimikh
\by L. V. Sabinin and P. O. Mikheev
\paper On the infinitesimal theory of local analytic loops
\jour Soviet Math. Dokl.
\vol  36 
\yr 1988
\pages 545--548
\endref

\ref \no \sarkhone
\by K. S. Sarkharia
\paper The de Rham cohomology of foliated manifolds
\paperinfo Ph. D. thesis, Stony Brook, 1974
\endref

\ref \no \sarkhtwo
\by K. S. Sarkharia
\paper Non-degenerescence of some spectral sequences
\jour Annales de l'Institut Fourier
\vol 34
\yr 1984
\pages 39--46
\endref

\ref \no \stashnin
\by J. D. Stasheff
\paper Deformation theory and the Batalin-Vilkovisky 
master equation
\paperinfo in: Deformation Theory and Symplectic Geometry,
Proceedings of the Ascona meeting, June 1996,
D. Sternheimer, J. Rawnsley, S.  Gutt, eds.,
Mathematical Physics Studies, Vol. 20
\publ Kluwer Academic Publishers
\publaddr Dordrecht/Boston/London
\yr 1997
\pages 271--284
\endref

\ref \no \stashset
\by J. D. Stasheff
\paper Poisson homotopy algebra. An idiosyncratic survey of
homotopy algebraic topics related to Alan's interests
\paperinfo Preprint, 2003
\endref

\ref \no \vanesthr
\by W. T. Van Est
\paper Alg\`ebres de Maurer-Cartan et holonomie
\jour Ann. Fac. Sci. Toulouse Math.
\vol (5) (suppl.)
\yr 1989
\pages 93--134
\endref

\ref \no \weinomni
\by A. Weinstein
\paper Omni-Lie algebras
\paperinfo Microlocal analysis of the Schr\"odinger equation
and related topics, Kyoto, 1999
\jour S\hataccent urikaisekikenky\hataccent usho 
K\hataccent oky\hataccent uroku
(RIMS)
\vol 1176
\yr 2000
\pages 95--102
\endref

\ref \no \yamaguti
\by K. Yamaguti
\paper On the Lie triple systems and its generalizations
\jour J. Sci. Hiroshima Univ. Ser. A
\vol  21 
\yr 1957/58
\pages 155-160
\endref

\enddocument